\documentclass[a4paper,11pt]{amsart}

\usepackage{amssymb}
\usepackage{xy}
\usepackage{bm}
\xyoption{all}
\xyoption{poly}
\xyoption{arc}

\makeindex

\newtheorem{theo}{Theorem}[section]
\newtheorem{coro}[theo]{Corollary}
\newtheorem{prop}[theo]{Proposition}
\newtheorem{lemma}[theo]{Lemma}

\newtheorem{conj}[theo]{Conjecture}
\theoremstyle{definition}
\newtheorem{defi}[theo]{Definition}
\newtheorem{exercise}[theo]{Exercise}

\newtheorem{remark}[theo]{Remark}
\newtheorem{question}[theo]{Question}
\newtheorem{example}[theo]{Example}
\newtheorem{convention}[theo]{Convention}

\def\cq{{\backslash}} 

\def\CC{{\mathcal{C}}}
\def\CS{{\mathcal{S}}}

\def\CO{{\mathcal{O}}}

\def\CL{{\mathcal{L}}}
\def\CN{{\mathcal{N}}}

\def\CH{{\mathcal{H}}}

\def\CV{{\mathcal{V}}}
\def\CU{{\mathcal{U}}}
\def\ACA{{\mathcal{A}}}

\def\CG{{\mathcal{G}}}

\def\CT{{\mathcal{T}}}
\def\CO{{\mathcal{O}}}
\def\CL{{\mathcal{L}}}

\def\CS{{\mathcal{S}}}
\def\CK{{\mathcal{K}}}
\def\CD{{\operatorname{D}}}
\def\CDbullet{{\operatorname{D}_{\bullet}}}

\def\CU{{\mathcal{U}}}
\def\CVW{{\widehat{\CU}_1}}
\def\CVWr{\CV_g}

\def\BC{{\mathbf{C}}}

\def\BC{{\mathbb{C}}}
\def\BR{{\mathbb{R}}}
\def\BZ{{\mathbb{Z}}}

\def\BM{{\mathbf{M}}}

\def\Braid{B}

\def\wt{\operatorname{wt}\nolimits}

\def\Id{\operatorname{Id}\nolimits}

\def\GL{\operatorname{GL}\nolimits}

\def\UC{\operatorname{UniCover}\nolimits}

\def\cc{\operatorname{cc}\nolimits}
\def\lbl{\operatorname{lbl}\nolimits}
\def\rlbl{\operatorname{rlbl}\nolimits}
\def\clbl{\operatorname{clbl}\nolimits}
\def\LL{LL}
\def\LLL{\overline{LL}}
\def\En{E_n}

\def\Red{\operatorname{Red}\nolimits}
\def\Aut{\operatorname{Aut}\nolimits}

\def\Cat{\operatorname{Cat}\nolimits}

\def\Hom{\operatorname{Hom}\nolimits}
\def\bar{\operatorname{bar}\nolimits}
\def\gar{\operatorname{gar}\nolimits}

\def\Ext{\operatorname{Ext}\nolimits}
\def\Diag{\operatorname{Diag}\nolimits}

\def\reg{\operatorname{reg}\nolimits}
\def\gen{\operatorname{gen}\nolimits}
\def\obj{\operatorname{obj}\nolimits}

\def\im{\operatorname{im}\nolimits}
\def\re{\operatorname{re}\nolimits}

\def\Poin{\operatorname{Poin}\nolimits}

\def\Spec{\operatorname{Spec}\nolimits}

\def\Disc{\operatorname{Disc}\nolimits}

\def\ie{{\em i.e.}}

\title{Finite complex reflection arrangements are $K(\pi,1)$}
\author{David Bessis}
\address{DMA, \'Ecole normale sup\'erieure, 
45 rue d'Ulm, 75230 Paris cedex 05, France}
\email{david.bessis@ens.fr}

\address{tinyclues, 15 rue du Caire, 75002 Paris, France}
\email{db@tinyclues.com}

\begin{document}

\setcounter{tocdepth}{1}

\begin{abstract}
Let $V$ be a finite dimensional complex vector space and
$W\subseteq \GL(V)$ be a finite complex reflection group.
Let $V^{\reg}$ be the complement in $V$ of the reflecting
hyperplanes.
We prove that $V^{\reg}$ is a $K(\pi,1)$ space. This
was predicted by a classical conjecture, originally stated by Brieskorn
for complexified real reflection groups.
The complexified real case
follows from a theorem of Deligne and, after contributions
by Nakamura and Orlik-Solomon, only six exceptional
cases remained open. In addition to solving this six cases,
our approach is applicable to most previously known cases,
including complexified real groups for which we obtain a new
proof, based on new geometric objects.
We also address a number of questions about $\pi_1(W\cq V^{\reg})$,
the braid group of $W$.
This includes a description of periodic elements in terms
of a braid analog of Springer's theory of regular elements.
\end{abstract}

\maketitle

\tableofcontents

\printindex

\section*{Notation}

We save the letter $i$ for indexing purposes and denote
by $\sqrt{-1}$ a complex square root of $-1$ fixed once for all.
If $n$ is a positive integer, we denote by $\zeta_n$ the standard $n$-th
root of unity $\exp(2\sqrt{-1}\pi /n)$.

Many objects depend on a complex reflection group $W$, e.g., the braid
group $B(W)$. We often drop the explicit mention of $W$, and write
$B$ for $B(W)$. When $n$ is an integer, we denote by $B_n$ the braid
group on $n$ strings, together with its standard generating set
$\bm{\sigma}_1,\dots,\bm{\sigma}_{n-1}$; it is isomorphic to the braid group
of $\mathfrak{S}_n$ in its permutation reflection representation
(see Section \ref{sectionconfig}). The groups $B$ and $B_n$ appear
simultaneously and should not be confused.
 
\section*{Introduction}

Let $V$ be a finite dimensional complex vector space and
$W\subseteq \GL(V)$ be a complex reflection group
(all reflection groups considered here are assumed to be finite).

Let $V^{\reg}$ be the complement in $V$ of the reflecting
hyperplanes.
In the case when $W$ is a type $A$ reflection group,
Fadell and Neuwirth proved in the early 1960s that
$V^{\reg}$ is a $K(\pi,1)$ space (this is an elementary use of fibration
exact sequences, see \cite{fadell}).
Brieskorn conjectured in 1971, \cite{brieskorn2},
that the $K(\pi,1)$ property holds
when $W$ is a complexified real reflection group.
It is not clear who first stated the conjecture in the context
of arbitrary complex reflection groups. It may be found, for example,
in Orlik-Terao's book:

\begin{conj}[\cite{orlikterao}, p. 163 \& p. 259]
The universal cover of $V^{\reg}$ is contractible.
\end{conj}

Our main result is a proof of this conjecture.
It is clearly sufficient to consider the case when $W$ is irreducible,
which we assume from now on.
Irreducible complex reflection groups have been classified by Shephard-Todd,
\cite{shephard}.

The complexified real case (\ie, Brieskorn's conjecture)
was quickly settled by Deligne, \cite{deligne}.
The rank $2$ complex case is trivial.
The case of the infinite family $G(de,e,n)$ 
was solved in 1983 by Nakamura, \cite{nakamura}
(here again, the monomiality of the group allows an efficient use of
fibrations).
A few other cases immediately follow from the observation by Orlik-Solomon,
\cite{orliksolomon}, that certain discriminants of non-real
 complex reflection groups
are isomorphic to discriminants of complexified real reflection groups.

Combining all previously known results, the conjecture remained open for
six exceptional types: $G_{24}$, $G_{27}$, $G_{29}$, $G_{31}$,
$G_{33}$ and $G_{34}$. We complete the proof of the conjecture by dealing
with these cases.

Let $d_1 \leq \dots \leq d_n$ be the degrees of $W$.
Let $d_1^* \geq \dots \geq d_n^*=0$ be the codegrees of $W$.
We say that $W$ is a \emph{duality} group if
$d_i+d_i^*=d_n$ for all $i$ (by analogy with the real case, we then say
that $d_n$ is the
\emph{Coxeter number} of $W$, denoted by $h$).
We say that $W$ is \emph{well-generated} if it may be generated 
by $n$ reflections.
\index{well-generated reflection groups}
Orlik-Solomon observed, by inspecting the classification of Shephard-Todd,
that
$$W \text{ is a duality group } \Leftrightarrow W \text{ is well-generated}.$$

The first ten sections of
 this article are devoted to the proof of the following theorem:
\begin{theo}
\label{theokapiun}
Let $W$ be a well-generated complex reflection group. The universal
cover of $V^{\reg}$ is contractible.
\end{theo}

The proof relies on combinatorial and geometric objects
which are specific to well-generated groups. It is essentially ``case-free'',
although a few combinatorial lemmas still require some limited use of the Shephard-Todd classification.

Five of the six open cases are well-generated:
$G_{24}$, $G_{27}$, $G_{29}$, $G_{33}$ and $G_{34}$.
The theorem also applies to the complexified 
real case, for which we obtain a new proof,
not relying on \cite{deligne}.

The remaining case, $G_{31}$, is not well-generated: it is
an irreducible complex reflection group
of rank $4$ which cannot be generated by less than $5$ reflections.
Fortunately, we may view it as the centralizer of a $4$-regular
element (in the sense of Springer, \cite{springer}) in the group $G_{37}$
(the complexification of the real group of type $E_8$).
By refining the geometric and combinatorial tools introduced in the study
of the duality case, one obtains a relative version of
Theorem \ref{theokapiun}:

\begin{theo}
\label{theo31}
Let $W$ be a well-generated 
complex reflection group. Let $d$ be a Springer regular number,
let $\zeta$ be a primitive complex $d$-th root of unity,
let $w\in W$ be a $\zeta$-regular element.
Let $V':=\ker(w-\zeta)$, let $W'$ be the centralizer of $w$ in $W$,
viewed as a complex reflection group acting on $V'$ (see \cite{springer}).
Let $V'^{\reg}$ be the associated hyperplane complement.
The universal
cover of $V'^{\reg}$ is contractible.
\end{theo}

This in particular applies to $G_{31}$ and, based on earlier results,
completes the proof of the $K(\pi,1)$ conjecture.

Note that the case $d=1$ in Theorem \ref{theo31} is precisely Theorem
\ref{theokapiun}. We state the two results separately because it better
reflects the organisation of the paper.

As by-products of our construction, we obtain new cases of
several standard conjectures about the braid group of $W$, defined
by  $$B(W) := \pi_1(W \cq V^{\reg}).$$

\smallskip
{\flushleft \bf Proving Theorem \ref{theokapiun}: general strategy.}
The general architecture of our proof is borrowed from
Deligne's original approach but the details are quite different.
Every construction here is an analogue
of a construction from \cite{deligne} but relies on different combinatorial
and geometric objects.
Like in \cite{deligne}, one studies a certain braid monoid $M$, whose
structure expresses properties of reduced decompositions in $W$,
and one proves
that it is a lattice for the divisibility order
(this amounts to saying that the monoid is \emph{Garside}). Like in
\cite{deligne}, one uses semi-algebraic geometry to construct
an open covering of the universal cover of $V^{\reg}$, with the property
that non-empty intersections are contractible. This implies that
the universal cover is homotopy equivalent to the nerve of the covering.
Like in \cite{deligne}, one interprets this nerve as a certain flag complex
obtained from $M$. Like in \cite{deligne}, the contractibility of the nerve
follows from the lattice property for $M$.
However, our proof does not use the classical braid monoid, but
a dual braid monoid (\cite{dualmonoid}, \cite{BC}), whose construction
is generalized to all well-generated complex reflection groups.
The construction of the open covering is the most problematic step:
by contrast with the real case, one cannot rely on the notions of \emph{walls}
and \emph{chambers}. The idea here is to work in $W\cq V^{\reg}$ and to use
a generalization of the Lyashko-Looijenga morphism. This morphism
allows a description of $W\cq V^{\reg}$ by means of a
 ramified covering of a type
$A$ reflection orbifold.
Classical objects like \emph{walls}, \emph{chambers} and \emph{galleries}
can somehow be ``pulled-back'',
via the Lyashko-Looijenga morphism, to give semi-algebraic objects
related to the dual braid monoid.

\smallskip
{\flushleft \bf Lyashko-Looijenga coverings.}
The
quotient map $\pi:V^{\reg} \twoheadrightarrow W \cq V^{\reg}$ is a regular
covering.
Once a system of basic invariants $(f_1,\dots,f_n)$ is chosen, the quotient
space
$W \cq V^{\reg}$ identifies with the complement in $\BC^n$ 
of an algebraic hypersurface
$\CH$, the \emph{discriminant},
\index{$\CH$ (discriminant locus)}
\index{discriminant}
of equation $\Delta\in\BC[X_1,\dots,X_n]$.
If $W$ is an irreducible duality complex reflection group,
it is possible to choose $(f_1,\dots,f_n)$ such that 
$$\Delta= X_n^n + \alpha_2 X_n^{n-2} + \dots 
+ \alpha_n,$$
where $\alpha_2,\dots,\alpha_n\in \BC[X_1,\dots,X_{n-1}]$.
Let $Y := \Spec \BC[X_1,\dots,X_{n-1}]$, together with the natural
map $p:W \cq V^{\reg} \rightarrow Y$.
We have an identification  $W\cq V \simeq \BC^n \stackrel{\sim}{\rightarrow}
 Y\times\BC$ sending the orbit $\overline{v}$ of $v\in V$ to
$(p(\overline{v}),f_n(v))$.
The fiber of $p$ over $y\in Y$ is a line $L_y$ which
intersects $\CH$ at $n$ points (counted with multiplicities).
Generically,
the $n$ points are distinct. Let $\CK$ be the \emph{bifurcation locus},
\ie, the algebraic hypersurface of $Y$ consisting of points $y$ such
that the intersection has cardinality $<n$.
Classical results from invariant theory of complex reflection groups
make possible (and very easy) to generalize a
 construction by Looijenga and Lyashko: the map $\LL$ (for
``Lyashko-Looijenga'') sending
$y\in Y-\CK$ to the subset $\{x_1,\dots,x_n\} \subseteq
\BC$ such that $p^{-1}(y)\cap \CH  = \{(y,x_1),\dots,(y,x_n)\}$ is a regular
covering of degree $n!h^n/|W|$ of the (centered)
configuration space of $n$ points in $\BC$. In particular,
$Y-\CK$ is a $K(\pi,1)$.
This observation, which is apparently new in the non-real case, already
allows a refinement of our earlier results (\cite{zariski}, \cite{BM})
on presentations for the braid group of $W$.

\smallskip
{\flushleft \bf The dual braid monoid.}
When $W$ is complexified real, a \emph{dual braid monoid} was constructed in
\cite{dualmonoid} (generalising the construction of Birman-Ko-Lee, \cite{bkl};
similar partial results were
independently obtained by Brady-Watt, \cite{bw}).
The construction was later generalized in \cite{BC}
to the complex reflection group $G(e,e,n)$.
Let $R$ be the set of all reflections in a well-generated group
$W$. The idea is that the pair $(W,R)$ has some ``Coxeter-like'' features.
Instead of looking at relations of the type
 $$\underbrace{sts...}_{m_{s,t}}=\underbrace{tst...}_{m_{s,t}},$$ one
considers relations of the type $$st=tu$$ where $s,t,u\in R$.
Let $\mathcal{S}$ be the set of all relations of this type holding in $W$.
In general, $B(W) \not \simeq \left< R \left|  \mathcal{S} \right.
\right>$, but it is possible to find natural subsets $R_c \subseteq R$ and
$\mathcal{S}_c \subseteq
\mathcal{S}$ such that $B(W) 
\simeq \left< R_c \left|  \mathcal{S}_c \right.
\right>$
(if $W$ is complexified real, $R_c=R$).
The elements of $\mathcal{S}_c$ are called \emph{dual braid relations}.
The choices of $R_c$ and $\mathcal{S}_c$ are  
natural once a \emph{Coxeter element} $c$ has been chosen
(the notion of \emph{Coxeter element} generalizes to the non-real
well-generated groups in terms of Springer's theory of regular elements).
Since the relations are positive, one may view the presentation as
a monoid presentation, defining a monoid $M(W)$.
The crucial property of this monoid is that it is a lattice for the
divisibility order or, more precisely, a \emph{Garside monoid}.
Following Deligne, Bestvina, T. Brady and Charney-Meier-Whittlesey,
the Garside structure provides a convenient simplicial Eilenberg-McLane
$K(B(W),1)$ space (\cite{deligne}, \cite{bestvina}, \cite{brady}, \cite{cmw}).
The earlier results on the dual braid monoid are improved here in two
directions:
\begin{itemize}
\item The construction is generalized to the few exceptional cases 
($G_{24}$, $G_{27}$, $G_{29}$,
$G_{33}$ and $G_{34}$) not
covered by \cite{dualmonoid} and \cite{BC}.
\item A new geometric interpretation is given, via the Lyashko-Looijenga
covering. This interpretation is different from the one given in
\cite[Section 4]{dualmonoid}.
\end{itemize}

The second improvement is the most important. It relies (so far) 
on a counting argument,
following and extending a property which, for the complexified real case,
was conjectured by Looijenga and proved
in a letter from Deligne to Looijenga, \cite{deligneletter}.

\smallskip
{\flushleft \bf Tunnels.}
The classical theory of real reflection groups combines a
``combinatorial'' theory (Coxeter systems) and a ``geometric''
theory (expressed in the language, invented by Tits, of
walls, chambers, galleries, buildings...)
We expect the dual braid monoid approach to eventually provide
effective substitutes for much of this classical theory.
A first step in this direction is the
notion of \emph{tunnel}, which is a rudimentary
geometric object replacing
the classical notion of \emph{minimal gallery between two chambers}.
An important difference with the classical geometric language is
that tunnels are naturally visualized in $W\cq V$ (instead
of $V$).
A tunnel $T$ 
is a path in $W\cq V^{\reg}$ drawn inside a single
line $L_y$ (for some $y\in Y$) and with
constant imaginary part. It represents an element $b_T$ of the dual braid
monoid $M$. An element of $M$ is \emph{simple} if it is represented by
a tunnel. This notion coincides with the notion of simple element 
associated with the Garside structure.
In the classical approach, for any chamber $\mathcal{C}$, there are as
many equivalence classes of 
minimal galleries starting at $\mathcal{C}$ as simple elements (this number
is $|W|$). Here the situation is different: in a given $L_y$, not all
simple elements are represented. The simple elements represented in different $L_y$'s
may be compared thanks to a huge ``fat basepoint'' $\CU$ which
is both dense in $W\cq V^{\reg}$ and contractible.

\smallskip
{\bf \flushleft Proving Theorem \ref{theo31}.}
The strategy is the same as for Theorem \ref{theokapiun}.
With the notations of the theorem, the quotient space $W'\cq V'^{\reg}$
may be identified with $(W\cq V^{\reg})^{\mu_d}$, for the natural
action of the cyclic group $\mu_d$. 
This action induces an automorphism of $B(W)$ which, unfortunately,
does not preserve the dual braid monoid.
However, it is possible to replace $B(W)$ by a sort of categorical
barycentric subdivision, its \emph{$d$-divided Garside category
$M_d$}, on which $\mu_d$ acts by \emph{diagram automorphisms}.
This construction is explained in my separate article \cite{cyclic}
and recalled in Appendix \ref{garsidependix}.
The fixed subcategory
$M_d^{\mu_d}$ is again a Garside category. It should be thought of
a \emph{dual braid category for $B(W')$} and gives rise
to a natural simplicial space, whose realisation is an
Eilenberg-MacLane space.
As before, one shows that
$(W\cq V^{\reg})^{\mu_d}$ is homotopy equivalent to this simplicial model,
by studying the nerve of a certain open covering of a certain model 
of the universal
cover of $(W\cq V^{\reg})^{\mu_d}$, very similar to the one used
for $W\cq V^{\reg}$ (except that one has to replace the contractible
``basepoint'' $\CU$ by a family of non-overlapping contractible ``basepoints'',
one for each object of $M_d^{\mu_d}$.) This involves 
replacing tunnels by a suitable notion of \emph{circular tunnels}.
Section \ref{section11} focuses the geometric aspects of the proof of
Theorem \ref{theo31} -- it is probably fair to say that the
true explanation lies in the properties of $M_d$ and in the
general theorems about periodic elements in Garside groupoids that
are explained in \cite{cyclic}.
\smallskip

{\bf \flushleft By-products.}

\begin{theo}
\label{theointrogarside}
Braid groups of well-generated complex reflection groups
are Garside groups.
\end{theo}

In the situation of Theorem \ref{theo31}, we prove that the braid group
$B(W')$ is a weak Garside group, which is almost as good.

In particular, $B(W)$ is torsion-free, admits nice solutions to the word
and conjugacy problems, is biautomatic, admits a finite $K(\pi,1)$
(our construction provides an explicit one), and much more -- see
\cite{dehgar} for a quite complete reference. None of this was known
for the six exceptional groups mentioned above.

\begin{theo}[Theorem \ref{theoZZ}]
\label{theocenter}
The center of the braid group of an irreducible
complex reflection group is cyclic.
\end{theo}

Again, the cases of
$G_{24}$, $G_{27}$, $G_{29}$, $G_{31}$, $G_{33}$ and $G_{34}$ are new.
This settles a conjecture by Brou\'e-Malle-Rouquier, \cite{bmr}.

For $B(G_{29})$, $B(G_{31})$, $B(G_{33})$ and $B(G_{34})$, no presentations
were known until now, although some conjectures made in \cite{BM} were
supported by strong evidences.

\begin{theo}
\label{theointropres}
The conjectural presentations
for $B(G_{29})$, $B(G_{31})$, $B(G_{33})$ and $B(G_{34})$ given in \cite{BM}
are correct.
\end{theo}

Combined with \cite{BDM} and \cite{BM}, this completes the longstanding
task of finding presentations for all generalized braid groups associated
with finite complex reflection groups.
Theorem \ref{theointropres} is much easier than the
previously mentioned results and only relies on a minor improvement over 
\cite{zariski} and \cite{BM}. However,
the material presented here allows for a more conceptual proof.

\smallskip

{\bf \flushleft Periodic elements in braid groups.}
In connection with their work on Deligne-Lusztig varieties
(see \cite{boston} for more details), Brou\'e-Michel predicted the
existence of an analog for braid groups
of Springer's theory of regular elements.
This amounts to a conjectural description
of periodic elements (elements with a central
power) and their centralizers. When $W$ is the symmetric group,
periodic elements in $B(W)$ may be understood thanks to
Ker\'ekj\'art\'o's theorem on periodic homeomorphisms of the disk. In the
more general setting of spherical type Artin
groups, finding a simple description of periodic elements was an open
question.
We are able to solve these problems when $W$ is well-generated:
Theorem \ref{theoroots} contains a complete description of
the roots of the generator of the center of the pure braid group $P(W)$
and of their centralizers.

As for Theorem \ref{theo31}, the main conceptual ingredient
towards the proof of Theorem \ref{theoroots} is a general property of
Garside categories, explained in our separate paper \cite{cyclic}.
What is done here is the minor step consisting of re-interpreting the
general \emph{Ker\'ekj\'art\'o theorem for Garside categories} from
\cite{cyclic} in terms of the $S^1$-structure on the regular
orbit space $W\cq V^{\reg}$.

\smallskip

{\bf \flushleft Non-crossing partitions.}
A combinatorial by-product of our approach is a general
construction of \emph{generalized non-crossing partitions},
associated to each type of well-generated complex reflection groups.

In the classical cases $A_n$, $B_n$, $D_n$ and, more generally,
$G(e,e,n)$, the structure of $M(W)$ is understood in terms
of suitable notions of \emph{non-crossing partitions}
(\cite{athanasiadisreiner}, \cite{BC}, \cite{BDM}). The dual braid
monoid of an irreducible well-generated complex reflection $W$ gives rise
to a lattice of \emph{generalized non-crossing partitions},
whose cardinal is the \emph{generalized Catalan number}
$$\Cat(W) := \prod_{i=1}^n \frac{d_i+d_n}{d_i}.$$
(the term ``partition'' should not be taken too seriously: except for
the classical types,
lattice elements do not have natural interpretations as actual 
set-theoretic partitions.)
It is likely that this combinatorial object has some representation-theoretic
interpretation. In the ``badly-generated'' case, $\Cat(W)$ may fail to be
an integer, and the natural substitute for $NCP(W)$ is the
graph of simple elements of the dual braid category.

\section{Complex reflection groups, discriminants, braid groups}

Let $V$ be a vector space of finite dimension $n$.
A \emph{reflection group} in $\GL(V)$ is a
subgroup $W$ generated by \emph{(generalized) reflections},
\ie, elements whose fixed subspace is a hyperplane.
When the base field of $V$ is $\BC$, we say that $W$ is a \emph{complex
reflection group}. We are only interested in finite reflection groups and
will always assume finiteness, unless otherwise specified.

Let $W\subseteq \GL(V)$ be a complex reflection group.
A \emph{system of basic invariants} for $W$ 
\index{$f$ (system of basic invariants)}
is an $n$-tuple  $f=(f_1,\dots,f_n)$
of homogeneous generators of $\CO_V^W$, the algebra
of $W$-invariant polynomial functions on $V$.
A classical theorem of Shephard-Todd, \cite{shephard}, asserts
that such tuples exist, and that they consist of algebraically independent
terms.
Set $d_i:=\deg f_i$; these numbers are the \emph{degrees} of $W$.
Up to reordering, we may assume that
$d_1\leq d_2 \leq \dots \leq d_n$. The sequence $(d_1,\dots,d_n)$ is
then independent of the choice of $f$.

Choosing a system of basic invariants $f$ amounts to choosing
a graded algebra isomorphism
$\CO_V^W \simeq\BC[X_1,\dots,X_n], f_i \mapsto X_i$,
where the indeterminate
$X_i$ is declared homogeneous with degree $d_i$. Geometrically,
this isomorphism identifies the categorical quotient $W\cq V$
with the affine space $\BC^n$.

Further features of the invariant theory of complex reflection
groups involve invariant vector fields and invariant differential forms
on $V$. 

\begin{theo}[\cite{orlikterao}, Lemma 6.48]
The $\CO_V^W$-modules $(\CO_V\otimes V)^W$ and $(\CO_V\otimes V^*)^W$ are
free of rank $n$.
\end{theo}

If $f=(f_1,\dots,f_n)$ is a system of basic invariants, $df:=(df_1,\dots,df_n)$
is a $\CO_V^W$-basis for $(\CO_V\otimes V^*)^W$.
Being homogeneous, the module $(\CO_V\otimes V)^W$ admits an homogeneous
basis. 

\begin{defi}
A system of basic derivations for $W$ is an homogeneous $\CO_V^W$-basis
$\xi=(\xi_1,\dots,\xi_n)$ of $(\CO_V\otimes V)^W$,
with $\deg(\xi_1) \geq \deg(\xi_2) \geq \dots \geq \deg(\xi_n)$.

The sequence $(d_1^*,\dots,d_n^*):=
(\deg(\xi_1),\dots,\deg(\xi_n))$ is the sequence of \emph{codegrees} of $W$
(it does not depend on the choice
of $\xi$).
\end{defi}

Note that, as in \cite{zariski}, we label codegrees in decreasing order,
which is slightly unusual.
When $W$ is a complexified real reflection group, we have $V\simeq V^*$
as $W$-modules, thus $d_i^* = d_{n-i+1}-2$ for all $i$. This relation
is specific to the real situation and is not relevant here.

The Euler vector field on $V$ is invariant and of degree $0$.
Thus $d_n^*=0$.

Invariant vector fields define vector fields on the quotient variety.
Let $f$ be a system of basic invariants and $\xi$ be a system of
basic derivations.
For $j\in\{1,\dots,n\}$, the vector field $\xi_j$ defines a vector
field $\overline{\xi_j}$ on $W \cq V$. Since $\frac{\partial}{\partial f_1},
\dots,\frac{\partial}{\partial f_n}$ is a $\CO_V^W$-basis of the module
of polynomial vector fields on $W \cq V$, we have
$$\overline{\xi_j} = \sum_{i=1}^n m_{i,j} \frac{\partial}{\partial f_i},$$
where the $m_{i,j}$ are uniquely defined elements of $\CO_V^W$.

\begin{defi}
The discriminant 
matrix of $W$ (with respect to $f$ and $\xi$) is $M:=(m_{i,j})_{i,j}$.
\end{defi}

By weighted homogeneity, one has:
\begin{lemma}
\label{wtmij}
For all $i,j$, $\wt(m_{i,j})= d_i + d_j^*$.
\end{lemma}

The vector space $V$ decomposes as a direct
sum $\bigoplus_i V_i$ of irreducible representations of $W$.
Denote by $W_i$ the irreducible reflection group in $\GL(V_i)$
generated by (the restriction of) the reflections in $W$
whose hyperplanes contain $\bigoplus_{j\neq i}V_j$.
We have $W \simeq \prod_{i} W_i$.
Viewing $W$ and the $W_i$'s as
\emph{reflection groups}, \ie, groups endowed with a reflection representation,
it is natural to actually write $W= \bigoplus_iW_i$.

We denote by $\ACA$
\index{$\ACA$ (reflection arrangement)}
\index{reflection arrangement}
 the arrangement of $W$, \ie, the
set of reflecting hyperplanes of reflections in $W$.
We set $$V^{\reg}:=V-\bigcup_{H\in \ACA}H.$$
Denote by $p$ the quotient map $V\twoheadrightarrow W \cq V$.
Choose a basepoint $v_0\in V^{\reg}$.

\begin{defi}[\cite{bmr}] The \emph{braid group} of $W$
\index{$\Braid$ (braid group of $W$)}
\index{braid group}
is $B(W):=\pi_1(W\cq V^{\reg},p(v_0))$.
\end{defi}
Later on, when working with well-generated reflection groups,
we will slightly upgrade this definition, by replacing the basepoint
by a convenient contractible subspace of $W\cq V^{\reg}$ (see Definition
\ref{newbraidgroup}).

To write explicit equations, one chooses a system of 
basic invariants $f$. The \emph{discriminant}
$\Delta(W,f)\in\BC[X_1,\dots,X_n]$
\index{$\Delta$ (discriminant equation)}
\index{discriminant}
is
the reduced equation of $p(\bigcup_{H\in \ACA}H)$, via the 
identification $\CO_V^W \simeq\BC[X_1,\dots,X_n]$.

One easily sees that $B(W)\simeq \prod_{i} B(W_i)$. More generally,
all objects studied here behave ``semi-simply'', and we may restrict
our attention to irreducible complex reflection groups.

Since $B(W)$ is the fundamental group of the complement of an 
algebraic hypersurface, it is generated by particular elements
called \emph{generators-of-the-monodromy} or \emph{meridians}
(see, for example, \cite{bmr} or \cite{zariski}). They map
to reflections under the natural epimorphism $B(W)\rightarrow W$.
The diagrams given in \cite{bmr} symbolise presentations whose generators
are generators-of-the-monodromy (except for the six exceptional types
for which no presentation was known).
\begin{defi}
The generators-of-the-monodromy of $B(W)$ are called \emph{braid reflections}.
\index{braid reflections}
\end{defi}
This terminology was suggested by Brou\'e. It is actually tempting
to simply call them \emph{reflections}: since they generate $B(W)$, the
braid group appears to be some sort of (infinite) ``reflection group''.
This guiding intuition is quite effective.

Another natural feature of $B(W)$ is the existence of a natural
length function, \index{$l$|see{length function}} \index{length function! on $B(W)$}
which is the unique group morphism
$$l:B(W) \rightarrow \BZ$$
such that, for all braid reflection $s\in B(W)$, $l(s)=1$.

Consider the 
intersection lattice $\CL(\ACA):=\{ \bigcap_{H\in A} H | A\subseteq \ACA\}$.
Elements of $\CL(\ACA)$ are called \emph{flats}.
It is standard to endow $\CL(\ACA)$ with the reversed-inclusion partial
ordering:
$$ \forall L,L' \in \CL(\ACA), L \leq L' :\Leftrightarrow L \supseteq L'.$$

For $L \in \CL(\ACA)$, we denote by $L^0$ the complement in $L$ of
the flats strictly included in $L$.
The $(L^0)_{L \in \CL(\ACA)}$ form a stratification $\CS$ of $V$.
We consider the partial ordering on $\CS$ defined by 
$L^0 \leq L'^0 :\Leftrightarrow L \leq L'$. This is a degeneracy relation:
$$\forall L\in\CL(\ACA), \overline{L^0} = L= \bigcup_{S \in \CS, L^0 \leq S} S.$$

Since $W$ acts on $\ACA$, it acts on $\CL(\ACA)$ and we obtain 
a quotient stratification $\overline{\CS}$ of $W\cq V$ called
\emph{discriminant stratification}.

\begin{prop}[\cite{orlikterao}, Corollary 6.114]
\label{span}
Let $v\in V$. The vectors $\xi_1(v),\dots,\xi_n(v)$ span the
tangent space to the stratum of $\CS$ containing $v$.
The vectors $\overline{\xi}_1(v),\dots,\overline{\xi}_n(v)$ span the
tangent space to the stratum of $\overline{\CS}$ containing $\overline{v}$.
\end{prop}

Another chapter of the classical invariant
theory of complex reflection groups is Springer's theory of regular elements:

\begin{defi}
\label{defiregular}
\index{regular element|see{Springer theory}}
\index{Springer theory!regular element}
Let $W\subseteq \GL(V)$ be a complex reflection group.
Let $\zeta$ be a complex root of unity. An element
$w\in W$ is \emph{$\zeta$-regular} if
$\ker(w-\zeta)\cap V^{\reg}\neq \varnothing$. The eigenvalue $\zeta$ is
then called a \emph{regular eigenvalue} for $W$ and its order called a
\emph{regular number} for $W$.
\end{defi}

Note that, since $W$ acts freely on $V^{\reg}$, a $\zeta$-regular element
must have the same order as $\zeta$. The regularity of $\zeta$ only depends
on its order $d$, since the $k$-th power of a $\zeta$-regular element is
$\zeta^k$-regular.

The following theorem compiles some of the main features:

\begin{theo}
\label{theoregular}
\index{Springer theory}
Let $W\subseteq \GL(V)$ an irreducible complex reflection group,
with degrees $d_1,\dots,d_n$ and
codegrees $d_1^*,\dots,d_n^*$.

(1) Let $d$ be a positive integer. Set $$A(d):=\{i=1,\dots,n|\; d|d_i\}
\quad \text{and} \quad B(d):=\{i=1,\dots,n|\;d|d_i^*\}.$$
Then $|A(d)|\leq |B(d)|$, and $d$ is regular if and only if $|A(d)|=|B(d)|$.

(2) Let $w$ be a
$\zeta$-regular element of order $d$. Let $V':=\ker(w-\zeta)$.
The centralizer $W':=C_W(w)$, viewed in its natural representation in $\GL(V')$,
is a complex reflection group with degrees $(d_i)_{i\in A(d)}$ and
codegrees $(d_i^*)_{i\in B(d)}$. 

(3) Let $w$ be $\zeta$-regular element of order $d$.
Then $W'\cq V'\simeq (W\cq V)^{\mu_d}$ and
$W'\cq V'^{\reg}\simeq (W\cq V^{\reg})^{\mu_d}$ (where the $\mu_d$-action
is the quotient action of the scalar multiplication on $V$).
\end{theo}

Statement (2) was proved by Springer in his seminal paper \cite{springer}
(except for the part about codegrees, first observed by Denef-Loeser, \cite{denefloeser},
then conceptually proved by Brou\'e, \cite[5.19 (4)]{brouelnm}).
Statement (3) was
proved independently by Lehrer and Denef-Loeser. Statement (1) was
initially observed by Lehrer-Springer on a case-by-case basis, a
conceptual proof was given by Lehrer-Michel, \cite{lehrermichel}.

\begin{example}
\label{exampleg31}
Let $W:=G_{37}=W(E_8)$. It is a well-generated
complex reflection group in $\GL_8(\BC)$,
whose degrees are $$2,8,12,14,18,20,24,30.$$
By duality, the codegrees are $$0,6,10,12,16,18,22,28.$$ The integer $4$
is regular. The centralizer $W'$ is a complex reflection group of type
$G_{31}$ \index{$G_{31}$ (exceptional group)}.
Its degrees are $$8,12,20,24,$$
while the codegrees are $$0,12,16,28.$$
It is not a duality group, and it is not well-generated (see next section).
\end{example}

\section{Well-generated complex reflection groups}
\label{sectionwellgen}

Irreducible complex reflection groups
were classified fifty years ago by Shephard and Todd, \cite{shephard}.
There is an infinite family $G(de,e,n)$, where $d,e,n$ are positive
integers, and $34$ exceptions $G_4,\dots,G_{37}$.
Let us distinguish three subclasses of complex reflection groups:
\begin{itemize}
\item \emph{(complexified) 
real reflection groups}, obtained by scalar extension from
reflection groups of real vector spaces;
\item \emph{$2$-reflection groups},
 generated by reflections of order $2$;
\item \emph{well-generated reflection groups}, complex
reflection groups  $W\subseteq \GL(V)$ which can be generated
by $\dim_{\BC} V/V^{W}$ reflections, where $V^W:=\{v\in V | \forall w\in W,
wv=v\}$.
\end{itemize}
Real reflection groups are both $2$-reflection groups and well-generated.
For non-real groups, any combination of the other two properties may hold.

As far as the $K(\pi,1)$ conjecture and properties of 
braid groups are concerned, it is enough to restrict
one's attention to $2$-reflection groups:

\begin{defi}
\label{defisodisc}
Let $W\subseteq \GL(V)$ and $W'\subseteq \GL(V')$ be complex reflection groups.
We say that $W$ and $W'$ are \emph{isodiscriminantal} if one may
find systems of basic invariants $f$ (resp. $f'$) for $W$ (resp. $W'$)
such that $\Delta(W,f)=\Delta(W',f')$.
\end{defi}

When this happens, $W\cq V^{\reg}\simeq W'\cq V'^{\reg}$ and
$B(W)\simeq B(W')$.

\begin{theo}
\label{shephard}
\index{discriminant!it suffices to consider $2$-reflection groups}
Any complex reflection group is isodiscriminantal to a complex
$2$-reflection group.
\end{theo}

\begin{proof}
This may be observed on the classification and
was certainly known to experts.
In \cite{bmr}, Brou\'e-Malle-Rouquier associate to each complex
reflection group a diagram symbolizing a presentation by generators and relations;
they notice that the degrees and codegrees are invariants of the underlying braid diagram
(removing torsion relations from the presentation); actually, the braid diagrams
are invariants of isodiscriminantality classes (compare \cite{bmr} with \cite{orliksolomon}).
So the theorem can be rephrased as: for any diagram in the tables of
Brou\'e-Malle-Rouquier, the
diagram with the same braid relations but where all torsion relations have order $2$ is also 
in the tables —- this is an easy check.
\end{proof}

\begin{remark}
The work of Couwenberg-Heckman-Looijenga, \cite{chl}, coud possibly
be adapted
to provide a direct
argument. Let us sketch a conjectural way to proceed.
All references and notations are from \cite{chl}.
Assume that $W$ is not a $2$-reflection group.
For each $H\in \ACA$, let $e_H$ be the order of the pointwise stabilizer
$W_H$ and set $\kappa_H:=1 - e_H/2$. Consider the Dunkl connection
$\nabla$ with connection form $\sum_{H\in \ACA} \omega_H\otimes \kappa_H
\pi_H$, as in Example 2.5. Since $e_H\geq 2$, we have $\kappa_H\leq 0$.
In particular, $\kappa_0=1/n\sum_{H\in\ACA}\kappa_H \leq 0$ (Lemma 2.13)
and we are in the situation of \emph{loc. cit.} Section 5. 
In many cases, this suffices to conclude.
The problem is that, even though at least
some $e_H$ have to be $>2$, it is possible
that $\ACA$ contains several orbits, some of them with $e_H=2$. To handle
this, one has to enlarge the ``Schwarz symmetry group'' of
\emph{loc. cit.}, Section 4.
\end{remark}

The importance of the distinction between well-generated and
``badly-generated'' groups was first pointed out by Orlik-Solomon,
who observed in \cite{orso} 
a coincidence with invariant-theoretical aspects. Their observations
may be refined and completed as follows:
\begin{theo}
\label{monic}
Let $W$ be an irreducible complex reflection group. The following
assertions are equivalent:
\begin{itemize}
\item[(i)] $W$ is well-generated.
\item[(ii)] For all $i\in \{1,\dots,n\}$,
$d_i+d_i^* = d_n$.
\item[(iii)] For all $i\in \{1,\dots,n\}$, $d_i+d_i^* \leq d_n$.
\item[(iv)] For any system of basic invariants $f$, there exists
a system
of basic derivations $\xi$ such that the discriminant matrix decomposes as
$M= M_0 + X_n M_1$, where $M_0,M_1$ are matrices with coefficients
in $\BC[X_1,\dots,X_{n-1}]$ and $M_1$ is lower triangular with non-zero
scalars
on the diagonal.
\item[(v)] For any system of basic invariants $f$, we have
$\frac{\partial^n\Delta(W,f)}{(\partial X_n)^n}\in \BC^{\times}$
(in other words,
$\Delta(W,f)$, viewed as a polynomial in $X_n$ with coefficients
in $\BC[X_1,\dots,X_{n-1}]$, is monic of degree $n$).
\end{itemize}
\end{theo}

The matrix $M_1$ from assertion (iv) is an analogue of the matrix $J^*$
from \cite{saitoorbifold}, p. 10. Assertion (iv) itself generalizes the
non-degeneracy argument for $J^*$, which is an 
important piece of the construction of Saito's ``flat structure''.

\begin{proof}
(i) $\Rightarrow$ (ii) was observed in \cite{orso}
inspecting the classification. We still have no good explanation.

(ii) $\Rightarrow$ (iii) is trivial.

(iii) $\Rightarrow$ (v).
Let $h:=d_n$. A first step is to observe that,
under assumption (iii), $h$ is a regular number.
Indeed,
condition (iii) implies, for any $i=1,\dots,n-1$,
that $0<d_i< h$ and $0<d_i^*<h$, thus that
exactly one degree ($d_n$) and one codegree ($d_n^*$) are multiples of $h$,
thus $h$ is regular (Theorem \ref{theoregular} (1)).
Since $h$ is regular and divides only one degree, we may use
\cite{zariski}, Lemma 1.6 (ii) to obtain assertion (v):
the discriminant is $X_n$-monic, and by weighted-homogeneity it must be
of degree $n$.

(iii) $\Rightarrow$ (iv) is a refinement of the previous discussion.
Each entry $m_{i,j}$ of the matrix $M$ is weighted-homogeneous of
weight $d_i+d_j^* \leq d_n+d_1^*= 2h - d_1 < 2h$; since
$X_n$ has weight $h$, $\deg_{X_n}m_{i,j} \leq 1$. This explains
the decomposition $M = M_0 + X_n M_1$, where $M_0$ and $M_1$ have
coefficients in $\CO_Y$.
If $i<j$ and $d_i< d_j$, then $d_i+d_j^* < d_j + d_j^* = h$, thus $\deg_{X_n}
m_{i,j} = 0$. The matrix $M_1$ is \emph{almost lower triangular},
\ie, lower triangular except that
there could be non-zero terms above the diagonal in square diagonal
blocks corresponding to successive equal degrees (successive degrees
may indeed be equal, as in the example of type $D_4$, where the degrees
are $2,4,4,6$).

Let $i_0< j_0$ such that $d_{i_0} =d_{i_0+1} = \dots =d_{j_0}$
(looking at the classification, one may observe
 that this forces $j=i+1$; this observation is not used in the argument
below).
For all $i,j \in
\{i_0,\dots,j_0\}$, we have $d_i+d_j^* = h$. By weighted homogeneity,
this implies that the corresponding square block of $M_1$ consists of scalars.
The basic derivations $\xi_{i_0},\dots,\xi_{j_0}$ all have the same degree,
thus one is allowed to perform Gaussian elimination on the corresponding
columns of $M$. Thus, up to replacing $\xi$ by another system of basic
derivations $\xi'$, we may assume that $M_1$ is lower triangular.

The diagonal terms of $M_1$ must be scalars, once again by weighted homogeneity.
Assuming (iii), we already know that (v) holds. The determinant
of $M$ is $\Delta(X,f)$; (v) implies that the coefficient of $X_n^n$
is non-zero.
This coefficient is the product of the diagonal terms of $M_1$. We have
proved (iv).

(iv) $\Rightarrow$ (v) is trivial.

(v) $\Rightarrow$ (i) follows from the main result in \cite{zariski}.
\end{proof}

The following notion was considered in \cite{saito} for real
reflection groups. 

\begin{defi}
A system of basic derivations is \emph{flat} (with respect to 
$f$) if
the discriminant matrix may be written $M=M_0+X_n\Id$,
where $M_0$ is a matrix with coefficients in $\BC[X_1,\dots,X_{n-1}]$.
\end{defi}

\begin{coro}
\label{flat}
\index{well-generated reflection groups!admit flat systems of basic derivations}
Let $W$ be a well-generated irreducible reflection group, together
with a system of basic invariants $f$.
There exists a flat system of basic derivations.
\end{coro}

\begin{proof}
Let $\xi$ be any system of basic derivations.
Write $M=M_0+X_nM_1$, as in characterisation $(iv)$ from the above theorem.
The matrix $M_1$ is invertible
in $\GL_n(\BC[X_1,\dots,X_{n-1}])$.
The matrix $M_1^{-1} M = M_1^{-1}M_0 + X_n \Id$ represents
a flat system of basic derivations (weighted homogeneity is preserved
by the Gaussian elimination procedure).
\end{proof}

Contrary to what happens with real reflection groups, we may
not use the identification $V\simeq V^*$ to obtain a ``flat system
of basic invariants''.

Irreducible groups which are not well-generated may always
be generated by $\dim V + 1$ reflections. This fact has been observed
long ago, by case-by-case inspection, but no
general argument is known. In some sense, these badly-generated groups
should be thought of as affine groups.
The simplest example of a non-well-generated group is the group $G(4,2,2)$,
generated by
\index{well-generated reflection groups!example of a badly-generated group}
$$\begin{pmatrix} -1 & 0 \\ 0 & 1 \end{pmatrix}, \;
\begin{pmatrix} 0 & 1 \\ 1 & 0 \end{pmatrix}, \;
\begin{pmatrix} 0 & -i \\ i & 0 \end{pmatrix}.$$
Among high-dimensional exceptional complex reflection groups,
only $G_{31}$ is badly-generated (see Example \ref{exampleg31}).

In the sequel, several arguments are case-by-case.
The list of irreducible well-generated complex reflection groups
is given in Table 1.
The actual number of cases to consider depends on the type of result:
\begin{itemize}
\item Thanks to Theorem \ref{shephard}, for any statement about
$W\cq V^{\reg}$ and its topology, one may restrict one's attention to
groups generated by involutive reflections, which are listed in the first
column of the table.
\item Some results involve the actual structure of $W$ (e.g., Proposition
\ref{propcardhurwitz}), and groups with higher order reflections have to
be considered. These cases
are listed in the second column of the table, on the same
line as the corresponding $2$-reflection groups
(see Theorem \ref{shephard}).
\end{itemize}

\begin{table}
\label{tablecas}
\begin{tabular}{|c|c|c|}
\hline
\multicolumn{2}{|c|}{$2$-reflection groups}
& Other well-generated groups \\
\hline
Real  & $G(1,1,n)$ ($A_{n-1}$) & $G_4$, $G_8$, $G_{16}$, $G_{25}$, $G_{32}$ \\
\cline{2-3}
series  & $G(2,1,n)$ ($B_{n-1}$) & $G(d,1,n)$, $G_5$, $G_{10}$,
$G_{18}$, $G_{26}$  \\
\hline
\begin{tabular}{c} Real/complex \\ series \end{tabular} &
$G(e,e,n)$ (contains $D_n$ and $I_2(e)$) & $G_6$, $G_9$, $G_{14}$, $G_{17}$,
$G_{20}$, $G_{21}$ \\
\hline
Real   & $G_{23}$ ($H_3$), $G_{28}$ ($F_4$), $G_{30}$ ($F_4$) & \\
exceptions & $G_{35}$ ($E_6$), $G_{36}$ ($E_7$), $G_{37}$ ($E_8$) & \\
\hline
\begin{tabular}{c} Complex \\ exceptions \end{tabular} &
 $G_{24}$, $G_{27}$, $G_{29}$, $G_{33}$, $G_{34}$ & \\
\hline
\end{tabular}
\smallskip
\caption{Irreducible well-generated complex reflection groups.}
\index{well-generated reflection groups!classification}
\end{table}

\begin{lemma}
\label{lemmaparawellgen}
\index{well-generated reflection groups!parabolic subgroups are well-generated}
Let $W\subseteq \GL(V)$ be a well-generated complex reflection group.
Let $v\in V$.
Let $$V_v:=\bigcap_{H\in \ACA, \; v\in H} H,$$
$$W_v:=\{w\in W | wv=v\}.$$
Then $W_v$ may be generated by $\dim_{\BC}  V/V_v$ reflections.
In particular, $W_v$ is again a well-generated complex reflection group.
\end{lemma}

\begin{proof}
The fact that $W_v$ is again a complex reflection group is a classical
theorem due to Steinberg. The fact that $W_v$ is again well-generated
is easy to check on the classification (it follows for example
from Brou\'e-Malle-Rouquier's observation that their diagrams
in \cite{bmr} provide generating systems for representatives of
all conjugacy classes of $W_v$; when $W$ is real, no case-by-case
is needed, since $W_v$ is again real thus well-generated).
\end{proof}

\section{Symmetric groups, configurations spaces and classical braid groups}
\label{sectionconfig}

This section introduces some basic terminology and notations.
Everything here is classical and elementary.

Let $n$ be a positive integer.
The symmetric group $\mathfrak{S}_n$ may be viewed as reflection
group, acting on $\BC^n$ by permuting the canonical basis.
This representation is not irreducible. Let $H$ be the hyperplane
of equation $\sum_{i=1}^n X_i = 0$ (where $X_1,\dots,X_n$ is the dual canonical
basis of $\BC^n$). It is preserved by $\mathfrak{S}_n$, which acts
on it as an irreducible complex reflection group.

We have
$\BC[X_1,\dots,X_n]^{\mathfrak{S}_n} = \BC[\sigma_1,\dots,\sigma_n]$ and
$\CO_H^{\mathfrak{S}_n} = \BC[\sigma_1,\dots,\sigma_n]/\sigma_1$,
where $\sigma_1,\dots,\sigma_n$ are the elementary symmetric functions
on $X_1,\dots,X_n$. Set
$$\overline{E}_n:= \mathfrak{S}_n \cq \BC^n = \Spec \BC[\sigma_1,\dots,\sigma_n]$$
\index{$\En$! non-centered configuration space $\overline{E}_n$}
 and 
$$E_n := \mathfrak{S}_n \cq H = \Spec \BC[\sigma_1,\dots,\sigma_n]/\sigma_1.$$
\index{$\En$! centered configuration space $E_n$}
These spaces have more convenient descriptions in terms of multisets.

Recall that a \emph{multiset}
is a set $S$ (the \emph{support} of the multiset) together with
a map $m:S \rightarrow \BZ_{\geq 1}$ (the \emph{multiplicity}).
The \emph{cardinality} of such a multiset is $\sum_{s\in S} m(s)$
(it lies in $\BZ_{\geq 0} \cup \{\infty\}$).
If $(S,m)$ and $(S',m')$ are two multisets and if $S,S'$ are subsets
of a common ambient set, then we may define a \emph{multiset (disjoint)
union}
$(S,m)\cup (S',m')$, whose support is $S\cup S'$ and whose multiplicity
is $m+m'$ (where $m$, resp. $m'$, is extended by $0$ outside $S$, resp
$S'$).
If $(s_1,\dots,s_n)$ is a sequence of elements of a given set $S$,
we use the notation $\pmb{\{ } s_1,\dots,s_n\pmb{\} }$ (with brackets
in bold font) to refer
to the multiset $\bigcup_{i=1}^n (\{s_i\},1)$, \ie, the multiset
consisting of the $s_i$'s ``taken with multiplicities''.

Let $(x_1,\dots,x_n)\in\BC^n$. The associated $\mathfrak{S}_n$-orbit
is uniquely determined by 
$\pmb{\{ }x_1,\dots,x_n \pmb{\} }$.
This identifies $E'_n$ with the set of multisets of cardinality $n$
with support in $\BC$ (such multisets are called \emph{configurations
of $n$ points in $\BC$}).
The subvariety $E_n$, defined by $\sigma_1=0$, consists of
\emph{centered} configurations, \ie, configurations
 $\pmb{\{ }x_1,\dots,x_n\pmb{\} }$
satisfying $\sum_{i=1}^nx_i = 0$.
The natural inclusion $E_n \subseteq E'_n$ admits the retraction $\rho$
defined by $$\rho(\pmb{\{ }x_1,\dots,x_n\pmb{\} }):=
\pmb{\{ }x_1 - \sum_{i=1}^nx_i/n,\dots,x_n - \sum_{i=1}^nx_i/n\pmb{ \}}.$$
Algebraically, this corresponds to the identification of $\BC[\sigma_1,\dots,
\sigma_n]/\sigma_1$ with $\BC[\sigma_2,\dots,\sigma_n]$.

We find it convenient to use
configurations in $E'_n$ to represent
elements of $E_n$, implicitly working through $\rho$. 
E.g., in the proof of Proposition \ref{YCcontractible},
it makes sense to describe
a deformation retraction of a subspace of $E_n$ to a point in terms
of arbitrary configurations because the construction, which only
implies the relative values of the $x_i$'s,
is compatible with $\rho$. We adopt this viewpoint from now on,
without further justifications (compatibility will always be obvious).

Consider the lexicographic total ordering of $\BC$:
if $z,z'\in \BC$, we set
$$z \leq z' :\Leftrightarrow \left\{ \begin{matrix}
\re(z) < \re(z') & \text{or} \\
\re(z)=\re(z') \text{ and } \im(z) \leq \im(z'). \end{matrix} \right.$$

\begin{defi}
\label{defisupport}
\index{support (ordered support of a configuration)}
The \emph{ordered support} of an element of $E_n$ is the unique sequence
$(x_1,\dots,x_k)$ such that the set $\{x_1,\dots,x_k\}$ is the support and
$x_1 < x_2 < \dots < x_k$.
\end{defi}

We may uniquely represent an element of
$E_n$ by its ordered support $(x_1,\dots,x_k)$ and the sequence
$(n_1,\dots,n_k)$ of multiplicities at $x_1,\dots,x_k$.

The regular orbit space $E_n^{\reg}:=\mathfrak{S}_n\cq H^{\reg}$
\index{$\En$! regular centered configuration space $E_n^{\reg}$} 
consists of those multisets whose support
has cardinality $n$ (or, equivalently, whose multiplicity 
is constantly equal to $1$).
More generally, the strata of the discriminant stratification of $E_n$
are indexed by partitions of $n$: the stratum $S_{\lambda}$
associated with a
partition $\lambda=(\lambda_1,\lambda_2,\dots,\lambda_k)$, where
the $\lambda_j$'s are integers with
$\lambda_1\geq \lambda_2 \geq \dots \geq \lambda_k > 0$
and  $\sum_{j=1}^k\lambda_j=n$,
consists of configurations whose supports have cardinality $k$
and whose multiplicity functions take the
values $\lambda_1,\dots,\lambda_k$ (with multiplicities).

The braid group $B_n$ associated with $\mathfrak{S}_n$ is the usual
braid group on $n$ strings. We need to be more precise about our
choice of basepoint. For this purpose, we define
$$E_n^{\gen}$$ as
\index{$\En$! generic centered configuration space $E_n^{\gen}$} 
the subset of $E_n^{\reg}$ consisting of configurations of $n$
points with distinct real parts (this is the first in a series of definitions
of semi-algebraic nature). It is clear that:
\begin{lemma}
\label{lemmaEngen}
$E_n^{\gen}$
is contractible.
\end{lemma}

Using our topological conventions, we set $$B_n:=\pi_1(E_n^{\reg},
E_n^{\gen}).$$
This group admits a \emph{standard generating set} (the one
considered by Artin), consisting of
braid reflections $\bm{\sigma}_1,\dots,\bm{\sigma}_{n-1}$ 
defined as follows (we used bold fonts to avoid confusion with the elementary
symmetric functions).
Let $(x_1,\dots,x_n)$ be the ordered support of a point in $E_n^{\gen}$.
Then $\bm{\sigma}_i$
is represented by the following motion of the support:
$$
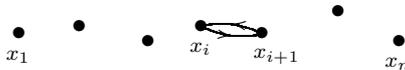
\xy
(-10,0)="1", (-2,1)="2", (7,-1)="3", 
(14,1)="4", (22,0)="5", (32,3)="6", (40,-1)="7",
"1"*{\bullet},"2"*{\bullet},"3"*{\bullet},"4"*{\bullet},
"5"*{\bullet},"6"*{\bullet},"7"*{\bullet},
(14,-2)*{_{x_i}},(24,-3)*{_{x_{i+1}}},
(-10,-3)*{_{x_1}},(40,-4)*{_{x_{n}}},
"4";"5" **\crv{(18,-1)}
 ?(.5)*\dir{>},
"5";"4" **\crv{(18,2)}
 ?(.5)*\dir{>}
\endxy $$

Artin's presentation for $B_n$ is
$$B_n  = \left< \bm{\sigma}_1,\dots,\bm{\sigma}_{n-1}
 \left| \bm{\sigma}_i\bm{\sigma}_{i+1}
\bm{\sigma}_i = \bm{\sigma}_{i+1}\bm{\sigma}_i
\bm{\sigma}_{i+1}, \bm{\sigma}_i\bm{\sigma}_j = 
\bm{\sigma}_j\bm{\sigma}_i \text{ if } |i-j| > 1
 \right. \right>.$$
The following definition requires a compatibility condition which
is a classical elementary consequence of the above presentation.

\begin{defi}
Let $G$ be a group. The \emph{(right)
Hurwitz action of $B_n$ on $G^n$}
is defined by
$$(g_1,\dots,g_{i-1},g_i,g_{i+1},g_{i+2},\dots,g_n) \cdot \bm{\sigma}_i
:= (g_1,\dots,g_{i-1},g_{i+1},g_{i+1}^{-1}g_ig_{i+1},g_{i+2},\dots,g_n),$$
for all $(g_1,\dots,g_n)\in G^n$ and all $i\in \{1,\dots,n-1\}$.
\end{defi}

This action preserves the fibers of the product
map $G^n\rightarrow G,(g_1,\dots,g_n)
\mapsto g_1\dots g_n$.

\section{Affine Van Kampen method}
\label{sectionvankampen}

This sections contains some generalities about Zariski-Van Kampen
techniques.
Let $P \in \BC[X_1,\dots,X_n]$ be a reduced polynomial.
What we have in mind is that $P$ is the discriminant of a complex reflection
group, but it does not cost more to work in a general context.
Let $\CH$ be the hypersurface of $\BC^n$ defined by $P=0$. Van Kampen's method
is a strategy for computing a presentation of $$\pi_1(\BC^n - \CH).$$

We assume that $P$ actually involves $X_n$ and write
$$P = \alpha_0 X_n^d +\alpha_1 X_n^{d-1} + 
\alpha_2 X_n^{d-2} + \dots + \alpha_d,$$
where $d$ is a positive integer (the \emph{degree in $X_n$}) and
 $\alpha_0,\dots,\alpha_d\in \BC[X_1,\dots,X_{n-1}]$, with
$\alpha_0\neq 0$.
We say that $P$ is \emph{$X_n$-monic}
if $\alpha_0$ is a scalar; when one is only
interested in the hypersurface defined by $P$, it is convenient to then
renormalise $P$ to have $\alpha_0=1$.
We denote by $\Disc_{X_n}(P)$ the \emph{discriminant} of $P$ with respect
to $X_n$, \ie, the resultant of $P$ and $\frac{\partial P}{\partial X_n}$.

This discriminant is a non-zero element of $\BC[X_1,\dots,X_{n-1}]$ and defines
an hypersurface $\CK$ in $Y:=\Spec\BC[X_1,\dots,X_{n-1}] \simeq \BC^{n-1}$.
Let $p$ be the natural projection
$\Spec\BC[X_1,\dots,X_n] \twoheadrightarrow Y$.
\index{$Y$|see{Saito quotient}}
\index{Saito quotient}

\begin{defi}
\index{$\CK$ (bifurcation locus)}
\index{bifurcation locus}
The \emph{bifurcation locus} of $P$ (with respect to the projection $p$)
is the algebraic hypersurface $\CK\subseteq Y$
defined by the equation $\Disc_{X_n}(P)=0$.
\end{defi}

\begin{defi}
A point $y\in Y$ is said to be \emph{generic} if it is not in $\CK$;
a \emph{generic line of direction $X_n$} is the fiber $L_y$ of $p$ over
a generic point $y\in Y$. 
\end{defi}

We represent points of $\BC^n$ by pairs $(y,z)$ where $z$ is the value of
the $X_n$ coordinate and $y$ is the image of the point under $p$.
Let $E:=p^{-1}(Y-\CK) \cap (\BC^n - \CH)$. The projection $p$ restricts
to a locally trivial fibration $E \rightarrow Y-\CK$, whose fibers 
are complex lines with $d$ points removed.
Choose a basepoint $(y,z)\in E$, let $F$ be the fiber containing $(y,z)$.

The following basic lemma was brought to my attention by Deligne and should
certainly have been included in my earlier paper \cite{zariski}.

\begin{lemma}
If $P$ is $X_n$-monic,
the fibration $p:E \rightarrow Y-\CK$ is split.
\end{lemma}

\begin{proof}
Assume $\alpha_0=1$.
Consider $\phi: Y \rightarrow \BC, y\mapsto 1+ \sum_{i=1}^{d} | \alpha_i(y)|$.
Then the map $Y-\CK \rightarrow E, y\mapsto (y,\phi(y))$ is a splitting
since by construction $P(y, \phi(y))$ is always non-zero.
\end{proof}

As a consequence, the fibration long exact sequence breaks into split short
exact sequences. Consider the commutative diagram:
$$\xymatrix{ 1 
\ar[r]  & \pi_1(F,(y,z)) \ar[r]^{\iota_*} \ar[dr]_{\alpha} 
& \pi_1(E,(y,z)) \ar[r]^{p_*} \ar[d]_{\beta} & \pi_1(Y -\CK,y) \ar[r] & 1 \\
   & & \pi_1(\BC^n-\CH,(y,z)).}$$
whose first line is a split exact sequence (the end of the fibration
exact sequence) and $\beta$ comes from the inclusion
of spaces.

We are in the context of \cite[Theorem 2.5]{zariski}, from which we
conclude that $\alpha$ is surjective.
We can actually be more precise and write a presentation for
$\pi_1(\BC^n-\CH,(y,z))$.
The semi-direct product structure of $\pi_1(E,(y,z))$ defines a 
morphism $\Phi: \pi_1(Y -\CK,y) \rightarrow \Aut ( \pi_1(F,(y,z)))$.

\begin{theo}[Van Kampen presentation]
\label{theovankampen}
Let $f_1,\dots,f_d$ be generators of
$\pi_1(F,(y,z)) \simeq F_d$ (the free group on $d$ generators). Let $g_1,\dots,g_m$ be generators
of $\pi_1(Y -\CK,y)$ with associated automorphisms $\phi_j:=\Phi(g_j)$.
We have
$$\pi_1(\BC^n-\CH,(y,z)) \simeq
\left< f_1,\dots,f_d \left|
f_i = \phi_j(f_i), \; 1\leq i \leq d,\; 1\leq j \leq m \right. \right>.$$
\end{theo}

\begin{proof}
Using the semi-direct product structure, we have the presentation
$$ \pi_1(E,(y,z)) \simeq 
\left< f_1,\dots,f_d,g_1,\dots,g_m \left|
g_jf_i = \phi_j(f_i)g_j, \; 1\leq i \leq d,\; 1\leq j \leq m \right. \right>.$$
One concludes by observing that $g_1,\dots,g_m$ may be chosen to be
meridians (``generators-of-the-monodromy'') around the irreducible components
of $p^{-1}(\CK)$: by \cite[Lemma 2.1.(ii)]{zariski}, $\ker \beta$
is generated as a normal subgroup by those meridians.
\end{proof}

\begin{coro}[Explicit Zariski $2$-plane section]
Let $P$ be a reduced polynomial in $\BC[X_1,\dots,X_n]$.
Let $a_1,\dots,a_{n-2}\in \BC$. Assume that:
\begin{itemize}
\item[(i)] $P$ is $X_n$-monic,
\item[(ii)] the coefficients of $\Disc_{X_n}(P)$ viewed as polynomials with
variable $X_{n-1}$ and coefficients in $\BC[X_1,\dots,X_{n-2}]$ are
altogether coprime (this in particular holds when 
$\Disc_{X_n}(P)$ is $X_{n-1}$-monic),
\item[(iii)] $\Disc_{X_{n-1}}(\Disc_{X_n}(P))(a_1,\dots,a_{n-2}) \neq 0$.
\end{itemize}
Let $\CH$ be the hypersurface of $\BC^n$ with equation $P=0$,
let $\Pi$ be the affine $2$-plane in $\BC^n$ defined by $X_1=a_1,\dots,
X_{n-2}=a_{n-2}$.
Then the map $\Pi\cap (\BC^n - \CH) \hookrightarrow \BC^n - \CH$ is a
$\pi_1$-isomorphism.
\end{coro}

\begin{proof}
Consider the affine line $p(\Pi) \subseteq Y$. Condition (iii) expresses
that $p(\Pi)$ is a \emph{generic line of direction $X_{n-1}$} in $Y$ and,
under condition (ii), \cite[Theorem 2.5]{zariski}, the injection
$p(\Pi)\cap (Y-\CK) \hookrightarrow Y -\CK$ is $\pi_1$-surjective. In
particular, in Theorem \ref{theovankampen}, we can choose loops drawn in
$\Pi \cap E$ to represent the lifted generators of $\pi_1(Y-\CK)$.
\end{proof}

Combined with our earlier work with Jean Michel, this corollary suffices
to prove Theorem \ref{theointropres}, since 
the $2$-plane sections described in \cite{BM} satisfy conditions (i)--(iii).
Presentations for $\pi_1(\Pi\cap(\BC^n-\CH))$ were obtained in \cite{BM},
using our software package VKCURVE.
At this stage, Theorem \ref{theointropres} relies on brutal
computations.
The sequel will provide a much more
satisfying approach at least for $W \neq G_{31}$.

\section{Lyashko-Looijenga coverings}
\label{sectionlyashkolooijenga}

Let $W$ be an irreducible well-generated complex reflection group, together with
a system of basic invariants $f$ and a flat system of basic
derivations $\xi$ with discriminant matrix $M_0+X_n\Id$.
Expanding the determinant, we observe that
$$\Delta_f = \det ( M_0+X_n \Id)= 
 X_n^n +\alpha_2 X_n^{n-2} + 
\alpha_3 X_n^{n-3} + \dots + \alpha_n,$$
where $\alpha_i\in \BC [X_1,\dots,X_{n-1}]$.
Since $\Delta_f$ is weighted homogeneous of total weight $nh$ for the system
of weights $\wt(X_i)=d_i$, each $\alpha_i$ is weighted homogeneous of
weight $ih$. 

\begin{defi}
\index{Lyashko-Looijenga morphism}
\index{$\LL$ (Lyashko-Looijenga morphism)}
The (generalized) \emph{Lyashko-Looijenga morphism} is the
morphism $\LL$ from $Y = \Spec\BC[X_1,\dots,X_{n-1}]$ to $E_n\simeq
\Spec\BC[\sigma_2,\dots,
\sigma_n]$ defined by $\sigma_i \mapsto (-1)^i\alpha_i$.
\end{defi}

This is of course much better understood in the following
geometric terms.
Via the choice of a system of basic invariants,
we have chosen an isomorphism
$$W\cq V \simeq  Y\times\BC.$$ Let $v\in V$. The orbit $\overline{v} \in W\cq V$ is represented by a pair
$$(y,z) \in  Y\times\BC,$$
where $z=f_n(v)$ and $y$ is the point in $Y$ with coordinates
$(f_1(v),\dots,f_{n-1}(v))$.
This encoding of points in $W\cq V$ will be used throughout this article.
As in the previous section,
we study the space $W\cq V\simeq Y\times\BC$ according to the fibers of the
the projection  $p:W\cq V \rightarrow Y, (y,z) \mapsto y$.
\begin{defi} For any point $y$ in $Y$, we denote
by $L_y$ the fiber of the projection
$p:W\cq V \rightarrow Y$ over $y$.
\end{defi}
For any $y\in Y$, the affine line $L_y$
intersects the discriminant $\CH$ in $n$ points
(counted with multiplicities), whose coordinates are
$$(y,x_1),\dots,(y,x_n),$$
where $\pmb{\{ }x_1,\dots,x_n\pmb{\} }$ is the multiset of solutions in $X_n$
of $\Delta_f=0$
where each $\alpha_i$ has been replaced by its value at $y$.
We have $$\LL(y) = \pmb{\{ }x_1,\dots,x_n\pmb{\} }.$$

The bifurcation locus $\CK\subseteq Y$ (see previous section) corresponds
precisely to those $y$ such that $\LL(y)$ contains multiple points.

The main theorem
of this section generalizes earlier results from \cite{looijenga}:

\begin{theo}
\label{theocovering}
The polynomials $\alpha_2,\dots,\alpha_{n}\in\BC[X_1,\dots,X_{n-1}]$ are algebraically independent and $\BC[X_1,\dots,X_{n-1}]$ is a
free graded $\BC[\alpha_2,\dots,\alpha_{n}]$--module of rank
$n!h^n/|W|$.
As a consequence, $\LL$ is a finite morphism.
It restricts to an unramified covering 
$Y-\CK \twoheadrightarrow E_n^{\reg}$ of degree $n!h^n/|W|$.
\end{theo}

I thank Eduard Looijenga for precious help with the theorem.
A prior version of this text contained a gap (the key Lemma
\ref{lemmaquasifinite}) and the very nice argument below is due to him.

\begin{lemma}
\label{lemmalooij}
Let $v\in V$ with image $\overline{v}\in W\cq V$.
Let $$V_v:=\bigcap_{H\in \ACA, \; v\in H}H.$$
The multiplicity of $\CH$ at $\overline{v}$ is 
$\dim_{\BC} V/V_v$.
\end{lemma}

\begin{proof}
At $v=0$, the multiplicity is the valuation of $\Delta_f$,
which is indeed $n$.

When $v\neq 0$, we consider the parabolic subgroup $W_v:=\{w\in W| wv=v\}$.
By Lemma \ref{lemmaparawellgen}, $W_v$ is again a well-generated complex
reflection group. The quotient map $\bigcup_{H\in \ACA}H \twoheadrightarrow
\CH = W \cq (\bigcup_{H\in \ACA}H)$ factors through
$\bigcup_{H\in \ACA} H \twoheadrightarrow
W_v\cq (\bigcup_{H\in \ACA}H)$. Because $W_v\cq (\bigcup_{H\in \ACA}H)\twoheadrightarrow
W \cq (\bigcup_{H\in \ACA}H)$
is unramified over $\overline{v}$, the multiplicity
of $\CH$ at $\overline{v}$ coincides with the multiplicity
of $W_v\cq(\bigcup_{H\in \ACA}H)$ at the image $\tilde{v}$ of $v$.
Around $\tilde{v}$, $W_v\cq(\bigcup_{H\in \ACA}H)$ is the same
as $W_v\cq(\bigcup_{H\in \ACA,\; v\in H}H)$. This hypersurface
is a direct product of $V_v$ with the discriminant of $W_v$.
The multiplicity at $\tilde{v}$ is the multiplicity at the origin
of the discriminant of $W_v$. After reduction to the irreducible
case, we apply the already solved case: the
multiplicity is the rank of $W_v$ or, in other words, $\dim_{\BC} V/V_v$.
\end{proof}

\begin{remark}
As it was suggested by Referee \#4, Lemma \ref{lemmalooij} is actually a 
characterization of well-generated reflection groups:
as pointed out in \cite[Proposition 4.2]{zariski}
the minimum
number of reflections needed to generate an irreducible complex reflection
group is equal to the valuation of the discriminant (\ie, the degree of the smallest degree
monomial), which is the same as the multiplicity at $0$: so when $W$ isn't well-generated
the multiplicity at $0$ is greater than $\dim_{\BC}(V)$.
\end{remark}

\begin{lemma}
\label{lemmaquasifinite}
$\LL^{-1}(0) = \{0\}$.
\end{lemma}

In this statement, the $0$ on the left denotes the multiset with $n$ copies of $0$.

\begin{proof}
Let $C$ be the tangent cone to $\CH$ in $\BC^n\simeq  Y\times\BC$.
It is a closed subvariety
of the tangent bundle to $\BC^n$. Because 
$\CH$ is quasi-homogeneous and $\forall i\in\{1,\dots,n-1\}, \wt(X_n)> \wt(X_i)$,
the cone $C$ is ``horizontal'' at $0$. In particular, the (fiber over $0$
of the tangent cone to the) ``vertical'' line
$L_0$ of equation $X_1=\dots = X_{n-1}=0$ is not in $C$.

Let $y= (x_1,\dots,x_{n-1})\in \LL^{-1}(0)$, \ie, such
that the line
 $L_y$ with $X_1=x_1,\dots,X_{n-1}=x_{n-1}$ intersects $\CH$ in only one point, $(y,0)\in Y\times\BC$. We want to prove
that $y=0$. Because $\CH$ is quasi-homogeneous, it is enough
to work in a neighborhood of the origin. In particular, we may
assume that $y$ is close enough to $0$ for $L_y$ to still be outside $C$.
Using a refined B\'ezout theorem
(\cite[Corollary 12.4]{fulton}), we  have
$$i((y,0),L_y.\CH; \BC^n) = m_{(y,0)}(\CH)$$
where
\begin{itemize}
\item $i((y,0),L_y.\CH; \BC^n)$ is the intersection multiplicity of $\CH$ and
$L_y$ at $(y,0)$. This is the order of $0$ as a root of the polynomial
$\Delta_f|_{X_1=x_1,\dots,X_{n-1}=x_{n-1}}$. By assumption, this is $n$.
\item $m_{(y,0)}(\CH)$ is the multiplicity of $\CH$ at $(y,0)$. Let
$v$ be a preimage of $(y,0)$ in $V$. By Lemma \ref{lemmalooij},
$m_{(y,0)}(\CH) =\dim_{\BC} V/V_v$.
\end{itemize}
Thus $\dim_{\BC}V/V_v = n$, $V_v=0$, $W_v=W$ and $v=0$.
\end{proof}

Because $\LL$ is quasi-homogeneous, Lemma \ref{lemmaquasifinite} implies
that $\LL$ is a finite (= quasi-finite and proper) morphism,
or, in other words, that 
$\BC[X_1,\dots,X_{n-1}]$ is a
finite graded $\BC[\alpha_2,\dots,\alpha_{n}]$--module.
In particular, $\alpha_2,\dots,\alpha_{n}$ are algebraically independent.

Because $\BC[X_1,\dots,X_{n-1}]$ is Cohen-Macaulay and
is finite over $\BC[\alpha_2,\dots,\alpha_{n}]$, it is a free
$\BC[\alpha_2,\dots,\alpha_{n}]$--module.
The rank may be computed by comparing Hilbert series.
Since each $X_i$ has weight $d_i$,
the Hilbert series of
$\BC[X_1,\dots,X_{n-1}]$ is $\prod_{i=1}^{n-1} \frac{1}{1-t^{d_i}}$.
Since each $\alpha_i$ has weight $ih$,
the Hilbert series of 
$\BC[\alpha_2,\dots,\alpha_n]$ is $\prod_{i=2}^n \frac{1}{1-t^{ih}}$.
The rank is the limit at $t \to 1$ of the quotient of these series, equal to
$$\prod_{i=1}^{n-1} \frac{(i+1)h}{d_i} = \frac{n!h^{n-1}}{d_1\dots d_{n-1}}=
\frac{n!h^n}{|W|}.$$

Theorem \ref{theocovering} now follows from the following generalization of
\cite[Theorem 1.4]{looijenga}:

\begin{lemma}
$\LL$ is \'etale on $Y-\CK$.
\end{lemma}

\begin{proof}
As mentioned in \cite[(1.5)]{looijenga}, the result will follow
if we prove that for all $y\in Y-\CK$ with $\LL(y)=\{x_1,\dots,x_n\}$,
the hyperplanes $H_1,\dots,H_n$ tangent to $\CH$ at the $n$
distinct points $(y,x_1),\dots,(y,x_n)$ are in general position.

To prove this, we use Proposition \ref{span}.
Each $H_i$ is spanned by $\overline{\xi}_1(y,x_i),\dots,
\overline{\xi}_n(y,x_i)$. Let $(\varepsilon_1,\dots,\varepsilon_n)$
be the basis of $(W \cq V)^*$ dual to
 $(\frac{\partial}{\partial X_1},\dots,
\frac{\partial}{\partial X_n})$. Let $l_i = \sum_j \lambda_{j,i} \varepsilon_j$
be a non-zero vector in $(W \cq V)^*$ orthogonal
to $H_i$. This amounts to taking a non-zero
column vector $(\lambda_{j,i})_{j=1,\dots,n}$
in the kernel of $M(y,x_i)$, or equivalently an eigenvector of $M_0(y)$
associated to the eigenvalue $-x_i$. By assumption, the $x_i$ are distinct.
The eigenvectors are linearly independent.
\end{proof}

The theorem has the following corollary (which we won't use).

\begin{coro}
\index{bifurcation locus!complement is $K(\pi,1)$}
The space $Y-\CK$ is a $K(\pi,1)$.
\end{coro}

For the sake of clarity, let us also mention:

\begin{coro}
\label{coromumumumu}
Let $v\in V$. Denote by $(y,z)\in Y\times\BC$ (identified with $W\cq V$)
the image of $v$. The following integers coincide:
\begin{itemize}
\item[(i)] The multiplicity of $z$ in $\LL(y)$.
\item[(ii)] The intersection multiplicity of $L_y$ with $\CH$ at $(y,z)$.
\item[(iii)] The multiplicity of $\CH$ at $(y,z)$.
\item[(iv)] The rank $\dim_{\BC} V/V_v$ of the parabolic subgroup $W_v$.
\end{itemize}
\end{coro}

\begin{proof}
The integer defined by (i) and (ii) are the same by their very definition.
The identity between (iii) and (iv) is Lemma \ref{lemmalooij}.
The argument used to prove Lemma \ref{lemmaquasifinite} also
shows the identity between (ii) and (iii).
\end{proof}

The discriminant stratification of $E_n$ yields a natural stratification
of $Y$: when $\lambda$ is a partition of $n$, the stratum
$Y_{\lambda}$ consists of points $y$ such that the multiplicities
of $\LL(y)$ are distributed according to $\lambda$.
Applying the corollary, one sees that the stratification
$Y=\bigsqcup_{\lambda\models n} Y_\lambda$ is the ``shadow'' of the
discriminant stratification restricted to $\CH$.

\section{Tunnels, labels and Hurwitz rule}

Let $W$ be an irreducible well-generated complex reflection group.
We keep the notations from the previous section.
Let $y\in Y$.
Let $U_y$ be the complement in $L_y$ of the vertical imaginary half-lines below
the points of $\LL(y)$, or in more formal terms:
$$U_y := \{(y,z) \in L_y | \forall x \in \LL(y),
\re(z)=\re(x) \Rightarrow \im(z)>\im(x)\}.$$

Here is an example where the support of $\LL(y)$ consists of $4$ points,
and $U_y$ is the complement of three half-lines:
$$
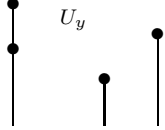
\xy
(-5,-10)="11",
(-5,0)="1", (-5,6)="2", (7,-4)="3", (7,-10)="33",
(14,2)="4", (14,-10)="44",
(3,4)*{_{U_y}},
"1"*{\bullet},"2"*{\bullet},"3"*{\bullet},"4"*{\bullet},
"11";"1" **@{-},
"2";"1" **@{-},
"33";"3" **@{-},
"44";"4" **@{-}
\endxy $$

``Generically'', $U_y$ is the complement of $n$ vertical half-lines.
We have to be careful about what ``generically'' means here:
a prerequisite
is that $\LL(y)$ should consist of $n$ distinct points, which amounts to
$y\in Y-\CK$ or equivalently $$\LL(y) \in E_n^{\reg},$$
but this is not enough: one needs these points
to be on distinct vertical lines, or equivalently that
$$\LL(y) \in E_n^{\gen}.$$

\begin{defi}
We set $Y^{\gen}:= \LL^{-1}(E_n^{\gen})$.
\end{defi}
This space, being the ``fiber'' of the covering $\LL$ over the ``basepoint''
of the basespace $E_n^{\reg}$,
is equipped with a Galois action of $B_n:=\pi_1(E_n^{\reg},E_n^{\gen})$.

\begin{defi}
\label{fatty}
\index{$\CU$ (fat basepoint)}
\index{fat basepoint}
The \emph{fat basepoint} of $W\cq V^{\reg}$ is the subset $\CU$ defined by
$$\CU:= \bigcup_{y\in Y} U_y,$$
or, equivalently,
by
$$\CU := \{ (y,z) \in Y \times \BC | \forall x \in \LL(y),
\re(z)=\re(x) \Rightarrow \im(z)>\im(x)\}.$$
\end{defi}

This definition calls for a few comments:

\begin{itemize}
\item First, as we will see if the next lemma and Definition
\ref{newbraidgroup} below, the subset $\CU$ can and will be used
``as if'' it was a genuine basepoint.
\item Note also that the definition implicitly relies on the choice
of a preferred direction in the complex line. The fat basepoint
is a truly semi-algebraic object and, geometrically, all constructions
below depend on the choice of a preferred element in the unit circle $S^1$.
\end{itemize}

\begin{lemma}
\index{fat basepoint!is dense, open and contractible}
The fat basepoint $\CU$ is dense in $W\cq V^{\reg}$, open and contractible.
\end{lemma}

\begin{proof}
The first two statements are clear.

Define a continuous function $\beta:Y\rightarrow \BR$ by
 $$\beta(y):=\max\{\im(x)|x \in \LL(y)\}+1.$$
Points of $W\cq V$ are represented by pairs $(y,z)\in Y\times\BC$,
or equivalently
by triples $(y,a,b)\in Y \times\BR\times \BR$,
where $a=\re(z)$ and $b=\im(z)$.
For $t\in [0,1]$, define $\phi_t: W\cq V \rightarrow W\cq V$ by
$$\phi_t(y,a,b) := \left\{ 
\begin{matrix} (y,a,b) & \text{if } b\geq \beta(y), \\
(y,a,b+t(\beta(y)-b)) & \text{if } b \leq \beta(y). \end{matrix}
\right.$$
Each $\phi_t$ preserves $\CU$ and
the homotopy $\phi$ restricts to a deformation retraction of
$\CU$ to $$\bigcup_{y\in Y} \{(y,z) \in L_y | \im(z) \geq \beta(y)\}.$$
The latter is a locally trivial bundle over the contractible space $Y$,
with contractible fibers (the fibers are half-planes). Thus it is contractible.
\end{proof}

As explained in Appendix \ref{appendixA}, we may
(and will) use $\CU$ as ``basepoint'' for $W \cq V^{\reg}$ and refine
our definition of $B(W)$:

\begin{defi}
\label{newbraidgroup}
The \emph{braid group of $W$} is $B(W):=\pi_1(W \cq V^{\reg},\CU)$.
\end{defi}

As for other notions actually depending on $W$,
 we often write $B$ instead of $B(W)$,
since most of the time we implicitly refer to a given $W$.

\begin{remark}
\label{remarkoverU}
We will need to consider a natural projection $\pi: B \rightarrow W$.
Recall that such a morphism is part of the fibration exact sequence
$$\xymatrix@1{ 1 \ar[r] & \pi_1(V^{\reg})
\ar[r] & \pi_1(W\cq V^{\reg} ) \ar[r]^{\phantom{mmm}\pi} & W \ar[r] & 1.}$$
For this exact sequence to be well-defined, one
has to make consistent choices of basepoints in $V^{\reg}$
and in $W\cq V^{\reg}$.
We have already described our ``basepoint''
$\CU$ in $W\cq V^{\reg}$. Choose $u\in \CU$ and choose a 
preimage $\tilde{u}$ of $u$ in $V^{\reg}$. If $u'\in \CU$ is another choice and
if $\gamma$ is a path in $\CU$ from $u$ to $u'$, $\gamma$ lifts
to a unique path $\tilde{\gamma}$ starting at $\tilde{u}$; since
$\CU$ is contractible, the fixed-endpoint 
homotopy class of $\tilde{\gamma}$ (and, in particular, its final point)
does not depend on $\gamma$. In other words, once we have chosen
a preimage of one point of $\CU$, we have a natural section $\tilde{\CU}$
of $\CU$ in $V^{\reg}$, as well as a transitive system of isomorphisms
between $(\pi_1(V^{\reg},\tilde{u}))_{\tilde{u}\in\tilde{\CU}}$.
From now on, we assume we have made
such a choice, and we define the pure braid group
as $\pi_1(V^{\reg},\tilde{\CU})$. This selects one particular morphism
$\pi$ (the $|W|$ possible choices yield conjugate morphisms).
\end{remark}

\begin{defi}
\label{defitunnel}
A \emph{semitunnel} is a
triple $T=(y,z,L)\in Y\times \BC \times \BR_{\geq 0}$ such
that $(y,z) \in \CU$ and the affine segment $[(y,z),(y,z+L)]$ lies in 
$W\cq V^{\reg}$.
The path $\gamma_T$ associated with $T$ is the path $t\mapsto (y,z+tL)$. 
The semitunnel $T$ is a \emph{tunnel} if in
addition $(y,z+L) \in \CU$.
\index{tunnels}
\index{tunnels!semitunnels}
\end{defi}

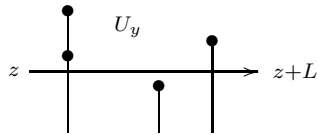
$$\xy
(-5,-10)="11",
(-5,0)="1", (-5,6)="2", (7,-4)="3", (7,-10)="33",
(14,2)="4", (14,-10)="44",
(3,4)*{_{U_y}},
"1"*{\bullet},"2"*{\bullet},"3"*{\bullet},"4"*{\bullet},
"11";"1" **@{-},
"2";"1" **@{-},
"33";"3" **@{-},
"44";"4" **@{-},
(-10,-2);(20,-2) **@{-}, *\dir{>},
(-12,-2)*{_{z}}, (25,-2)*{_{z+L}}
\endxy $$

The distinction between tunnels and semitunnels should be understood
in light of our topological conventions: if $T$ is a tunnel,
$\gamma_T$ represents
an element $$b_T\in \pi_1(W\cq V^{\reg},\CU),$$ 
while semitunnels will be
used to represent points of the universal cover
$(\UC(W \cq V^{\reg},\CU))_\CU$ (see Section \ref{sectionuniversalcover}).
 
\begin{defi}
\index{$S$ (set of simple elements in $B$)}
\index{simple elements}
\index{tunnels!represent simple elements}
An element $b\in B$ is \emph{simple} if $b=b_T$ for some tunnel $T$.
The set of simple elements in $B$ is denoted by $S$.
\end{defi}

\begin{remark}
Later on (Corollary \ref{Sisfinite}), we will show that $S$ is finite.
This might be disconcerting at first sight, as
$S$ is stable under particular conjugacy operations and one might be
misled into believing that $S$ is a union of conjugacy classes, which
it is not. The finiteness of $S$ and its direct description
in terms of the combinatorics of the finite group $W$
(Proposition \ref{posetiso}) is a key ingredient
of this paper.
\end{remark}
 
Each tunnel lives in a given $L_y$, where it may be represented by
an horizontal (= constant imaginary part) segment avoiding $\LL(y)$
and with endpoints in $U_y$.
The triple $(y,z,L)$ may be uniquely recovered
from $[(y,z),(y,z+L)]$. A frequent abuse of terminology will
consist of using the term \emph{tunnel} (or \emph{semitunnel})
to designate either the triple $(y,z,L)$, or the segment
$[(y,z),(y,z+L)]$, or the pair $(y,[z,z+L])$, depending on the context
(in particular, when intersecting tunnels with geometric objects,
the tunnels should be understood as affine segments).

Let $y\in Y$.
Let $(x_1,\dots,x_k)$ be the ordered support of $\LL(y)$.
The space $p_x((L_y \cap W\cq V^{\reg}) - U_y)$ is an union of
$k$ disjoint open affine intervals $I_1,\dots,I_k$,
where $$I_i := \left\{ \begin{matrix} (x_i-\sqrt{-1}\infty,x_i)
 & \text{if } 
i=1 \text{ or } (i> 1 \text{ and } \re(x_{i-1}) <\re(x_{i})), \\
(x_{i-1},x_i) & \text{otherwise.}
\end{matrix} \right.$$
(by $(x_i-\sqrt{-1}\infty,x_i)$, we mean the open vertical half-line below
$x_i$).
In the first case (when $I_i$ is not bounded),
we say that $x_i$ is \emph{deep}. In the picture below, there
are three deep points, $x_1$, $x_3$ and $x_4$.

$$
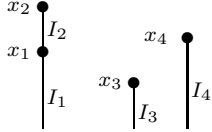
\xy
(-5,-10)="11",
(-5,0)="1", (-5,6)="2", (7,-4)="3", (7,-10)="33",
(14,2)="4", (14,-10)="44",
"1"*{\bullet},"2"*{\bullet},"3"*{\bullet},"4"*{\bullet},
(-8,0)*{_{x_1}},(-8,6)*{_{x_2}},(4,-4)*{_{x_3}},(10,2)*{_{x_4}},
(-3,-6)*{_{I_1}},(-3,3)*{_{I_2}},(9,-8)*{_{I_3}},(16,-5)*{_{I_4}},
"11";"1" **@{-},
"2";"1" **@{-},
"33";"3" **@{-},
"44";"4" **@{-}
\endxy $$

Choose a \emph{system of elementary tunnels} for $y$. By this, we mean the 
choice, for each $i=1,\dots,k$, of a small tunnel $T_i$ in $L_y$
crossing $I_i$
and not crossing the other intervals; let $s_i:=b_{T_i}$ be the associated 
element of $B$.

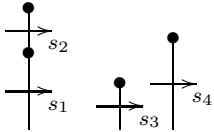
$$\xy
(-5,-10)="11",
(-5,0)="1", (-5,6)="2", (7,-4)="3", (7,-10)="33",
(14,2)="4", (14,-10)="44",
"1"*{\bullet},"2"*{\bullet},"3"*{\bullet},"4"*{\bullet},
"11";"1" **@{-},
"2";"1" **@{-},
"33";"3" **@{-},
"44";"4" **@{-},
(-8,-5);(-2,-5) **@{-}, *\dir{>}, (-1,-7)*{_{s_1}},
(-8,3);(-2,3) **@{-}, *\dir{>}, (-1,1)*{_{s_2}},
(4,-7);(10,-7) **@{-}, *\dir{>}, (11,-9)*{_{s_3}},
(11,-4);(17,-4) **@{-}, *\dir{>}, (18,-6)*{_{s_4}}
\endxy $$

These elements depend only on $y$ and not on the explicit choice
of elementary tunnels.

\begin{defi}
\label{defilabel}
\index{$\lbl$ (label map)}
\index{label}
The sequence $\lbl(y):=(s_1,\dots,s_k)$ is the \emph{label of $y$}.
Let ${i_1},{i_2},\dots,{i_l}$ be the indices of the successive 
deep points of $\LL(y)$.
The \emph{deep label} of $y$ is the subsequence $(s_{i_1},\dots,s_{i_l})$.
\end{defi}

In the above example, the deep label is $(s_1,s_3,s_4)$.
The length of the label is $n$ if and only if $y\in Y-\CK$. In this 
case, the
deep label coincides with the label if and only if $y \in Y^{\gen}$.

Later on, it will appear
that the pair $(\LL(y),\lbl(y))$ uniquely determines $y$
(Theorem \ref{theolabel}). 

\begin{remark}
\label{remarkygen}
When $y$ is generic, $\lbl(y)$ is an $n$-tuple of braid reflections
(because an elementary tunnel crossing the interval below a point
in $\LL(y)$ is essentially the same as a small circle around this point;
when $y$ is generic, the points in $\LL(y)$ correspond to smooth points
of the discriminant, and the elementary tunnels represent generators
of the monodromy).
\end{remark}

Consider the case $y=0$ (given by the equations $X_1=0,\dots,X_{n-1}=0$).
The multiset $\LL(y)$ has support $\{0\}$ with multiplicity $n$.

\begin{defi}
\label{defidelta}
\index{Garside element}
\index{$\delta$ (dual Garside element)}
We denote by $\delta$ the simple element such that $\lbl(0)=(\delta)$.
\end{defi}

This element plays the role of Deligne's element $\Delta$.
Choose $v\in V^{\reg}$ such that the $W$-orbit $\overline{v}$ lies in 
$L_0$. Brou\'e-Malle-Rouquier consider the element (denoted by $\bm{\pi}$, 
\cite[Notation 2.3]{bmr}) in the pure braid group $P(W)$ represented by the loop
\begin{eqnarray*}
[0,1] & \longrightarrow & V^{\reg}\\
t & \longmapsto & v\exp(2\sqrt{-1}\pi t).
\end{eqnarray*}
We prefer a different notation:

\begin{defi}
We call \emph{full-twist} this element of $P(W)$, and denote it by $\tau$.
\end{defi}

They observe that this element lies in the center of $B$
(\cite[Theorem 2.24]{bmr}) and conjecture that it generates the center
of $P$.

Since $X_n$ has weight $h$, $\delta^h$ coincides with $\tau$.
More precisely, $\delta$ is represented by the loop which is 
the image in
$W\cq V^{\reg}$ of the path in $V^{\reg}$
\begin{eqnarray*}
[0,1] & \longrightarrow & V^{\reg} \\
t & \longmapsto & v\exp(2\sqrt{-1}\pi t/h)
\end{eqnarray*}
In particular:

\begin{lemma}
\label{pideltaiscox1}
The element
$\tau=\delta^h$ is central in $B(W)$ and lies in $P(W)$. The image of
$\delta$ in $W$ is $\zeta_h$-regular, in the sense of Springer
(see Definition \ref{defiregular}).
\index{Springer theory!yields torsion elements in $B/ZB$}
\end{lemma}

See Theorem \ref{theoZ} for a description of the center of $B(W)$.

Each tunnel lives in a single fiber $L_y$. Let us now explain how
simple elements represented by
tunnels living in different fibers may be compared. The idea is
that, since being a tunnel is an
open condition, one may perturb $y$ without affecting the simple
element:

\begin{defi}
\label{defiTn}
Let $T=(y,z,L)$ be a tunnel. A \emph{$T$-neighbourhood} of $y$
is a path-connected neighbourhood $\Omega$ of $y$ in $Y$ such that,
for all $y'\in \Omega$, $T':=(y',z,L)$ is a tunnel.
\end{defi}

Such neighbourhoods clearly exist for all $y\in Y$.

\begin{lemma}[Hurwitz rule]
\index{Hurwitz rule}
Let $T=(y,z,L)$ be a tunnel, representing a simple element $s$.
Let $\Omega$ be a $T$-neighbourhood of $y$.
For all $y'\in \Omega$, $T':=(y',z,L)$ represents $s$.
\end{lemma}

\begin{proof}
This simply expresses that the tunnels $(y',z,L)$ and
$(y,z,L)$ represent homotopic paths, which is clear by definition of
$\Omega$.
\end{proof}

\begin{remark}
\label{remarkbasictunnel}
Let $T_1,\dots,T_k$
is a system of elementary tunnels for $y$. Let $\Omega_i$ be
a $T_i$-neighbourhood for $y$. A \emph{standard neighbourhood} of $y$
could be defined as a path-connected neighbourhood $\Omega$
of $y$ inside
$\cap_{i=1}^k \Omega_i$. These standard neighbourhoods form a basis
for the topology of $Y$.
A consequence of Hurwitz rule is that the label of $y$ may uniquely
be recovered once we know the label of a single $y'\in \Omega$.
Distinct $y'\in \Omega$ correspond to different ``desingularizations'' of $y$,
and their labels are obtained by further factorizing terms in the label of $y$.
Among them are full desingularizations (corresponging to factorizations in $n$ terms),
corresponding to choosing $y'$ in the non-empty intersection $\Omega \cap Y^{\gen}$.
\end{remark}

The remainder of this section consists of various consequences of
Hurwitz rule.

\begin{coro}
\label{corolengthlabel}
Let $y\in Y$. Let $(x_1,\dots,x_k)$ be the ordered support
of $\LL(y)$, and $(n_1,\dots,n_k)$ be the multiplicities.
For any $i$, the natural length of $s_i$ is
given by  $$l(s_i)=\sum_{j \text{ s.t. } \re(x_i)=\re(x_j) \text{ and }
\im(x_i) \leq \im(x_j) } n_j.$$
\end{coro}

\begin{proof}
The case $y\in Y^{\gen}$ is a consequence of Remark \ref{remarkygen}.
The general case follows by perturbing and applying Hurwitz rule.
\end{proof}

\begin{coro}
\label{corodeeplabel}
Let $y\in Y$.
Let $(s_{i_1},\dots,s_{i_l})$ be the deep label of $y$.
We have $s_{i_1}\dots s_{i_l}=\delta$.
\end{coro}

\begin{proof}
Any tunnel $T$ 
deep enough and long enough represents $s_{i_1}\dots s_{i_l}$.
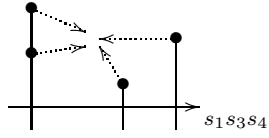
$$\xy
(-5,-10)="11",
(-5,0)="1", (-5,6)="2", (7,-4)="3", (7,-10)="33",
(14,2)="4", (14,-10)="44",
"1"*{\bullet},"2"*{\bullet},"3"*{\bullet},"4"*{\bullet},
"11";"1" **@{-},
"2";"1" **@{-},
"33";"3" **@{-},
"44";"4" **@{-},
"1";(2,1) **@{.}, *\dir{>},
"2";(2,3) **@{.}, *\dir{>},
"3";(4,1) **@{.}, *\dir{>},
"4";(4,2) **@{.}, *\dir{>},
(-8,-7);(17,-7) **@{-}, *\dir{>}, (22,-9)*{_{s_1s_3s_4}}
\endxy $$
The origin $0\in Y$ lies in
a $T$-neighbourhood of $y$. To conclude, apply the Hurwitz rule.
\end{proof}

Let us recall the following standard notion:

\begin{defi}[Hurwitz action]
Let $G$ be a group, let $B_n$ be the braid group on $n$ strings with its usual system of
generators $\bm{\sigma}_1,\dots,\bm{\sigma}_{n-1}$.
The \emph{Hurwitz action} of $B_n$ on $G^n$, denoted as a multiplication on the right,
is the unique group right action such that
$$(g_1,\dots,g_n) \cdot \bm{\sigma}_i = (g_1,\dots,g_{i-1},g_{i}g_{i+1}g_{i}^{-1},g_{i},g_{i+2},\dots,
g_n).$$
\end{defi}

In the following corollary, the notation $y \cdot \beta$ refers to the covering action
of $\pi_1(E_n^{\reg},x)$ on $\LL^{-1}(x)$.

\begin{coro}
\index{Hurwitz action}
Let $x \in E_n^{\gen}$.
Let $\beta\in \pi_1(E_n^{\reg},x)$, $y\in \LL^{-1}(x)$ and $y':= y \cdot
\beta$.
Let $(b_1,\dots,b_n)$ be the label of $y$ and 
$(b_1',\dots,b_n')$ be the label of $y'$.
Then $$(b_1',\dots,b_n') = (b_1,\dots,b_n) \cdot \beta,$$
where $\beta$ acts by right Hurwitz action. \index{Hurwitz action}
\end{coro}

\begin{proof}
It is enough to prove this for a standard generator $\bm{\sigma}_i$.
Let $(x_1,\dots,x_n)$ be the ordered support of $x$.
By Hurwitz rule, we may adjust the imaginary parts of the $x_i$'s without
affecting the label; in particular we may assume that $\im(x_i) < \im(x_{i+1})$.
We may find tunnels $T_-=(y,z_-,L)$ and $T_+=(y,z_+,L)$ as in the picture
below:
$$
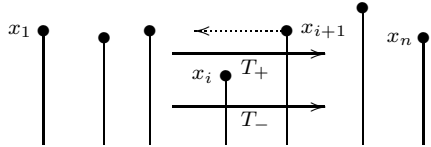
\xy
(-10,0)="1", (-2,-1)="2", (4,0)="3", 
(14,-6)="4", (22,0)="5", (32,3)="6", (40,-1)="7",
(-10,-15)="11", (-2,-15)="22", (4,-15)="33", 
(14,-15)="44", (22,-15)="55", (32,-15)="66", (40,-15)="77",
"1";"11" **@{-}, "2";"22" **@{-}, "3";"33" **@{-}, "4";"44" **@{-},
"5";"55" **@{-}, "6";"66" **@{-}, "7";"77" **@{-},
"1"*{\bullet},"2"*{\bullet},"3"*{\bullet},"4"*{\bullet},
"5"*{\bullet},"6"*{\bullet},"7"*{\bullet},
(11,-6)*{_{x_i}},(27,0)*{_{x_{i+1}}},
(-13,0)*{_{x_1}},(37,-1)*{_{x_{n}}},
"5";(10,0) **@{.}, *\dir{>},
(7,-3);(27,-3) **@{-}, *\dir{>},
(7,-10);(27,-10) **@{-}, *\dir{>},
(18,-12)*{_{T_-}},
(18,-5)*{_{T_+}}
\endxy $$
The path in $E_n^{\reg}$ where $x_{i+1}$ moves along the dotted arrow
and all other points are fixed represents $\bm{\sigma_i}$.
Applying Hurwitz rule to $T_+$, we obtain $b'_i=b_{i+1}$; applying
Hurwitz rule to $T_-$, we obtain $b'_ib'_{i+1}= b_ib_{i+1}$. The result
follows.
\end{proof}

\begin{coro}
\label{corocardhurwitz}
Let $y\in Y^{\gen}$. 
The cardinality of the Hurwitz orbit $\lbl(y) \cdot B_n$ is at most $n!h^n/|W|$,
and there is an equivalence between:
\begin{itemize}
\item[(i)] $|\lbl(y)\cdot B_n| = n!h^n/|W|$.
\item[(ii)] The orbits $y\cdot B_n$ and $\lbl(y)\cdot B_n$ are 
isomorphic as $B_n$-sets.
\item[(iii)] The map
$Y^{\gen} \rightarrow E_n\times B^n, y \mapsto (\LL(y),\lbl(y))$
is injective.
\end{itemize} 
\end{coro}

In the next section, we will prove that conditions (i)--(iii) actually hold.
This is not a trivial statement.

\begin{proof}
Let $G$ be a group and $\Omega,\Omega'$ be two $G$-sets, together
with a $G$-set morphism $\rho:\Omega\rightarrow \Omega'$.
Assume that $\Omega'$ is transitive. Then $\rho$ is surjective.
Assume in addition that $\Omega$ is finite. Then $\Omega'$ is finite,
 $|\Omega'| \leq |\Omega|$ and $$|\Omega'| = |\Omega| \Leftrightarrow
\text{$\rho$ is injective} \Leftrightarrow \text{$\rho$ is an isomorphism}.$$

By the previous corollary, one has $\lbl(y\cdot\beta)=\lbl(y)\cdot\beta$.
In other words, the map $\lbl$ extends to a $B_n$-set morphism $y\cdot B_n
\rightarrow \lbl(y) \cdot B_n$.
We apply our discussion to $G:=B_n$, $\Omega:=y\cdot B_n$ and
$\Omega':= \lbl(y) \cdot B_n$. Both $B_n$-sets are clearly transitive.
Since $y\in Y-\CK$, we have $|\LL^{-1}(\LL(y))|=|y\cdot B_n| = n!h^n/|W|$
(Theorem
\ref{theocovering}). We deduce that $|\lbl(y) \cdot B_n| \leq n!h^n/|W|$
and $\text{(i)} \Leftrightarrow \text{(ii)}$.

Assertion (iii) amounts to saying that, for all $y\in Y^{\gen}$,
 $\lbl$ is injective on the fiber
of $\LL$ containing $y$.
This fiber is precisely the orbit $y\cdot B_n$.
Under this rephrasing, it is clear that
$\text{(ii)} \Leftrightarrow \text{(iii)}$.
\end{proof}

\begin{coro}
\label{coros1si}
Let $s$ be a simple element. There exists $y\in Y^{\gen}$
and $i\in\{1,\dots,n\}$ such that
$s=s_1\dots s_i$,  where
$(s_1,\dots,s_n):=\lbl(y)$.
\end{coro}

\begin{proof}
Let $T$ be a tunnel representing $s$.
Any $T$-neighbourhood of $y$ contains generic points.
Up to perturbing $y$, we may assume that $y\in Y^{\gen}$. The picture below
explains, on an example, how
to move certain points (following the dotted paths) of
the underlying configuration to reach a suitable $y'$:
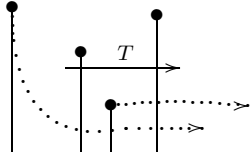
$$\xy
(-5,-13)="11",
(-5,0)="1", (-5,6)="2", (4,0)="3", (4,-13)="33",
(14,5)="4", (14,-13)="44", (8,-7)="5", (8,-13)="55",
"5"*{\bullet},"2"*{\bullet},"3"*{\bullet},"4"*{\bullet},
"33";"3" **@{-},
"11";"2" **@{-},
"55";"5" **@{-},
"44";"4" **@{-},
(2,-2);(17,-2) **@{-}, *\dir{>}, (10,0)*{_{T}},
"2";(10,-10) **\crv{~*=<4pt>{.} (-5,0),(-5,-13)},
(10,-10);(20,-10) **\crv{~*=<4pt>{.} (13,-10),(17,-10)},
*\dir{>},
"5";(26,-7) **\crv{~*=<4pt>{.} (16,-6)},
*\dir{>}
\endxy $$
This path in $E_n^{\reg}$
lifts, via $\LL$, to a path in a $T$-neighbourhood of $y$
whose final point $y'$ satisfies the conditions of the lemma.
\end{proof}

\section{Reduced decompositions of Coxeter elements}

\subsection{From braids to elements of $W$}

\begin{defi}
\index{well-generated reflection groups!generalized Coxeter element}
\index{$c$ (Coxeter element)}
\index{Coxeter element}
Let $W$ be an irreducible well-generated complex reflection group.
An element $c\in W$ is a
\emph{(generalized) Coxeter element} if it is $\zeta_h$-regular.
More generally, if $W$ is a well-generated complex reflection group
decomposed as a sum $W=\bigoplus_i W_i$ of irreducible groups,
a Coxeter element in $W$ is a product
$c=\prod_i c_i$ of Coxeter elements in each $W_i$.
\end{defi}

\begin{lemma}
\label{eigencox}
When $W$ is irreducible, a Coxeter element $c$ in $W$ has no non-trivial fixed point.
\end{lemma}

\begin{proof}
As shown by Springer, \cite{springer}, the eigenvalues of a $\zeta$-regular element
are $\zeta^{1-d_1},\dots,\zeta^{1-d_n}$.
Applying this to $\zeta=\zeta_h$ and noting that $0\leq d_1 \leq \dots \leq d_n=h$,
we obtain the desired result.
\end{proof}

When $W$ is irreducible, we may use the constructions from the previous
section and the morphism $\pi:B(W) \twoheadrightarrow W$ to obtain
 typical Coxeter elements (see 
Definition \ref{defidelta}):

\begin{lemma}
\label{pideltaiscox2}
When $W$ is irreducible,
the element $c:=\pi(\delta)$ is a Coxeter element in $W$.
\end{lemma}

\begin{proof}
This is a rephrasing of Lemma \ref{pideltaiscox1}.
\end{proof}
The other Coxeter elements, which are conjugates of $c$,
appear when considering other basepoints over $\CU$ (see Remark
\ref{remarkoverU}).

More generally, we have:

\begin{lemma}
\label{lemmaparacoxnew}
Let $y\in Y$. Let $(x_1,\dots,x_k)$ be the ordered support
of $\LL(y)$, let $(n_1,\dots,n_k)$ be the multiplicities.
Assume that $\re(x_1)< \dots < \re(x_k)$.
Let $(s_1,\dots,s_k):=\lbl(y)$.

For all $i$, set $c_i:=\pi(s_i)$. Then
there exists a preimage $v_i\in V$ of $(y,x_i)\in  Y\times\BC \simeq W\cq V$
such that $c_i$ is a Coxeter element in the parabolic subgroup $W_{v_i}$. 

In particular, if $n_i=1$, then $c_i$ is a reflection. 
\end{lemma}

\begin{proof}
When $n_i=n$ (thus $i=k=1$), the result is Lemma \ref{pideltaiscox2}.

When $n_i=1$, as pointed out in Remark \ref{remarkygen}, the element
$s_i$ is represented by a
 small loop around a smooth point in the discriminant, thus
is a braid reflection, and maps to a reflection in $W$.

The general case is similar: locally near $(y,x_i)$, the discriminant
is a direct product of $V_{v_i}$ with the discriminant of $W_{v_i}$
(see the proof of Lemma \ref{lemmalooij}). This local structure
provides a specific morphism $B(W_{v_i}) \to B(W)$, such that
the element ``$\delta_i$'' in $B(W_{v_i})$ (the product of the $\delta$'s
associated with each irreducible components of $W_{v_i}$) maps to $s_i$.
The lemma follows.
\end{proof}

The assumption $\re(x_1)< \dots < \re(x_k)$ may be removed at
the cost of replacing  $\pi(s_i)$ by $\pi(s_{i-1})^{-1}\pi_(s_i)$
when $\re(x_{i-1})=\re(x_i)$. This is behind
Definition \ref{defireducedlabel} below.

Let $R$ be the set of all reflections in $W$.
As in \cite{dualmonoid}, for all $w\in W$, we denote by $\Red_R(w)$ the
set of \emph{reduced $R$-decompositions of $w$}, \ie, minimal
length sequences of elements of $R$ with product $w$.
Since $R$ is closed under conjugacy, $\Red_R(c)$ is stable under Hurwitz
action.
We also consider the length function $l_R:W\rightarrow \BZ_{\geq 0}$,
whose value at $w$ is the common length of the elements of $\Red_R(w)$,
and two partial orderings of $W$
defined as follows: for all $w,w' \in W$, we set
\index{$\preccurlyeq_R$ ($R$-prefix partial ordering in $W$)}
$$w \preccurlyeq_R w' :\Leftrightarrow l_R(w) + l_R(w^{-1}w')  = l_R(w')$$
and 
$$w' \succcurlyeq_R w :\Leftrightarrow l_R(w'w^{-1})
+ l_R(w)  = l_R(w').$$ 
Since $R$ is invariant by conjugacy, we have
$w\preccurlyeq_R w' \Leftrightarrow w' \succcurlyeq_R w$.

Let $y\in Y^{\gen}$, let $(s_1,\dots,s_n):=\lbl(y)$, let
$(r_1,\dots,r_n):=(\pi(s_1),\dots,\pi(s_n))= \pi^n(\lbl(y))$.
By Lemma \ref{lemmaparacoxnew}, since all multiplicities are $1$,
the factorization $(r_1,\dots,r_n)$ expressed $c$ as a product of $n$ reflections.
Because the fixed-point set of $c$ is trivial (Lemma \ref{eigencox}), $c$ cannot
expressed as the product of less than $n$ reflections, so the factorization
has minimal length and it lies in $\Red_R(c)$.

The key result of this section is:

\begin{theo}
\label{theolooijenga}
Let $y\in Y^{\gen}$,
The maps
$$y\cdot B_n \stackrel{\lbl}{\longrightarrow} \lbl(y) \cdot B_n
\stackrel{\pi^n}{\longrightarrow} \pi^{n}(\lbl(y)) \cdot B_n$$
are isomorphisms of $B_n$-sets, where $y\cdot B_n$ is the
Galois orbit of $y$, and
$\lbl(y) \cdot B_n$ and $\pi^{n}(\lbl(y)) \cdot B_n$ are
Hurwitz orbits.
\end{theo}

The theorem implies that conditions (i)--(iii) from Corollary
\ref{corocardhurwitz} actually hold.

In the real case, this was initially conjectured by Looijenga, \cite[(3.5)]{looijenga}, and proved in a letter from Deligne to Looijenga
(crediting discussions with Tits and Zagier),
\cite{deligneletter};
an equivalent property (\cite[Fact 2.2.4]{dualmonoid}) was independently
used in our earlier construction of the dual braid monoid.

These proofs for the real case are based on case-by-case numerology:
because the $B_n$-sets are transitive, it suffices to prove
that the cardinality of $y\cdot B_n$ (which is by construction the 
degree of $\LL$) coincides with that of $(r_1,\dots,r_n) \cdot B_n$.
This is an enumeration problem in $W$ and may be tackled by 
case-by-case analysis (the infinite family are easy to deal with,
computers can take care of the exceptional types).

This enumerative approach carries on to our setting: Theorem
\ref{theolooijenga} immediately follows from Proposition
\ref{propcardhurwitz} below.

\begin{prop}
\label{propcardhurwitz}
\index{Hurwitz action!is transitive on $\Red_R(c)$}
Let $W$ be a well-generated complex reflection group.
Let $c$ be a Coxeter element in $W$.
The Hurwitz action is transitive on $\Red_R(c)$.
When $W$ is irreducible,
one has $|\Red_R(c)| = n!h^n/|W|$.
\end{prop}

\begin{proof}
The proposition clearly reduces to the case when $W$ is irreducible: in the reducible case,
reduced decompositions of Coxeter elements are ``shuffles'' of reduced
decompositions of the Coxeter elements of the irreducible summands.

We prove the result case-by-case (see Table 1 in section \ref{sectionwellgen}
for the list of cases to be considered).
The complexified real case is studied in \cite{deligneletter}
(transitivity is easy and does not require case-by-case, see for
example \cite{dualmonoid}, Proposition 1.6.1).
The $G(e,e,r)$ case combines two results from \cite{BC}:
Proposition 6.1 (transitivity) and Theorem 8.1 (cardinality).

The case of $G(d,1,r)$ goes as follows.
For all integers $i,j$ with $1\leq i < j \leq n$, denote by $\tau_{i,j}$
the permutation matrix associated with the transposition
$(i \; j) \in \mathfrak{S}_n = G(1,1,n) \hookrightarrow G(d,1,r)$.
For all $\zeta\in\mu_d$ and all $i\in\{1,\dots,n\}$, denote
by $\rho_{i,\zeta}$ the diagonal matrix $\Diag(1,\dots,1,\zeta,1,\dots,
1)$, where $\zeta$ is in $i$-th position.
There are two types of reflections in $G(d,1,n)$:
\emph{long} reflections are 
elements of the form
$\rho_{i,\zeta}$, with $\zeta\neq 1$; \emph{short} reflections
are elements of the form
 $\tau_{i,j}^{\zeta}:=\rho_{i,\zeta}^{-1}\tau_{i,j}\rho_{i,\zeta}$,
with $\zeta\in \mu_d$ and $1\leq i < j \leq n$.
A typical Coxeter element is
$c:=\rho_{1,\zeta_d}\tau_{1,2}\tau_{2,3}\dots\tau_{n-1,n}$.
Since Coxeter elements form a single conjugacy class, and since
$R$ is invariant under conjugacy, it suffices to prove the claims for
this particular $c$. 
Let $(r_1,\dots,r_n)\in \Red_R(c)$.
Let us prove that it is Hurwitz equivalent to
$(\rho_{1,\zeta_d},\tau_{1,2},\tau_{2,3},\dots,\tau_{n-1,n})$.
Consider
the morphism $G(d,1,n)\rightarrow G(1,1,n),g \mapsto \overline{g}$
sending a monomial matrix
to the underlying permutation matrix.
This map sends $\tau_{i,j}^{\zeta}$ to $\tau_{i,j}$ and
$\rho_{i,\zeta}$ to $1$.
The element $\overline{c}$ is a Coxeter element in $G(1,1,n)$.
One deduces that there is a unique long reflection $r_{i_0}$
among $r_1,\dots,r_n$ and that $(\overline{r_1},\dots,\widehat{\overline{r_{i_0}}},\dots,
\overline{r_n})$
is a reduced $R$-decomposition of $\overline{c}$ in $G(1,1,n)$.  
Up to applying suitable Hurwitz moves, we may
assume that $i_0=1$. Using the transitivity result already known in
the type $G(1,1,n)$ case, we see that $(r_1,\dots,r_n)$ is
Hurwitz equivalent to
$(\rho_{i,\zeta},\tau_{1,2}^{\alpha_1},\dots,\tau_{n-1,n}^{\alpha_{n-1}})$,
where $i\in\{1,\dots,n\}$,
 $\zeta\in\mu_e-\{1\}$ and $\alpha_1,\dots,\alpha_{n_1}\in\mu_d$.
By considering the determinant, we see that $\zeta=\zeta_d$.
A direct computation shows that, if $i>1$,
$$\rho_{i,\zeta}\tau_{1,2}^{\alpha_1}\dots\tau_{i-1,i}^{\alpha_{i-1}}=
\tau_{1,2}^{\alpha_1}\dots\tau_{i-1,i}^{\alpha_{i-1}}\rho_{i-1,\zeta}.$$
One may use this relation to construct an explicit
sequence of Hurwitz moves showing that
$(r_1,\dots,r_n)$ is equivalent to 
$(\rho_{1,\zeta_d},\tau_{1,2}^{\alpha'_1},\dots,\tau_{n-1,n}^{\alpha'_{n-1}})$.
One concludes by observing that
$\tau_{1,2}^{\alpha'_1}\dots\tau_{n-1,n}^{\alpha'_{n-1}}=
\tau_{1,2}\dots\tau_{n-1,n}$ forces $\alpha'_1=\alpha'_2=\dots = \alpha'_{n-1}
=1$. The claim about cardinality is not difficult,
once it is observed that an element
of $\Red_R(c)$ is uniquely determined by: 1) the position of the long
reflection $\rho_{i,\zeta}$, 2) the integer $i$, 3) a reduced $R$-decomposition
of the Coxeter element $\overline{c}$. 

There remains a finite number of exceptional
types which are treated by computer.
\end{proof}

Until the end of this section, we assume that $W$ is irreducible. 

\begin{lemma}
\label{lemmalengthsimple}
Let $y\in Y^{\gen}$, with label $(s_1,\dots,s_n)$.
Let $i\in\{1,\dots,n\}$. Then $$l_R(\pi(s_1\dots s_i))=i.$$
\end{lemma}

\begin{proof}
Since $y\in Y^{\gen}$, each $s_j$ is a braid reflection, mapped
under $\pi$ to a reflection $r_j\in R$.
Thus $l_R(\pi(s_1\dots s_i))\leq i$ and $l_R(\pi(s_{i+1}\dots s_n))\leq n-i$.
Since $\pi(s_1\dots s_i)\pi(s_{i+1}\dots s_n)$ is a Coxeter element
(Lemma \ref{pideltaiscox2}), it has
length $n$. This forces both inequalities to be equalities.
\end{proof}

\begin{lemma}
\label{corosimples}
\index{simple elements!bijection from $B$ to $W$}
The restriction of $\pi$ to the set $S$ of simple elements is
injective.
\end{lemma}

\begin{proof}
Let $s$ and $s'$ be simple elements such that $\pi(s)=\pi(s')$.
By Corollary \ref{coros1si}, we may find $y,y'\in Y^{\gen}$,
with $\lbl(y)=(s_1,\dots,s_n)$ and 
$\lbl(y')=(s_1',\dots,s_n')$,
and $i,j\in\{0,\dots,n\}$ such that $s=s_1\dots s_i$
and $s'=s'_1\dots s'_j$.
Both $\LL(y)$ and $\LL(y')$ consists of $n$ distinct points with distinct real parts;
the naive affine homotopy from $\LL(y)$ to $\LL(y')$ can be lifted to an homotopy from
$y$ to $y''$ such that $\LL(y'') = \LL(y')$ and, thanks to the Hurwitz rule, $\lbl(y) = \lbl(y'')$.
So, up to replacing $y$ with $y''$, we can assume that $\LL(y)=\LL(y')$.
By Lemma \ref{lemmalengthsimple}, we must have $i=j$.
Applying Proposition \ref{propcardhurwitz} to $w=\pi(s)$,
we may find $\beta\in B_i$ such that
 $(\pi(s_1),\dots,\pi(s_i))\cdot \beta = (\pi(s'_1),\dots,\pi(s'_i))$.
Similarly, we find $\beta'\in B_{n-i}$ such that
 $(\pi(s_{i+1}),\dots,\pi(s_n))\cdot \beta' = (\pi(s'_{i+1}),\dots,
\pi(s'_n))$. View $B_i \times B_{n-i}$ as a subgroup of $B_n$ (the
first factor braids the first $i$ strings, the second factors braids
the $n-i$ last strings), and set $\beta'':=(\beta,\beta')\in B_n$.
We have $\pi^n(\lbl(y))\cdot \beta''=\pi^n(\lbl(y'))$.
Applying Theorem \ref{theolooijenga}, this implies that $\lbl(y) \cdot \beta'' = \lbl(y')$.
Clearly, $\beta''$ does not modify the product of
the first $i$ terms of the labels.
Thus $s=s'$.
\end{proof}

\begin{coro}
\label{Sisfinite}
The set of simple elements $S$ is finite.
\end{coro}

\subsection{Simplicial Hurwitz structures}

\begin{defi}
\label{deficdc}
\index{$\CDbullet$! decompositions of Coxeter element $c$}
Let $k$ be a positive integer.
We set:
$$\CD_k(\delta) := \{ (s_1,\dots,s_k)\in S^k | \delta = s_1\dots s_k\}.$$
$$\CD_{\bullet}(\delta) := (\CD_k(\delta))_{k\in \BZ_{\geq 0}}$$
$$\CD_k(c) := \{ (w_1,\dots,w_k)\in W^k | c = w_1\dots w_k \text{ and }
l_R(c) = \sum_i l_R(w_i)\}.$$
$$\CD_{\bullet}(c) := (\CD_k(c))_{k\in \BZ_{\geq 0}}$$

\end{defi}

The definition of $\CD_{\bullet}(\delta)$ is a particular
case of Definition \ref{deficddelta}, but this anticipates
on what will be discussed in the next section (dual braid monoid
Garside structure).

As often with graded objects, it is convenient to
view $\CD_{\bullet}(\delta)$
and $\CD_{\bullet}(c)$ as disjoint unions of their graded components.

Let $t = (t_1,\dots,t_k)$ be a sequence in either
$\CD_{\bullet}(\delta)$
or $\CD_{\bullet}(c)$. We may consider:
\begin{itemize}
\item \emph{faces} of $t$, sequences of the form $$(t_1,\dots,t_{i-1},t_it_{i+1},t_{i+2},\dots,t_k)$$
\item \emph{degeneracies} of $t$, sequences of the form
$$(t_1,\dots,t_i, 1, t_{i+1},\dots,t_k).$$
\end{itemize}

This equips both $\CD_{\bullet}(\delta)$
and $\CD_{\bullet}(c)$
with a \emph{simplicial set} structure (see Appendix \ref{garsidependix}).
But there is an additional structure on both sets, provided by
``graded Hurwitz action'': each $B_k$ acts on $\CD_k$.

\begin{remark}
There are obvious compatibility rules between graded Hurwitz action
and simplicial structure. The two structures combine in a
fantastic algebraic package -- I do not have any good name
for it (``simplicial Hurwitz structure'', ``stratified Hurwitz set''?) -- 
that faithfully encodes both the monodromy theory and the ramification
theory of the Lyashko-Looijenga
covering.
By just considering the action of the ``Coxeter element'' braid in $B_k$, combined with the simplicial structure,
we obtain the ``helicoidal'' structure of Remark \ref{helicoremark} in Appendix \ref{garsidependix} (a generalization
of \emph{cyclic} structures, in the sense of Connes, \cite{connes}).
But the simplicial Hurwitz structure provides more than that (in a way,
it is a ``parabolically helicoidal structure'').
\end{remark}

\begin{theo}
\label{bulletiso}
\index{$\CDbullet$!$\CD_{\bullet}(c) \simeq \CD_{\bullet}(\delta)$}
The projection
maps $(s_1,\dots,s_k)\mapsto (\pi(s_1),\dots,\pi(s_k))$ induce
an isomorphism of simplicial sets
$$\CD_{\bullet}(\delta) \stackrel{\sim}{\longrightarrow} \CD_{\bullet}(c).$$
\end{theo}

\begin{proof}
Clearly, the map is well-defined and compatible with both faces and
degeneracies. Injectivity is an obvious consequence of Lemma \ref{corosimples}.

Surjectivity: an $(c_1,\dots,c_k)\in \CD_{\bullet}(c)$ can be obtained
by applying a sequences of face maps and degeneracy maps starting
from an element in $\Red_R(c)$ (start by concatenating reduced decompositions of non-trivial $c_i$'s, it's obvious how to get from there 
back to $(c_1,\dots,c_k)$).
By compatibility, it is enough to show that elements of $\Red_R(c)$ are
in the image. As Hurwitz action is transitive on $\Red_R(c)$ (Proposition
\ref{propcardhurwitz}) and
compatible with projection, it is enough to show that \emph{one} 
element of $\Red_R(c)$ is in the image, which is obvious: just take
a generic $y$, the projection of $\lbl(y)$ lies in $\Red_R(c)$.
\end{proof}

\begin{remark}
\label{remarktempting}
It is very tempting, and very convenient too,
 to identify $\CD_{\bullet}(\delta)$ with $\CD_{\bullet}(c)$,
 and $S$ with its image in $W$.
In particular, we will say that $w\in W$ is a \emph{simple element} if it lies in
the image of $S$.
Very often, when considering labels, we will consider those as factorisations of $c$ in $W$.
This viewpoint
helps remembering that computations involving labels are about combinatorics in a finite group,
and that small examples can be worked out by hand.
\end{remark}

\subsection{Reduced labels and trivialization of $Y$}

This following variation on the notion of label was introduced by 
Vivien Ripoll, \cite{ripoll1}, after he noticed unnecessary complications
in earlier versions of the current paper:

\begin{defi}
\label{defireducedlabel}
\index{$\rlbl$ (reduced label)}
\index{label!reduced label}
Let $y \in Y$ with label $(s_1,\dots,s_k)$. Let $(x_1,\dots,x_k)$ be
the ordered support of $\LL(y)$.
The \emph{reduced label} of $x$ is the sequence
$\rlbl(y) = (s_1',\dots,s_k')$ defined
by
$$s'_i := \left\{ \begin{matrix} s_i  & \text{ if $i=k$ or
$\re(x_i) < \re(x_{i+1})$} \\
s_i s_{i+1}^{-1} & \text{ if $i<k$ and $\re(x_i) = \re(x_{i+1})$}.
\end{matrix}
\right.$$
\end{defi}

One way to geometrically
understand reduced labels is to see them as braids represented 
by small loops around points in the support:

$$
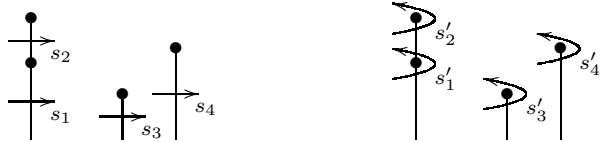
\xy
(-5,-10)="11",
(-5,0)="1", (-5,6)="2", (7,-4)="3", (7,-10)="33",
(14,2)="4", (14,-10)="44",
"1"*{\bullet},"2"*{\bullet},"3"*{\bullet},"4"*{\bullet},
"11";"1" **@{-},
"2";"1" **@{-},
"33";"3" **@{-},
"44";"4" **@{-},
(-8,-5);(-2,-5) **@{-}, *\dir{>}, (-1,-7)*{_{s_1}},
(-8,3);(-2,3) **@{-}, *\dir{>}, (-1,1)*{_{s_2}},
(4,-7);(10,-7) **@{-}, *\dir{>}, (11,-9)*{_{s_3}},
(11,-4);(17,-4) **@{-}, *\dir{>}, (18,-6)*{_{s_4}}
\endxy
\qquad
\qquad
\qquad
\xy
(-5,-10)="11",
(-5,0)="1", (-5,6)="2", (7,-4)="3", (7,-10)="33",
(14,2)="4", (14,-10)="44",
"1"*{\bullet},"2"*{\bullet},"3"*{\bullet},"4"*{\bullet},
"11";"1" **@{-},
"2";"1" **@{-},
"33";"3" **@{-},
"44";"4" **@{-},
(-8,-2);(-8,2) **\crv{(3,0)} ?>*\dir{>}, (-1,-2)*{_{s_1'}},
(-8,4);(-8,8) **\crv{(3,6)} ?>*\dir{>}, (-1,4)*{_{s_2'}},
(4,-6);(4,-2) **\crv{(15,-4)} ?>*\dir{>}, (11,-6)*{_{s_3'}},
(11,0);(11,4) **\crv{(22,2)} ?>*\dir{>}, (18,0)*{_{s_4'}}
\endxy
$$

Alternately, the reduced label of $y$ can be viewed as the label of a
generic $y'$ obtained by applying a small clockwise rotation to $y$:

\begin{lemma}
\label{lemmareducedlabel}
For all $y\in Y$, there exists a real number $\alpha >0$, such that
for all $\varepsilon$ such that $0 < \varepsilon < \alpha$,
the reduced label of $y$ coincides with the label of $e^{-\sqrt{-1}\pi \varepsilon}y$.

As a consequence, the reduced label of $y$ is a sequence of simple elements, with 
product $\delta$.
\end{lemma}

Note: this result is better understood in light of the general study of the $\BC^{\times}$-action,
see Section \ref{section11} and, in particular, the basic Lemma \ref{lemmaLLd}.

\begin{proof}
Obvious consequence of the Hurwitz rule.
\end{proof}

Rephrasing Corollary \ref{corodeeplabel} in terms of reduced labels, we also get:

\begin{lemma}
\label{lemmarlblcd}
For all $y\in Y$, $\rlbl(y)\in\CD_{\bullet}(c)$.
\end{lemma}

(As announced in Remark \ref{remarktempting}, we choose to work in $\CD_{\bullet}(c)$
rather than $\CD_{\bullet}(\delta)$.)

\begin{defi}
\label{deficompatible}
Let $x \in E_n$ with ordered support $(x_1,\dots,x_k)$
and ordered multiplicity $(n_1,\dots,n_k)$.
Let $\sigma=(s_1,\dots,s_l) \in \CD_{\bullet}(c)$ be a factorization
of $c$ into $l$ simple elements.

Both $(n_1,\dots,n_k)$ and $(l_R(s_1),\dots,l_R(s_k))$ are
\emph{compositions of $n$} (finite integral sequences that add up to $n$).
If $(n_1,\dots,n_k) = (l_R(s_1),\dots,l_R(s_k))$, we say that
$x$ and $\sigma$ are \emph{compatible}.

We denote by $E_n \boxtimes \CD_{\bullet}(c)$ 
the set of compatible pairs. In other words, it is the
pullback:

$$
\xymatrix{
E_n \boxtimes \CD_{\bullet}(c)
 \ar[r] \ar[d] &  \CD_{\bullet}(c) \ar[d]\\
E_n \ar[r] & \{\text{compositions of $n$}\}
}
$$
\end{defi}

Clearly, if $x$ and $\sigma$ are compatible,
then $\sigma$ must be non-degenerate.
So the pullback map $E_n \boxtimes \CD_{\bullet}(c) \to\CD_{\bullet}(c)$ is not surjective.
People finding this annoying may want to introduce a specific notation for the set
of non-degenerate factorizations (factorizations not containing the trivial element).

\begin{defi}
\label{defiface}
\index{$\sigma \vdash \tau$ ($\tau$ is a face of $\sigma$)}
\index{$F_{\sigma}$ (reduced decompositions having $\sigma$ as a face)}
If $\sigma,\tau\in \CD_{\bullet}(c)$,
we say that \emph{$\tau$  is a face of $\sigma$}, and
write $$\sigma \vdash \tau$$
if each term in $\tau$ is the partial product of
consecutive terms in $\sigma$ (in other words: $\tau$ is obtained
from $\sigma$ by consecutive simplicial face operators).

For all $\sigma\in \CD_{\bullet}(c)$, we set $$F_{\sigma}:=\{\rho\in \Red_R(c)|
\rho\vdash \sigma\}.$$
\end{defi}

\begin{lemma}
\label{lemmanotempty}
If $\sigma=(c_1,\dots,c_k)\in \CD_{\bullet}(c)$ is non-degenerate, then $F_{\sigma}\neq\varnothing$, and consists of a single Hurwitz orbit for the natural subgroup
$B_{l(c_1)}\times \dots \times B_{l(c_k)}$ of the braid group $B_n$.
\end{lemma}

\begin{proof}
Elements of $F_{\sigma}$ are obtained by concatenating
reduced decompositions of $c_1,\dots,c_k$. Each $c_i$ is a parabolic
Coxeter element (Lemma \ref{lemmaparacoxnew}),
whose factorizations form a single $B_{l(c_i)}$-orbit (Proposition \ref{propcardhurwitz}.)
\end{proof}

\begin{theo}[trivialization of $Y$]
\label{theolabel}
\index{Saito quotient!decomposition theorem}
The map $\LL\times \rlbl$ induces a bijection
\begin{eqnarray*}
\LL\times \rlbl:Y & \stackrel{\sim}{\longrightarrow} & E_n \boxtimes \CD_{\bullet}(c)
\end{eqnarray*}
\end{theo}

\begin{proof}
That the image of $\rlbl$ lies in $\CD_{\bullet}(c)$ is Lemma \ref{lemmarlblcd}.
That the image of $\LL\times \rlbl$ lies in $E_n \boxtimes \CD_{\bullet}(c)$
follows from Corollary \ref{corolengthlabel}.
By Theorem \ref{theolooijenga}, $\LL\times \rlbl$ restricts to a bijection
$Y^{\reg} \stackrel{\sim}{\rightarrow} E_n^{\reg} \boxtimes \CD_{\bullet}(c)$.
What remains at stake is the behavior in the singular part of the covering.

Surjectivity of $LL\times \rlbl$.
Let $(e,\sigma)\in E_n \boxtimes \CD_{\bullet}(c)$.
Let $(e_m)_{m\in \BZ_{\geq 0}}$ be a sequence of points in $E_n^{\gen}$ converging
to $e$ (because $E_n^{\gen}$ is dense in $E_n$, such a sequence exists).
Let $\sigma'\in F_{\sigma}$ (which is non-empty by Lemma \ref{lemmanotempty}).
By Theorem \ref{theolooijenga}, there exists a unique sequence
$(y_m(\sigma'))_{m\in \BZ_{\geq 0}}$
of points in $Y$ such that, for all $m$, $LL(y_m)=e_m$ and $\rlbl(y_m(\sigma'))=\sigma'$.
This sequence lies in $LL^{-1}(\{e\}\cup \bigcup_{m}\{e_m\})$, a compact subset of $Y$
(a finite morphism, $LL$ is proper; the pre-image of a compact subset under a proper morphism
is compact),
and admits an adherence value $y$ such that $LL(y)=e$. Applying the Hurwitz rule, one observes
that $\rlbl(y) = \sigma$.

Note that the above argument also shows that the ramification degree of $LL$ at $y$
is $|F_{\sigma}|$, as each $LL^{-1}(e_m)$ consists of the $|F_{\sigma}|$ distinct points
$\{y_m(\sigma') | \sigma' \in F_{\sigma}\}$.
 
Injectivity of $LL\times \rlbl$.
Let $e\in E_n$.
As $LL$ is a finite morphism of degree $|\Red_R(c)|$, the ramification formula over $e$ 
says that
$$|\Red_R(c)| = \sum_{y \in LL^{-1}(e)} d(y),$$
where $d(y)$ is the ramification degree at $y$ (see for example \cite[Example 4.3.7]{fulton}).
By grouping according to $\rlbl(y)$, we get
$$|\Red_R(c)| = \sum_{
\substack{\sigma \in \CD_{\bullet}(c) \\ \text{$\sigma$ compatible with $e$}}
}
|LL^{-1}(e) \cap \rlbl^{-1}(\sigma)|\cdot |F_{\sigma}|.$$
Clearly, each element of $\Red_R(c)$ lies in $F_{\sigma}$ for exactly one $\sigma\in\CD_{\bullet}(c)$
such that $\sigma$ is compatible with $e$ (this $\sigma$ is obtained by
multiplying consecutive terms according to the multiplicity pattern defined by $e$).
So $$|\Red_R(c)| = \sum_{
\substack{\sigma \in \CD_{\bullet}(c) \\ \text{$\sigma$ compatible with $e$}}
} |F_{\sigma}|.$$
To conclude, we observe that the surjectivity part of theorem implies that, for each $\sigma$ compatible with $e$, $|LL^{-1}(e) \cap \rlbl^{-1}(\sigma)| \geq 1$.
This forces each $|LL^{-1}(e) \cap \rlbl^{-1}(\sigma)|$ to be equal to $1$.
\end{proof}

We can equip $E_n \boxtimes \CD_{\bullet}(c)$ with a natural topology that
turns the bijection of Theorem \ref{theolabel} into an homeomorphism. This topology is
as follows. Let $(x, (c_1,\dots,c_m))\in E_n \boxtimes \CD_{\bullet}(c)$.
Choose a system of elementary tunnels $(T_1,\dots,T_m)$
(as introduced above Definition \ref{defilabel}). Neighborhoods for $(x, (c_1,\dots,c_m))$
is obtained by considering all compatible pairs $(x', (d_1,\dots,d_l))$ as follows:
\begin{itemize}
\item[(1)] we choose $\Omega$, a small enough neighborhood of $x$ in $E_n$ such that (the projections in $\BC$ of) the tunnels
$T_1,\dots,T_m$ do not intersect any point in any configuration $x''\in \Omega$; that
such neighborhoods exist and form a basis for the topology of $E_n$ is obvious (and
is the analog in $E_n$ of the notion of $T$-neighborhood from Definition \ref{defiTn}),
\item[(2)] we allow $x'$ to be any point in $\Omega$,
\item[(3)] combining (1) and (2), we get Hurwitz rule equations for $T_1,\dots,T_m$, expressing
relations between the $d_j$'s to the $c_i$'s that must be satisfied.
\end{itemize}

\begin{remark}
\label{remarkhomotopylifting}
In the previous section, we used the homotopy lifting property of the unramified
part of $LL$ to lift paths in $E_n^{\reg}$ to paths in $Y$. In general, if $\gamma$ is a path
$[0,1]\to E_n$, there may be more than one way to lift $\gamma$ to a path
$\tilde{\gamma}$ such that $LL\circ \tilde{\gamma} = \gamma$, even if one fixes
the initial point $\tilde{\gamma}(0)$. However, a consequence of Theorem
\ref{theolabel} and the above discussion
is that if $\gamma$ has \emph{non-decreasing ramification} (\ie,
if points can be merged but not unmerged when $t$ increases), then there exists
a unique continuous lift $\tilde{\gamma}$ once $\tilde{\gamma}(0)$ has been fixed.
This will be very useful for constructing explicit retractions in the next sections.
\end{remark}

\begin{remark}
\label{halfwaysummary}
Charting the various geometric constructions, we get: 

$$
\xymatrix{
\bigcup_{H\in \ACA} H \ar[r] \ar[d] & V  \ar[d]_{W \cq}\\
\CH \ar[r] & W\cq V \ar@{=}[r] & Y \times \BC \ar[d]^{\text{Saito}}
  & &   \\
& \CK \ar[d] \ar[r] & Y \ar[d]_{\LL} \ar[rr]^{\rlbl} \ar@{=}[drr]^{\ref{theolabel}}
 & & \CD_{\bullet}(c) \\
& \CH_{A_{n-1}} \ar[r] & E_n & & \ar[ll] \ar[u] E_n \boxtimes \CD_{\bullet}(c)\\
& \bigcup_{1\leq i < j \leq n} H_{i,j} \ar[r] \ar[u]
 & \BC^n \ar[u]^{\mathfrak{S}_n  \cq}
}
$$

Here is what happened so far. We set out to study the higher homotopy groups
of $V^{\reg} = V - \bigcup_{H\in \ACA} H$.
Because $W$ has no ramification on $V^{\reg}$, we may work in the quotient
$W\cq V^{\reg}$, which we view as a singular fibration
over $Y$. Individually, the generic fibers (outside $\CK$)
and degenerated fibers (above $\CK$) are fairly
easy to control: they all are punctured complex lines.

Most of the hard work happens in the base space $Y$, which controls
how fibers are glued together (degeneracy, monodromy, etc.).
To visualize $Y$ and perform computations in it, we compare it via
$\LL$ with a classical configuration
space $E_n$. The trivialization theorem \ref{theolabel} gives a neat
description, in terms of the combinatorics of $W$, of both generic and
singular fibers of $\LL$.

One way to see $\CD_{\bullet}(c)$ is think of it as the ``Galois
group'' of $\LL$ -- except that $\LL$ is \emph{not} a Galois covering.
Metaphorically, $\LL: Y\to E_n$ is a \emph{virtual reflection group}: like
the quotient map $V \to W\cq V$, it is a finite
algebraic morphism between two affine spaces; when such a map is Galois,
the theorem of Chevalley-Shephard-Todd says that it must be the quotient map of a
complex reflection group.

\end{remark}

\subsection{A variation: trivializing $W\cq V$ and $W\cq V^{\reg}$}

The way from $W\cq V^{\reg}$
to $E_n\boxtimes \CD_{\bullet}(c)$, as summarized in Remark \ref{halfwaysummary}, is a bit long and complicated. This can be simplified
thanks to variations on the definitions of $LL$ and $E_n$:

\begin{defi}[extended Lyashko-Looijenga morphism]
\index{$\LLL$ (extended Lyashko-Looijenga morphism)}
Let $(y,z)$ be a point in $W\cq V\simeq Y\times \BC$.
We denote by $\LLL((y,z))$ the configuration $LL(y) - z$, obtained by shifting
by $-z$ all points in $LL(y)$.
\end{defi}

Note that the image of $\LLL$ lies in $\overline{E}_n$ (the space of not necessarily centered
configurations) rather than just $E_n$ (the space of centered configurations).

\begin{defi}
We denote by $\overline{E}_{n}^{\circ}$ the subspace of $\overline{E}_n$
\index{$\En$! (configurations not containing $0$) $\overline{E}_{n}^{\circ}$} 
consisting of configurations not containing $0$.
\end{defi}

For $(y,z)\in W\cq V$, we also set $\rlbl((y,z)):= \rlbl(y)$. The notion of
\emph{compatible pairs} in $E_n\times \CD_{\bullet}(c)$ carries on to 
$\overline{E}_n\times \CD_{\bullet}(c)$, and we define
$\overline{E}_n\boxtimes \CD_{\bullet}(c)$.

\begin{theo}
\label{nicetriv}
The map
\begin{eqnarray*}
\LLL \times \rlbl: W\cq V & \longrightarrow & \overline{E}_n\times \CD_{\bullet}(c) \\
(y,z) & \longmapsto & (LL(y) - z, \rlbl(y))
\end{eqnarray*}
is an homeomorphism.
It induces by restriction an homeomorphism:
$$\LLL \times \rlbl: W\cq V^{\reg} \stackrel{\sim}{\longrightarrow} \overline{E}_{n}^{\circ}
\boxtimes \CD_{\bullet}(c).$$
\end{theo}

\begin{proof}
Consider the map
\begin{eqnarray*}
W\cq V & \longrightarrow & \BC\times (E_n \boxtimes \CD_{\bullet}(c))\\
(y,z) & \longmapsto & (z, LL(y), \rlbl(y))
\end{eqnarray*}
Using Theorem \ref{theolabel} and the subsequent discussion, we see that it is an
homeomorphism.
Now observe that
\begin{eqnarray*}
\BC\times E_n & \longrightarrow & \overline{E} \\
(z,\pmb{\{ }x_1,\dots,x_n\pmb{\} }) & \longmapsto & \pmb{\{ }x_1 - z,\dots,x_n - z\pmb{\} }
\end{eqnarray*}
is an homeomorphism: indeed, $z$ can be recovered from $\pmb{\{ }x_1 - z,\dots,x_n - z\pmb{\}}$
as its barycenter.
The theorem follows easily.
\end{proof}

The next 3 sections rely on the trivialization of Theorem \ref{theolabel}, but
it is possible to rephrase them using Theorem \ref{nicetriv} instead.
Theorem \ref{nicetriv} will be especially useful in Section \ref{section11}, where we
will focus on rotational motions (rather than the translational motions used in the next
3 sections).

\section{The dual braid monoid}

Here again, $W$ is an irreducible well-generated complex reflection
groups.

For simplicity, we further assume, in this section and in the following ones,
that $W$ is generated by $2$-reflections:
by Theorem \ref{shephard}, this suffices to address the $K(\pi,1)$ conjecture.
We restrict ourselves to $2$-reflection groups not because the construction
would otherwise fail (it does work),
but because some case-by-case arguments (especially Lemma \ref{lemmalattice})
would need to be more extensively detailed.

We keep the notations from the previous sections.
Recall that an element $b\in B$ is \emph{simple} if $b=b_T$ for some tunnel
$T$, and that the set of simple elements is denoted by $S$.

\begin{defi}
\index{$M$ (dual braid monoid)}
\index{dual braid monoid}
The \emph{dual braid monoid} is the submonoid $M$ of $B$ generated by $S$.
\end{defi}

\index{$\preccurlyeq$, $\prec$ (left divisibility in $S$ and $M$)}
Consider the binary relation $\preccurlyeq$ defined on $S$ as follows.
Let $s$ and $s'$ be simple elements. We write $s \preccurlyeq s'$ if and
only if there exists $(y,z)\in W\cq V^{\reg}$, $L,L'\in \BR_{\geq 0}$ with
$L \leq L'$ such that $(y,z,L)$ is a tunnel representing $s$ and
$(y,z,L')$ is a tunnel representing $s'$.

We write $s \prec s'$ when $s\preccurlyeq s'$ and $s\neq s'$.

This section is devoted to the proof of:

\begin{theo}
\label{theogarside}
\index{dual braid monoid!is Garside monoid}
The monoid $M$ is a Garside monoid, with set of simple elements $S$ and Garside
element $\delta$. The relation $\preccurlyeq$ defined above on $S$ is the
restriction to $S$ of the left divisibility order in $M$.

The monoid $M$ generates $B$, which inherits a structure of Garside group.
\end{theo}

A survey of Garside theory is provided in Appendix \ref{garsidependix}.

\begin{remark}
The theorem blends two results of distinct natures: one is about
the Garside structure of a certain monoid $\BM(P_c)$ (see
Lemma \ref{lemmadualbraidrelations} for a presentation of
$\BM(P_c)$); the other one identifies $M$ with $\BM(P_c)$. The latter
essentially amounts to writing a presentation for $B$. It
is a substitute for Brieskorn's presentation theorem, \cite{brieskorn1},
except that our presentation involves
\emph{dual braid relations} instead
of Artin-Tits braid relations.
\end{remark}

Several results mentioned in the introduction follow from this theorem:

Theorem \ref{theointrogarside} does not assume irreducibility,
but follows immediately from the irreducible case since
direct products of Garside groups are Garside groups.

Another important consequence of Theorem \ref{theogarside} is that
one obtains a nice simplicial complex, $\gar(G,\Sigma)$ (see Definition \ref{gar2}), that is both contractible (Theorem \ref{gar3}) and acted on
by $B$. This complex will serve as a simplicial model for the universal cover of $W \cq V^{\reg}$.

The strategy of proof of Theorem \ref{theogarside} is
very similar to the one in \cite{BC}.

\begin{prop}
\label{propiiiii}
For all $s \in S$, we have $l(s)=l_R(\pi(s))$. For all
 $s,s'\in S$, the following statements are equivalent:
\begin{itemize}
\item[(i)] $s \preccurlyeq s'$,
\item[(ii)] $\exists s''\in S, ss''=s'$,
\item[(iii)] $\pi(s) \preccurlyeq_R \pi(s')$.
\end{itemize}
\end{prop}

\begin{proof}
The first statement follows from Lemma \ref{lemmalengthsimple}.

(i) $\Rightarrow$ (ii). 
Assume that $s \preccurlyeq s'$ and choose $(y,z)\in W\cq V^{\reg}$ 
and $L \leq L'$ such that $(y,z,L)$ is a tunnel representing $s$ and
$(y,z,L')$ is a tunnel representing $s'$.
Then $(y,z+L,L'-L)$ is a tunnel representing $s''\in S$ such
that $ss''=s'$.

(ii) $\Rightarrow$ (iii). The natural length function $l$ is additive
on $B$. Thus, under (ii), we have
 $l(s)+l(s'')=l(s')$. On the other hand, for all $\sigma\in S$,
$l(\sigma)=l_R(\pi(\sigma))$. Thus $l_R(\pi(s))+l_R(\pi(s''))=l_R(\pi(s'))$.
Since $\pi(s'') = \pi(s)^{-1}\pi(s')$, this implies that $\pi(s)\preccurlyeq_R
\pi(s')$.

(iii) $\Rightarrow$ (i).
We may find $y,y'\in Y^{\gen}$, with $\lbl(y)=(s_1,\dots, s_n)$ and
$\lbl(y')=(s'_1, \dots,s'_n)$ such that
$s=s_1s_2\dots s_{l(s)}$ and $s'=s'_1s'_2\dots s'_{l(s')}$.
Set $w:=\pi(s)$, $w':=\pi(s')$ and, for $i=1,\dots,n$,
set $r_i:=\pi(s_i)$ and $r'_i:=\pi(s'_i)$.
Assuming (iii), we may find $r''_{1},\dots,r''_{l(s')-l(s)}\in R$ such
that $r_1 \dots r_{l(s)} r''_{1}\dots r''_{l(s')-l(s)} = 
r'_1\dots r'_{l(s')}$.
The sequences $$(r_1, \dots, r_{l(s)}, r''_{1},\dots, r''_{l(s')-l(s)})$$
and $$(r'_1,\dots ,r'_{l(s')})$$ both lie in $\Red_R(w')$.
Since $w'\preccurlyeq c$, both sequences lie in the same Hurwitz orbit
(Proposition \ref{propcardhurwitz}). Thus
$$(r'_1,\dots ,r'_{n})$$ and
$$(r_1, \dots, r_{l(s)}, r''_{1},\dots,
 r''_{l(s')-l(s)},r'_{l(s')+1},\dots,
r'_n)$$ are transformed one onto the other by Hurwitz action of
a braid $\beta\in B_n$ only braiding the first $l(s')$ strands.
The Hurwitz transform of
$$\lbl(y') = (s'_1,\dots,s'_n)$$
by $\beta$ is the label 
$$(s''_1,\dots,s''_n)$$ of some $y''\in Y^{\gen}$.
Since the braid only involves the first $l(s')$ strands,
$s''_1\dots s''_{l(s')}= s'_1\dots s'_{l(s')} = s'$.
One has $\pi(s''_i)= r_i$ for $i=1,\dots,l(s)$,
thus $\pi(s''_1\dots s''_{l(s)}) = \pi(s_1\dots s_{l(s)})$.
By Lemma \ref{corosimples}, this implies
$s''_1\dots s''_{l(s)} = s_1\dots s_{l(s)} =s$.
For $x$ with small enough real and imaginary parts, one may find real
numbers $L,L'$ with $0< L \leq L'$ such that $(y'',x,L)$ is a tunnel with
associated simple $s''_1\dots s''_{l(s)}=s$ and $(y'',x,L')$ is a tunnel with
associated simple $s''_1\dots s''_{l(s')} =s'$.
\end{proof}

\begin{prop}
\label{posetiso}
The map $\pi$ restricts to an isomorphism
$(S,\preccurlyeq)\stackrel{\sim}{\rightarrow}
([1,c],\preccurlyeq_R)$. In particular, $\preccurlyeq$ is an order
relation on $S$.
\end{prop}
 
\begin{proof}
The previous proposition, applied to $s'=\delta$, proves that
$\pi(S)\subseteq [1,c]$. It also proves that $\pi$ induces a morphism
of sets with binary relations $(S,\preccurlyeq){\rightarrow}
([1,c],\preccurlyeq_R)$.
The injectivity is Lemma \ref{corosimples}.

Surjectivity: Choose $y\in Y^{\gen}$.
Let $(s_1,\dots,s_n):=\lbl(y)$. Let $r_i:=\pi(s_i)$.
We have $(r_1,\dots,r_n) \in \Red_R(c)$
Let $w \in [1,c]$. We may find $(r'_1,\dots,r'_n) \in \Red_R(c)$ such
that $r'_1\dots r'_{l_T(w)} = w$.
By Proposition \ref{propcardhurwitz}, $(r'_1,\dots,r'_n)$ is a Hurwitz
transformed of $(r_1,\dots,r_n)$, thus there exists $y'\in Y^{\gen}$
such that $\pi_*(\lbl(y'))=(r'_1,\dots,r'_n)$. The simple element which
is the product of the first $l(w)$ terms of $\lbl(y')$ is in $\pi^{-1}(w)$.
\end{proof}

We have the key lemma:

\begin{lemma}
\label{lemmalattice}
The poset $([1,c],\preccurlyeq_R)$ is a lattice.
\end{lemma}

\begin{proof}
The real case is done in \cite{dualmonoid} -- a recent beautiful
case-free argument has been found by Brady-Watt, \cite{bw2}.
The case of $G(e,e,r)$ is done
in \cite{BC}. The remaining cases have been checked by computer.
\end{proof}

\begin{defi}
The \emph{dual braid relations} with respect to $W$ and $c$ are all the
formal relations of the form $$rr'=r'r'',$$
where $r,r',r''\in R$ are such that $r\neq r'$, $rr'\in [1,c]$, and the
relation $rr'=r'r''$ holds in $W$.
\end{defi}

Clearly, dual braid relations only involve reflections in $R\cap [1,c]$.
When $W$ is complexified real, $R\subseteq [1,c]$. This does not hold
in general (see the tables at the end of the article).

As in \cite{dualmonoid} (see also \cite[Section 1]{cyclic} and Appendix \ref{garsidependix}), one
endows $[1,c]$ with a partial product and obtains a monoid
$\BM(P_c)$.

From Lemma \ref{lemmalattice}, we will deduce
that $\BM(P_c)$ is a Garside monoid. We will also identify
$\BM(P_c)$ with $M$. The following lemma generalizes
\cite[Theorem 2.1.4]{dualmonoid}:

\begin{lemma}
\label{lemmadualbraidrelations}
Let $R_c:=R\cap [1,c]$.
The monoid $\mathbf{M}(P_c)$ admits the monoid presentation 
$$\mathbf{M}(P_c) \simeq \left< R_c \left|
\text{ dual braid relations} \right. \right>.$$
\end{lemma}

\begin{remark}
When viewed as a group presentation, the presentation of the lemma is 
a presentation for $\mathbf{G}(P_c)$. As soon as we prove
 Theorem \ref{theogarside}, Lemma \ref{lemmadualbraidrelations}
will give an explicit presentation for $B$. A way to reprove
Theorem \ref{theointropres} for groups different from
$G_{31}$ is by simplifying the (redundant) presentation
given by the lemma. This does not involve any computer-assisted monodromy computation.
\end{remark}

\begin{proof}[Proof of Lemma \ref{lemmadualbraidrelations}]
By definition, $\mathbf{M}(P_c)$ admits the presentation with
generators $R_c$ and a relation
$r_1\dots r_k= r'_1\dots r'_k$ for each pair $(r_1,\dots,r_k)$,
$(r'_1,\dots,r'_k)$ of reduced $R$-decompositions of the same element
$w\in [1,c]$.
Call these relations \emph{Hurwitz relations}.
By transitivity of the Hurwitz action on $\Red_R(c)$ (Proposition
\ref{propcardhurwitz}), the Hurwitz relations are consequences of the 
dual braid relations.
The dual braid relations clearly hold in $\mathbf{M}(P_c)$
(to see this, complete $(r,r')$ to an element
$(r,r',r_3,\dots,r_n)\in\Red_R(c)$). This proves the lemma.
\end{proof}

\begin{proof}[Proof of Theorem \ref{theogarside}]
Set $R_c := R \cap [1,c]$.
We are in the situation of subsection 0.4 in \cite{dualmonoid},
$(W,R_c)$ is a \emph{generated group} and $c$ is \emph{balanced}
(one first observes that $\{w | w \preccurlyeq_{R_c} c \} = 
\{w | w \preccurlyeq_{R} c \}$ and
$\{w | c \succcurlyeq_{R_c} w \} = 
\{w | c \succcurlyeq_{R} w \}$; one concludes noting that $c$ is balanced
with respect to $(W,R)$, which is immediate since $R$ is invariant by
conjugacy).
By Lemma \ref{lemmalattice} and \cite[Theorem 0.5.2]{dualmonoid},
the premonoid $P_c:=[1,c]$ (together with the natural partial product)
is a Garside premonoid.
We obtain a Garside monoid $\mathbf{M}(P_c)$.

By Proposition \ref{posetiso}, the restriction of $\pi$ is a bijection
from $S$ to $P_c$. Let $\phi$ be the inverse bijection.
Let $w,w'\in P_c$.
Assume that the product $w.w'$ is defined in $P_c$. Let $w''$ be the value
of this product. One has $w \preccurlyeq_R w''$, thus
$\phi(w) \preccurlyeq \phi(w'')$ (using again Proposition \ref{posetiso}),
and we may find $b'$ in $S$ such that $\phi(w) b' = \phi(w'')$ (Proposition
\ref{propiiiii}).
Claim: $b'= \phi(w')$. Indeed, $b'$ and $\phi(w')$ are two simple elements
whose image by $\pi$ is $w^{-1}w''$; one concludes using Lemma
\ref{corosimples}.
This proves that $\phi$ induces a premonoid morphism $P_c \rightarrow S$
(where $S$ is equipped with the restriction of the
monoid structure), thus induces a monoid
morphism $\mathbf{M}(P_c) \rightarrow M$ and 
a group morphism $\Phi:\mathbf{G}(P_c) \rightarrow B$.

Let use prove that $\Phi$ is an isomorphism.
Choose a basepoint $y\in Y^{\gen}$.
Let $\gamma_1,\dots,\gamma_n$ be generators of $\pi_1(Y-\CK,y)$.
Let $(s_1,\dots,s_n)$ be the label of $y$.
Let us reinterpret the presentation from
Theorem \ref{theovankampen}  in terms of Hurwitz action.
Since $\LL:  Y-\CK \rightarrow E_n^{\reg}$ is a covering,
$\pi_1(Y-\CK,y)$ may be identified with a subgroup $H\subseteq B_n$.
The generators $f_1,\dots,f_n$ in Theorem \ref{theovankampen} may
be chosen to be $s_1,\dots,s_n$, and the monodromy automorphism
$\phi_1,\dots,\phi_m$ are obtained by Hurwitz action on $s_1,\dots,s_n$.
Let $h\in H$. Let $(s'_1,\dots,s'_n):= h(s_1,\dots,s_n)$; by this,
we mean the Hurwitz action of $h$ on the free group generated by
$s_1,\dots,s_n$; the $s_i'$'s are words in the $s_i$'s.
Call \emph{Van Kampen relations} the relations of the
type $s'_i = s_i$, for any $i\in \{1,\dots,n\}$, $h\in H$, and $s'_i$
obtained as above.
We have
$$B\simeq \left< s_1,\dots,s_n \left|
\text{ Van Kampen relations} \right. \right>.$$
The map $\pi$ induces a bijection from 
 $A:=\pi^{-1}(R_c) \subseteq S$ to $R_c$.
Let $r_1,\dots,r_n$ be the images of $s_1,\dots,s_n$.
By transitivity of Hurwitz action on $\Red_R(c)$, the group
$\mathbf{G}(P_c)$ is generated by $r_1,\dots,r_n$, the remaining
generators in the presentation of Lemma \ref{lemmadualbraidrelations}
appearing as conjugates of $r_1,\dots,r_n$ (by successive use of dual
braid relations).
Our generating sets are compatible and the morphism
$$\Phi: 
\mathbf{G}(P_c) \simeq \left< R_c \left|
\text{ dual braid relations} \right. \right> \rightarrow B \simeq
\left< s_1,\dots,s_n \left|
\text{ Van Kampen relations} \right. \right>$$
is defined by $r_i\mapsto s_i$.
Add to the presentation of $B$
 formal generators indexed by $$A -\{s_1,\dots,s_n\}$$
as well as the dual braid relations $\pi(r)\pi(r')=\pi(r')\pi(r'')$.
Since the relations already hold in $\mathbf{G}(P_c)$, they hold in $B$,
and we obtain a new presentation
$$B\simeq \left<A  \left|
\text{ Van Kampen relations on } \{s_1,\dots,s_n\}, \text{ dual braid
relations on } A \right. \right>.$$
To conclude that $\Phi$ is an isomorphism, it is enough to observe
that the dual braid relations encode the full Hurwitz action of
$B_n$ of $(s_1,\dots,s_n)$, while the Van Kampen relations encode the
action of $H\subseteq B_n$: thus Van Kampen relations are consequences
of dual braid relations, and $\mathbf{G}(P_c)$ and $B$ are given by
equivalent presentations.

Since $\Phi$ is an isomorphism and  $\mathbf{M}(P_c)$ naturally
embeds in $\mathbf{G}(P_c)$ (this is a crucial property of Garside
monoids, \cite{dehgar}), $\mathbf{M}(P_c)$ is isomorphic to its
image $M$ in $B$. The rest of theorem is clear.
\end{proof}

As mentioned earlier, one may view $B$ as a ``reflection group'', generated
by the set $\mathcal{R}$ of all braid reflections.
An element of $M$ is in $\mathcal{R}$ if and only if it has length $1$ for
the natural length function, or equivalently if it is an \emph{atom}
(\ie, an element which has no strict divisor in $M$ except the unit).
By Proposition \ref{propcardhurwitz} and Theorem \ref{theolooijenga},
there is a bijection between $\Red_R(c)$ and the image of $\lbl:Y^{\gen}
\rightarrow \mathcal{R}^n$. This image is clearly contained
in $\Red_{\mathcal{R}}(\delta)$. The conjecture below claims an analogue
in $B$ of the transitivity of the Hurwitz action on $\Red_R(c)$
 ($\delta$ is the natural substitute for a Coxeter
element in $B$).
It implies that any
element in  $\Red_{\mathcal{R}}(\delta)$ is the label of some $y\in Y^{\gen}$.

\begin{conj}
\label{braidreflections}
\index{Hurwitz action!conjectural transitivity on braid factorizations}
The Hurwitz action of $B_n$ on $\Red_{\mathcal{R}}(\delta)$ is transitive.
\end{conj}

\section{Chains of simple elements}
\label{sectionchains}

Here again, $W$ is an irreducible well-generated complex reflection group,
and the notations from the previous section are still in use.
For $k = 0,\dots, n$, we denote by $\CC_k$ the set of (strict)
 chains in $S-\{1\}$ of
cardinal $k$, \ie, the set of $k$-tuples $(c_1,\dots,c_k)$ in $S^k$ such that
$$1 \prec c_1 \prec \dots \prec c_k,$$
or equivalently the set of $k$-tuples $(c_1,\dots,c_k)$ in $M^k$ such that
$$1 \prec c_1 \prec \dots \prec c_k \preccurlyeq \delta.$$
It is convenient to write $\{1 \prec c_1 \prec \dots \prec c_k \}$ instead
of $(c_1,\dots,c_k)$.

We set $\CC:=\bigsqcup_{k=0}^n \CC_k$.

Let $C:=\{1 \prec c_1 \prec \dots \prec c_k\} \in \CC$.
We say that $y\in Y$ \emph{represents}
$C$ if there exists $x\in U_y$ and real numbers $L_1,\dots,L_k$ such that
$0< L_1 < \dots < L_k$ and, for $i=1,\dots,k$, $(y,x,L_i)$ is
a tunnel representing $C_i$.

Example: if $\LL(y)$ is as in the illustration below and
$(s_1,\dots,s_5) = \lbl(y)$, $y$ represents
 $1\prec s_1\prec s_1s_3 \prec
s_1s_3 s_4\prec s_1s_3s_4s_5 $,
 $1\prec s_2 \prec s_2 s_4 \prec s_2s_4s_5$,
$1\prec s_2 \prec s_2s_5$
and their subchains,
but does not represent $1\prec s_2 \prec s_2s_3$
nor $1\prec s_3 \prec s_3s_5 $.
$$
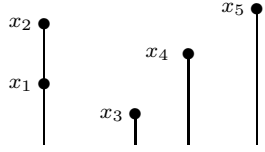
\xy
(-5,-10)="11",
(-5,-2)="1", (-5,6)="2", (7,-6)="3", (7,-10)="33",
(14,2)="4", (14,-10)="44",
(23,8)="5",(23,-10)="55",
"1"*{\bullet},"2"*{\bullet},"3"*{\bullet},"4"*{\bullet},
"5"*{\bullet},(20,8)*{_{x_5}}, "55";"5" **@{-},
(-8,-2)*{_{x_1}},(-8,6)*{_{x_2}},(4,-6)*{_{x_3}},(10,2)*{_{x_4}},
"11";"1" **@{-},
"2";"1" **@{-},
"33";"3" **@{-},
"44";"4" **@{-},
\endxy $$

\begin{defi}
For all $C\in \CC$, we set $Y_C := \{y\in Y | y \text{ represents }
C\}$.
\end{defi}

To illustrate this notion, we observe that implication (ii) $\Rightarrow$
(i) from Proposition \ref{propiiiii} expresses that for all $C\in \CC_2$,
$Y_C$ is non-empty. Based on the results from the
previous sections, this easily generalizes to:

\begin{lemma}
For all $C\in \CC$, the space $Y_C$ is non-empty.
\end{lemma}

The goal of this section is to prove:

\begin{prop}
\label{YCcontractible}
For all $C\in \CC$, the space
$Y_C$ is contractible.
\end{prop}

This technical result will be used in Section \ref{sectionuniversalcover},
when studying the nerve of an open covering of the universal cover of
$W\cq V^{\reg}$: we will need to prove that certain non-empty intersections of
open sets are contractible, and these intersections will appear as
fiber bundles over some $Y_C$, with contractible fibers.

The proposition
is not very deep nor difficult but
somehow \emph{inconvenient} to prove since the retraction will be described
via $\LL$, through ramification points. The following particular cases
are easier to obtain:
\begin{itemize}
\item If $C$ is the
chain $1 \prec \delta$, then $Y_C=Y \simeq
\BC^{n-1}$.
\item
More significantly, let $W$ be a complex reflection
group of type $A_2$. Up to renormalisation,
the discriminant is $X_2^2+X_1^3$. Identify $Y$ with $\BC$.
For all $y\in \BC$, $\LL(y)=\{\pm (-y)^{3/2}\}$.
In particular, $\LL(1) = \{\pm \sqrt{-1}\}$.
Let $s$ be the simple element represented by the tunnel
with $y=1$, $x=-1$ and $L=2$.
Let $C$ be the chain $1\prec s$.
Then $Y_C$ is the open cone consisting of non-zero elements
of $\BC$ with argument in the open interval $(-2\pi/3, 2\pi/3)$.
\item Assume that $C\in \CC_n$.
All points in $Y_C$ are
generic.
Consider the map $Y^{\gen} \rightarrow \CC_n\times E_n$
sending $y$ to the pair $(\{1\prec s_1
 \prec s_1s_2 \prec \dots \prec s_1s_2\dots s_n=\delta\},\LL(y))$,
where $(s_1,\dots,s_n)=\lbl(y)$. This map is an homeomorphism.
The $(Y_C)_{C\in \CC_n}$ are the connected components of $Y^{\gen}$. Each
of these components is homeomorphic to $E_n^{\gen}$, which is contractible.
These $(Y_C)_{C\in \CC_n}$ are some analogues of chambers.
\end{itemize}

In the following proposition, if $A\subseteq \LL(y)$ is a submultiset,
the \emph{deep label} of $A$ is the sequence $(t_1,\dots,t_p)$
of labels (with respect to $y$) of points in the support of
 $A$ which are deep
with respect to $A$ (since these points may not be deep in $\LL(y)$,
the deep label of $A$ is not necessarily a subsequence of the deep
label of $y$).

\begin{lemma}
\label{lemmaunique}
Let $y\in Y$. 
Let $T=(y,z,L)$ and $T'=(y,z',L')$ be two tunnels in $L_y$
such that $b_{T}=b_{T'}$.
Then $T$ and $T'$ cross the same intervals among $I_1,\dots,I_n$
(in other words, $T$ and $T'$ are homotopic as tunnels drawn in 
$L_y$).
\end{lemma}

\begin{proof}
Up to perturbing $y$, we may assume that $y\in Y^{\gen}$.
Let $(s_1,\dots,s_n):=\lbl(y)$, and $(r_1,\dots,r_n):=\pi^n(\lbl(y))$.
Let $i_1,\dots,i_l$ (resp. $j_1,\dots,j_m$) be the successive indices
of the intervals among $I_1,\dots,I_n$ crossed
by $T$ (resp. $T'$). We have $b_T=s_{i_1}\dots s_{i_l}$ and $b_{T'}=
s_{j_1}\dots s_{j_m}$.
Assuming that $b_T=b_{T'}$, we obtain
$s_{i_1}\dots s_{i_l} = s_{j_1}\dots s_{j_m}$
and $r_{i_1}\dots r_{i_l} = r_{j_1}\dots r_{j_m}$.
Let $w:=r_{i_1}\dots r_{i_l}$.
By Lemma \ref{lemmalengthsimple}, $l=m$; both $(r_{i_1},\dots,r_{i_l})$ and
$(r_{i_1},\dots,r_{i_l})$ are reduced decompositions of $w$ and we have 
$$(*) \qquad \ker(w-1) = \bigcap_{k=1}^l \ker(r_{i_k}-1) =
 \bigcap_{k=1}^l \ker(r_{j_k}-1).$$
Assume that $(i_1,\dots,i_l) \neq (j_1,\dots,j_l)$. We may find
$j\in \{j_1,\dots,j_l\}$ such that, for example,
$i_1 < \dots < i_k < j < i_{k+1} < \dots i_l$.
Noting that the element $r_{i_1}\dots r_{i_k} r_jr_{i_{k+1}} \dots r_{i_l}$
is a parabolic Coxeter element, we deduce that $\ker(r_j-1) \not \supseteq
\ker(w-1)$. This contradicts $(*)$.
\end{proof}

\begin{prop}
\label{inequalities}
Let $C= \{1\prec c_1 \prec c_2 \prec \dots \prec c_m \}$ be a chain
in $\CC_m$. Let $y\in Y$, let $(x_1,\dots,x_k)$ be the ordered support
of $\LL(y)$ and $(s_1,\dots,s_k)$ be the label of $y$.
The following assertions are equivalent:
\begin{itemize}
\item[(i)] $y\in Y_C$.
\item[(ii)] There exists a partition of $\LL(y)$ into
$m+1$ submultisets $A_0,\dots,A_m$ such that:
\begin{itemize}
\item[(a)] For all $i,j\in \{1,\dots,m\}$ with $i<j$, for all
$x\in A_i$ and all $x'\in A_j$, we have $\re(x) < \re(x')$.
\item[(b)] For all $x\in A_0$, one has
$$\re(x) < \min_{x'\in A_1\cup \dots \cup A_m} \re(x')$$ or
$$\im(x) < \min_{x'\in A_1\cup \dots \cup A_m} \im(x')$$ or
$$\re(x) > \max_{x'\in A_1\cup \dots \cup A_m} \re(x').$$
\item[(c)] For all $i\in \{1,\dots,m\}$, the product of
the deep label of $A_i$ is $c_{i-1}^{-1}c_i$ (where one sets $c_0=1$).
\end{itemize}
\end{itemize}
Moreover, in situation (ii), the partition $\LL(y) = A_0\sqcup\dots\sqcup A_m$
is uniquely determined by $y$ and $C$.
\end{prop}

The picture below illustrates the proposition for particular $y$ and $C$.
Here $\lbl(y)=(s_1,\dots,s_5)$ and the considered chain is
$$C= \{1 \prec s_2 \prec s_2 s_4\}.$$
We have chosen tunnels $T_1=(y,x,L_1)$ and $T_2=(y,x,L_2)$, with $L_1< L_2$,
such that $b_{T_1}=s_2$ and $b_{T_2}=s_2s_4$.
The dotted lines represent these tunnels, as well as the 
vertical half-lines above $x$, $x+L_1$ and $x+L_2$. They partition
the complex line into three connected components; the partition
$A_0\sqcup A_1 \sqcup A_2$ is the associated partition of $\LL(y)$.
It is clear the possibility of drawing the tunnels is subject precisely
to the conditions on $A_0\sqcup A_1 \sqcup A_2$ expressed in the
proposition.
$$
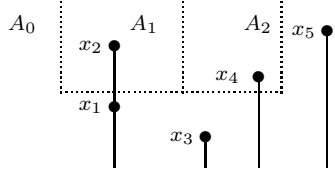
\xy
(-5,-10)="11",
(-5,-2)="1", (-5,6)="2", (7,-6)="3", (7,-10)="33",
(14,2)="4", (14,-10)="44",
(23,8)="5",(23,-10)="55",
(-17,9)*{_{A_0}},(-1,9)*{_{A_1}},(14,9)*{_{A_2}},
"1"*{\bullet},"2"*{\bullet},"3"*{\bullet},"4"*{\bullet},
"5"*{\bullet},(20,8)*{_{x_5}}, "55";"5" **@{-},
(-8,-2)*{_{x_1}},(-8,6)*{_{x_2}},(4,-6)*{_{x_3}},(10,2)*{_{x_4}},
"11";"1" **@{-},
"2";"1" **@{-},
"33";"3" **@{-},
"44";"4" **@{-},
(-12,0);(17,0) **@{.},
(-12,0);(-12,12) **@{.},
(4,0);(4,12) **@{.},
(17,0);(17,12) **@{.}
\endxy $$

\begin{proof}
(i) $\Rightarrow$ (ii):
$(y,z,L_1),\dots,(y,z,L_m)$ be tunnels representing the successive
non-trivial terms of $C=\{1\prec c_1 \prec \dots \prec c_m\}$.
 Set $z_0:=z$ and, for $i=1,\dots,m$, $z_i:=z+L_i$.
We have $$\re(z_0) < \re(z_1) < \dots < \re(z_m)$$ and
$$\im(z_0) = \im(z_1) = \dots = \im(z_m).$$
For $i\in\{1,\dots,m\}$, let $A_i$ be
the submultiset of $\LL(y)$ consisting of points $x$ such that
$\re(z_{i-1}) < \re(x) < \re(z_i)$ and 
$\im(x) > \im(z_0)$. Let $A_0$ be the complement in $\LL(y)$
of $A_1\cup \dots \cup A_m$. One easily checks (a), (b) and (c).

(ii) $\Rightarrow$ (i): Conversely, assume we are given a partition
$A_0\sqcup A_1\sqcup \dots \sqcup A_m$
satisfying conditions (a) and (b).
One may recover tunnels  $(y,z,L_1),\dots,(y,z,L_m)$ such that
the above construction yields the partition
$A_0\sqcup A_1\sqcup \dots \sqcup A_m$.
Condition (c) then implies that the tunnels
represent the elements of $C$, thus that $y\in Y_C$.

Uniqueness of the partition: this follows from condition (c) and
Lemma \ref{lemmaunique}.
\end{proof}

\begin{lemma}
\label{YC0}
Let $C= \{1\prec c_1 \prec c_2 \prec \dots \prec c_m \} \in \CC$.
Define $Y_C^0$ as the subspace of $Y_C$ consisting of points $y$ whose
associated partition $A_0,\dots,A_m$ (from Proposition \ref{inequalities}
(ii)) satisfies the following conditions:
\begin{itemize}
\item for $i=0,\dots,m$, the support of $A_i$ is a singleton $\{a_i\}$ and,
\item $\re(a_0) = \min_{i=1,\dots,m} \re(a_i) - 1$ and
$\im(a_0) = \min_{i=1,\dots,m} \im(a_i) - 1$.
\end{itemize}
Then $Y_C^0$ is contractible.
\end{lemma}

\begin{proof}
Let $y \in Y_C^0$. The support $(x_0,\dots,x_m)$
of $\LL(y)$ is generic, thus the label $(s_0,\dots,s_m)$ of $y$ coincides
with the deep label, and we have $s_0\dots s_m=\delta$ (Corollary
\ref{corodeeplabel}).
By Proposition \ref{inequalities}, Condition (ii)(c),
we have, for $i=1,\dots,m$, $s_i=c_{i-1}^{-1}c_i$ (where $c_0=1$).
Thus $s_0= \delta c_m^{-1}$.
We have proved that the label of any $y\in Y_C^0$ must
be $(\delta c_m^{-1},c_{0}^{-1}c_1,\dots,c_{m-1}^{-1}c_m)$.
A consequence of Theorem \ref{theolabel} is that
 the map $Y_C^0\rightarrow E_{m}^{\gen}$ sending
$y$ to $(x_1,\dots,x_m)$ is an homeomorphism.
One concludes with Lemma \ref{lemmaEngen}.
\end{proof}

We may now proceed to the proof of Proposition \ref{YCcontractible}.
Let $C=\{1\prec c_1 \prec \dots \prec c_m\} \in \CC$.
Let $y\in Y_C$. Let $(x_1,\dots,x_k)$ be the ordered support
of $\LL(y)$, and $(n_1,\dots,n_k)$ be the multiplicities.
Let $A_0,\dots,A_m$ be the partition of $\LL(y)$ described in 
Proposition \ref{inequalities} (ii).

The picture below gives an idea of the retraction of $Y_C$ onto $Y_C^0$ 
that will be explicitly constructed. It illustrates the motion of
a given point $y\in Y_C$; the black dots indicate the support of
$\LL(y)$ and the arrows how this support moves during the retraction.
$$
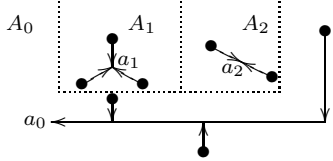
\xy
(-5,-4)="11",(-13,-4)="0",(-15,-4)*{_{a_0}},
(-5,-1)="1", (-5,7)="2", (-1,1)="22", (-9,1)="222",(-5,3.3)="20",
(-2.7,4)*{_{a_1}},
(7,-8)="3", (7,-4)="33",
(16,2)="4", (8,6)="44",(12,4)="40", (11,3)*{_{a_2}},
(23,8)="5",(23,-4)="55",
(-17,9)*{_{A_0}},(-1,9)*{_{A_1}},(14,9)*{_{A_2}},
"1"*{\bullet},"2"*{\bullet},"3"*{\bullet},"4"*{\bullet},
"22"*{\bullet},"222"*{\bullet},"44"*{\bullet},
"5"*{\bullet}, 
"2";"20" **@{-}, *\dir{>},
"22";"20" **@{-}, *\dir{>},
"222";"20" **@{-}, *\dir{>},
"4";"40" **@{-}, *\dir{>},
"44";"40" **@{-}, *\dir{>},
"1";"11" **@{-}, *\dir{>},
"3";"33" **@{-}, *\dir{>},
"5";"55" **@{-}, *\dir{>},
"55";"0" **@{-}, *\dir{>},
(-12,0);(17,0) **@{.},
(-12,0);(-12,12) **@{.},
(4,0);(4,12) **@{.},
(17,0);(17,12) **@{.}
\endxy $$

For $i=1,\dots,m$, consider the multiset mass center
$$a_i:= \frac{\sum_{x \in A_i} x}{|A_i|}$$
(in this expression, $A_i$ is viewed as a multiset: each $x_j$ in $A_i$ is taken
$n_j$ times, and $|A_i|$ is the multiset cardinal, \ie, the sum of the 
$n_j$ such that $x_j\in A_i$).

For each $t\in [0,1]$, let
$$\gamma_y(t):= A_0\cup \bigcup_{i=1}^m \{ (1-t) x + ta_i | x\in A_i\}$$
(here again, we consider the multiset union; in particular, the multicardinal
of $\gamma_y(t)$ is constant, equal to $n$ -- \ie, $\gamma_y(t)\in E_n$).
This defines a path in $E_n$.

As explained in Remark \ref{remarkhomotopylifting},
the path $\gamma_y$ uniquely lifts to a path
$\tilde{\gamma}_y$ in $Y$ such that $\tilde{\gamma}_y(0)=y$. An easy
consequence of Proposition \ref{inequalities} is that
$\tilde{\gamma}_y$ is actually drawn in $Y_C$.

Let $y':=\tilde{\gamma}_y(1)$.
Let $$R:=\min_{i=1,\dots,m} \re(a_i) = \re(a_1)$$ and
$$I:=\min_{i=1,\dots,m} \im(a_i).$$
For all $x \in A_0$, we have $\re(x) < R$ or $\im(x) < I$ or
$\re(x) > \re(a_m)$; denote by $x'$ the complex number with the same real part
as $x$ and imaginary part $I-1$; let $a_0:=R-1+\sqrt{-1}(I-1)$.
Consider the path $\beta_x: [0,1] \rightarrow \BC$ defined by:
$$\beta_x(t) := \left\{ \begin{matrix} (1-2t)x + 2t x' &
\text{if } t\leq 1/2 \\ 
(2-2t) x' + (2t-1) a_0 &  \text{if } t \geq 1/2. \end{matrix} \right.$$
For all $t\in[0,1]$, one has $\re(\beta_x(t)) < R$ or
$\im(\beta_x(t)) < I$ or
$\re(\beta_x(t)) > \re(a_m)$.

The path $\gamma_y':[0,1] \rightarrow E_n$ defined by 
$$\gamma_y'(t):= 
\{\underbrace{a_1,\dots,a_1}_{|A_1| \text{ times}}\} \cup \dots \cup
\{\underbrace{a_m,\dots,a_m}_{|A_m| \text{ times}}\} \cup
\bigcup_{x\in A_0} \beta_x(t) $$
lifts to
a unique path $\tilde{\gamma}_y'$ in $Y$ such that $\tilde{\gamma}_y'(0) = y'$.
Once again, a direct application of Proposition \ref{inequalities}
ensures that $\tilde{\gamma}_y'$ is actually in $Y_C$.
The endpoint $y'' := \tilde{\gamma}_y'(1)$ lies in the subspace
$Y_C^0$ of Corollary \ref{YC0}.

The map
\begin{eqnarray*}
\varphi:Y_C\times [0,1] & \longrightarrow & Y_C \\
(y,t) & \longmapsto & (\tilde{\gamma}_y' \circ \tilde{\gamma}_y)(t) 
\end{eqnarray*} is a retraction of $Y_C$ onto its
contractible subspace $Y_C^0$.
Thus $Y_C$ is contractible.

\section{The universal cover of $W \cq V^{\reg}$}
\label{sectionuniversalcover}

The main result of this section is:

\begin{theo}
\label{theouniversalcover}
\index{universal cover!of $W \cq V^{\reg}$}
The universal cover of $W \cq V^{\reg}$ is homotopy equivalent
to $\gar(B,S)$.
\end{theo}

Combined with Theorem \ref{gar3}, this proves our main result
Theorem \ref{theokapiun} in the irreducible case. The reducible case
follows.

To prove this theorem, we construct an open covering $(\widehat{\CU}_b)_{b\in B}$
such that
intersections of $(\widehat{\CU}_b)_{b\in B}$
are either empty or contractible (Proposition
\ref{intersections}). 
Under these assumptions, a standard theorem from algebraic topology
(\cite[4G.3]{hatcher})
shows that the universal cover is homotopy equivalent to the
\emph{nerve} \index{nerve!of an open covering}, \ie,
the simplicial space determined by non-empty intersections.
By showing that the nerve is $\gar(B,S)$ (Proposition \ref{nerve}),
we obtain the desired result.

As explained in Definition \ref{fatcover}, our ``basepoint'' $\CU$
provides us with a model denoted
$$\UC(W \cq V^{\reg},\CU)$$ for the universal cover of $W\cq V^{\reg}$.
Recall that, to any semitunnel $T$, one associates a path $\gamma_T$ whose
source is in $\CU$, thus a point in $\UC(W \cq V^{\reg},\CU)$.
Moreover, we have a left action of $B=\pi_1(W \cq V^{\reg}, \CU)$
on $\UC(W \cq V^{\reg},\CU)$.
With these conventions, our open covering is very easy to define:

\begin{defi}
\index{$\CVW$ (basic patch in the universal cover)}
The set $\widehat{\CU}_1$ 
is the subset of $\UC(W \cq V^{\reg},\CU)$
of elements represented by semitunnels.
For all $b\in B$, we set $\widehat{\CU}_b:= b \widehat{\CU}_1.$
\end{defi}

It is clear that two semitunnels represent the
same point in $\UC(W \cq V^{\reg},\CU)$ if and only if they are
equivalent in the following sense:

\begin{defi}
Two semitunnels $T=(y,x,L)$ and $T'=(y',x',L')$ are \emph{equivalent}
if and only if $y=y'$, $x+L=x'+L'$ and the affine segment $[(y,x),(y,x')]$ is
included in $\CU$.
\end{defi}

Let $T=(y,x,L)$ be a semitunnel.
The point of $\widehat{\CU}_1$ determined by $T$, or in other
words the equivalence class of $T$, is uniquely determined by $y$, $x+L$ and
$$\lambda(T):=\inf \left\{ \lambda'\in [0,L]
\left| [(y,x),(y,x+L-\lambda')] \subseteq \CU\right.
\right\}.$$
The number $\lambda$ is the infimum of the length of semitunnels in the
equivalence class of $T$. This infimum may not be a minimum, since
$(y,x+L-\lambda(T),\lambda(T))$ may not be a semitunnel (unless
$T$ is included in $\CU$).

We consider the following subsets of 
$\widehat{\CU}_1$:
\begin{itemize}
\item[($\CU_1$)] If 
$\lambda(T)=0 =\min \left\{ \lambda'\in [0,L] \left|
[x,x+L-\lambda'] \subseteq U_y\right.
\right\}$, then $T$ is a tunnel, equivalent to $(y,x+L,0)$. Elements
of $\widehat{\CU}_1$ represented by such tunnels of length $0$
form an open subset denoted by $\CU_1$. This subset is actually a sheet
over $\CU$ of the universal covering, corresponding to the trivial lift
of the ``basepoint''.
\item[($\overline{\CU}_1$)]
We denote by $\overline{\CU}_1$ the subset consisting of points represented
by semitunnels with $\lambda(T)=0$ (but without requiring that the
 ``$\inf$'' is actually a ``$\min$'').
We obviously have $\CU_{1} \subseteq \overline{\CU}_1$,
and $\overline{\CU}_1$ is contained in the closure of $\CU_1$.
\end{itemize}
Similarly, we set for all $b\in B$
$$\CU_b:=b\CU_1, \quad \overline{\CU}_b:=b\overline{\CU}_1.$$
Also, for any $y\in Y$, we denote by 
$\widehat{\CU}_{b,y}$
(resp. $\CU_{b,y}$, resp. $\overline{\CU}_{b,y}$) the intersection of 
$\widehat{\CU}_{b}$
(resp. $\CU_{b}$, resp. $\overline{\CU}_{b}$)
with the fiber over $y$ of the composed map
$\UC(W \cq V^{\reg},\CU) \rightarrow W \cq V^{\reg} \rightarrow Y$.

\begin{lemma}
The family $(\overline{\CU}_b)_{b\in B}$ is a partition of
$\UC(W \cq V^{\reg},\CU)$.
\end{lemma}

\begin{proof}
First, note that $b\neq 1 \Rightarrow \overline{\CU}_b \cap \overline{\CU}_1 = \varnothing$.
Indeed, let $p\in b\overline{\CU}_1 \cap \overline{\CU}_1$. Let $T=(y, x, L)$ be a semitunnel
representing $p$, let $T'=(y', x', L')$ be a semitunnel representing
$b^{-1}(p)$, let $\gamma$ be a path from $y$ to $y'$ representing $b$.
The paths $T$ and $\gamma T'$ represent the same point $p$ in the universal cover,
hence are homotopic.
By looking at the projection onto the base space, we see that $y+L = y'+L'$ and $x = x'$.
\begin{itemize}
\item 
If $L=0$ or $L'=0$, then $(y+L, x) = (y'+L', x')$ lies in $\CU$ and, as $\lambda(T)=\lambda(T')=0$,
we must have $L=L'=0$, hence $\gamma$ is homotopic to the trivial path.
\item 
If both $L> 0$ and $L'>0$, then for $\varepsilon$ small enough (\ie, such that $0< \varepsilon < \min(L,L')$),
\begin{itemize}
\item 
$T$ is the concatenation of
$T_0=(y,, x, L-\varepsilon)$ and $T_{\varepsilon}=(y+L-\varepsilon, x, \varepsilon)$,
\item $T'$ is the concatenation of
$T'_0=(y',, x', L'-\varepsilon)$ and $T'_{\varepsilon}=(y'+L'-\varepsilon, x', \varepsilon)=T_{\varepsilon}$.
\end{itemize}
Since $\lambda(T)=\lambda(T')=0$, $T_0$ and $T'_0$ are tunnels representing the trivial braid.
From $T\sim \gamma T'$, we deduce $T_0T_{\varepsilon}\sim \gamma T'_0T_{\varepsilon}$
and $T_0\sim \gamma T'_0$, which shows that $b=1$.
\end{itemize}

Thus it suffices to show that the projection $\overline{\CU}_1 \rightarrow
W \cq V^{\reg}$ is bijective.
Any point in $z\in W\cq V^{\reg}$ is the target of a unique equivalence
class of semitunnels $T$
with $\lambda(T)=0$; depending on whether $z\in \CU$ or not, the
associated point will be in $\CU_1$ or in $\overline{\CU}_1- \CU_1$.
\end{proof}

In particular, $(\widehat{\CU}_b)_{b\in B}$ is a covering of
$\UC(W \cq V^{\reg},\CU)$.

\begin{lemma}
\label{lemmacon}
For all $b\in B$, $\widehat{\CU}_b$ is open and contractible.
\end{lemma}

\begin{proof}
It is enough to deal with $b=1$.
That $\widehat{\CU}_1$ is open is easy.
Let $\CT$ be the space of semitunnels and $\sim$ the equivalence relation.
As a set, $\widehat{\CU}_1\simeq \CT/\sim$.
Consider the 
map
\begin{eqnarray*}
\phi:  \CT \times [0,1] & \longrightarrow  & \CT  \\
(T=(y,x,L),t) & \longmapsto & (y,x,L-t\lambda(T)) 
\end{eqnarray*}
If $T\sim T'$, then $\forall t, \phi(T,t)\sim \phi(T',t)$.
Thus $\phi$ induces a map $\overline{\phi}:\widehat{\CU}_1
\times [0,1] \rightarrow \widehat{\CU}_1$.
This map is continuous (this follows
from the fact that $\lambda$ induces a continuous function on
$\widehat{\CU}_1$) and $\overline{\phi} (T,t)=T$ if $T\in
\overline{\CU}_1$ or if $t=0$. We have proved that
$\widehat{\CU}_1$ retracts to $\overline{\CU}_1$.
We are left with having to prove that $\overline{\CU}_1$ is contractible.
Inside $\overline{\CU}_1$ lies $\CU_1$ which is contractible, since it
is a standard lift of the contractible space $\CU$.

We conclude by observing that
$\CU_1$ and $\overline{\CU}_1$ are homotopy equivalent. There
is probably a standard theorem from semialgebraic geometry applicable
here, but I was unable to
find a proper reference. Below is a ``bare-hand'' argument: it explains
how $\overline{\CU}_1$ may be ``locally retracted'' inside $\CU_1$
(constructing a global retraction seems difficult).

Both spaces have homotopy
type of $CW$-complexes and to prove homotopy equivalence it is enough
to prove that any continuous map $f:S^k \rightarrow \overline{\CU}_1$
(where $S^k$ is a sphere)
may be homotoped to a map $S^k \rightarrow {\CU}_1$.
Assume that there is a semitunnel $T=(y,z,L)$ such that $(T/\sim) \in f(S^k)
\cap (\overline{\CU}_1 - \CU_1)$.
We have $L>0$.
For any $\varepsilon > 0$, let $B_\varepsilon(y)$ be the open
ball of radius $\varepsilon$ in $Y$ around $y$, let $I_\varepsilon(z)$
be the affine interval $(z-\sqrt{-1}\varepsilon,z+\sqrt{-1}\varepsilon)$.
For $\varepsilon$ small enough, there is a unique continuous function
$L_{\varepsilon}: B_\varepsilon(y)\times I_\varepsilon(z) \rightarrow \BR$
such that $L_\varepsilon(y,z)=L$ and for all $(y',z')\in
 B_\varepsilon(y)\times I_\varepsilon(z)$, $(y',z',L_{\varepsilon}
(y',z'))$ represents a point in $\overline{\CU}_1 - \CU_1$.
The ``half-ball''
$$H_\varepsilon := \{(y',z',L') \in B_\varepsilon(y)\times
 I_\varepsilon(z) \times \BR | 0\leq L' \leq L_{\varepsilon}(y',z') \}$$
is a neighbourhood of $T/\sim$ in $\overline{\CU_1}$.
Working inside this neighbourhood, one may homotope $f$ to $f'$ such
that $f'(S^k)
\cap (\overline{\CU}_1 - \CU_1) \subseteq f(S^k)
\cap (\overline{\CU}_1 - \CU_1)-\{T/\sim\}$.
Compactness of $f(S^k)
\cap (\overline{\CU}_1 - \CU_1)$ guarantees that one can iterate this
process a finite number of times to get rid of all $f(S^k)
\cap (\overline{\CU}_1 - \CU_1)$.
\end{proof}

\begin{prop}
\label{nerve}
The nerve of $(\widehat{\CU}_b)_{b\in B}$ is
$\gar(B,S)$.
\end{prop}

\begin{proof}
Let $b,b'\in B$ such that
$\widehat{\CU}_b\cap \widehat{\CU}_{b'}\neq \varnothing$.
Let $T=(y,x,L)$ and $T'=(y',x',L')$ be semitunnels, representing
points $z$ and $z'$ in $\widehat{\CU}_1$, such that
$bz=b'z'$. The image of $z$ (resp. $z'$) in $W\cq V^{\reg}$ is
$(y,x+L)$ (resp. $(y',x'+L')$). Thus $x+L=x'+L'$ and $y=y'$.
Up to permuting $b$ and $b'$, we may assume that $L\geq L'$.
Since $x'=x+L-L'$ is in $U_y$, $T'':=(y,x,L-L')$ is a tunnel,
representing a simple element $b''$.
The tunnel $T$ is a concatenation of $T''$ and $T'$.
This implies that $z=b''z'$ and $b'z'=bz=bb''z'$.
By faithfulness of the $B$-action on the orbit of $z$,
we conclude that $bb''=b'$.

We have proved that the $1$-skeletons of the nerve and of
$\gar(B,S)$ coincide. To conclude, it remains to
check that the nerve is a flag complex. Let $C\subseteq B$ be such
that for all $b,b'\in C$, either $b^{-1}b'$ or $b'^{-1}b$ is simple.
Note that, unless $b=b'$, $b^{-1}b'$ and $b'^{-1}b$ cannot both be simple
(because Garside monoids are cancellative); in this setup, we write $b\preccurlyeq b'$ for 
$b^{-1}b'\in S$; this defines a total ordering on $C$ which, when both $b$ and $b'$
are in $M$, coincides
with the previously defined left prefix ordering.

We have to prove that $\bigcup_{b\in C} \widehat{\CU}_b\neq
\varnothing.$
Let $c_0,\dots,c_m$ be the elements of $C$, numbered according
to the total ordering on $C$ induced by $\preccurlyeq$:
$$c_0 \prec c_1 \prec c_2 \prec \dots \prec c_m \preccurlyeq c_0\delta.$$

Up to left-dividing each term by $c_0$, we may assume that $c_0=1$.
Let $y\in Y_C$. We may find $x,L_1,\dots,L_m$ such
that $(y,x,L_i)$ represents $c_i$ for all $i$.
The point of the universal cover represented by $(y,x,L_m)$
belongs to $\bigcup_{b\in C} \widehat{\CU}_b$.
\end{proof}

\begin{prop}
\label{intersections}
Let $C$ be a subset of $B$ such that
 $\bigcap_{b\in C}\widehat{\CU}_b \neq \varnothing $. Then 
 $\bigcap_{b\in C}\widehat{\CU}_b$ is contractible.
\end{prop}

\begin{proof}
As in the previous proof, we write $C=\{c_0,c_1,\dots,c_m\}$
with $$c_0 \prec c_1 \prec c_2 \prec \dots \prec c_m \preccurlyeq c_0\delta$$
and assume without loss of generality that $c_0=1$.

The case $m=0$ is Lemma \ref{lemmacon}.

Assume that $m\geq 1$. 
Let $T=(y,z,L)$. The point represented by
$T$ lies in $\bigcap_{b\in C}\widehat{\CU}_b$ if and only
if there exists $L_1,\dots,L_m$ with $0< L_1 < \dots < L_m < L$ such that,
for all $i$, $T_i:=(y,z,L_i)$ is a tunnel representing $c_i$.
Given $y\in Y$, it is possible to find such $L_1,\dots,L_m$ if and only
if $y\in Y_C$. This justifies:
$$\bigcap_{b\in C}\widehat{\CU}_{b,y} \neq \varnothing \Leftrightarrow
y \in Y_C.$$

For a given $y\in Y_C$, let us study
the intersection $\bigcap_{b\in C}\widehat{\CU}_{b,y}$.
Let $(x_1,\dots,x_k)$ be the ordered support of $\LL(y)$.
Let $A_0,\dots,A_m$ be the associated partition of $\LL(y)$,
defined in Proposition \ref{inequalities}.
Let $$R_+(y,C):=\max\{ \re(z) | z\in A_m\},$$
$$R_-(y,C):=\min\{ \re(z) | z\in A_1\},$$
 $$I_+(y,C):=\min\{\im(z) | z\in A_1\cup\dots \cup A_m\},$$
and
 $$I_-(y,C)=\sup\{\im(z) | z \in A_0 \text{ and } R_- \leq \re(z) \leq R_+\}.$$
We have $I_-(y,C) < I_+$. It may happen that $I_-(y,C)=-\infty$. 

We illustrate this on a picture. In the example,
the support of $\LL(y)$ is $x_1,\dots,x_6$
and $C=\{1\prec s_2 \prec s_2s_4\}$, where $(s_1,\dots,s_6)=\lbl(y)$.
The support of $A_1$ is $x_2$ and the support of $A_2$ is $x_4$.
The remaining points are in $A_0$.
The lines $\im(z)=I_-(y,C)$ and $\im(z)=I_+(y,C)$ are represented by full lines.
A semitunnel $(y,z,L)$
representing a point in $\bigcap_{b\in C}\widehat{\CU}_{b,y}$
must cross the intervals $I_2$ and $I_4$ represented by dotted lines.
One must have $\re(z) < R_-(y,C)$ and $I_-(y,C) < \im(z) < I_+(y,C)$; the final
point $z+L$ must satisfy $\re(z+L) > R_+(y,C)$. This final point
 may be any complex number
$z'$ in the rectangle $\re(z') > R_+(y,C)$ and $I_-(y,C) < \im(z')
< I_+(y,C)$, except
the points on the closed horizontal half-line to the right of $x_6$ (indicated
by a dashed line), which
cannot be reached.
$$
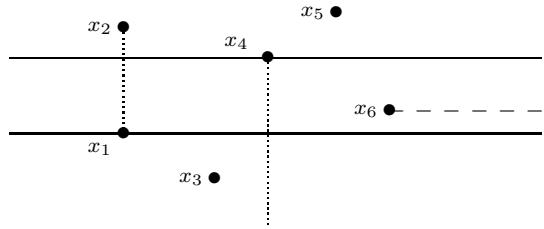
\xy
(-5,-10)="11",
(-5,-8)="1", (-5,6)="2", (7,-14)="3", (7,-10)="33",
(14,2)="4", (14,-20)="44",
(23,8)="5",(30,-5)="6",(50,-5)="66",
"1"*{\bullet},"2"*{\bullet},"3"*{\bullet},"4"*{\bullet},
"5"*{\bullet},(20,8)*{_{x_5}},
"6"*{\bullet},(27,-5)*{_{x_6}},
(-8,-10)*{_{x_1}},(-8,6)*{_{x_2}},(4,-14)*{_{x_3}},(10,4)*{_{x_4}},
"6";"66" **@{--},
(-20,-8);(50,-8) **@{-},
(-20,2);(50,2) **@{-},
"1";"2" **@{.},
"4";"44" **@{.}
\endxy $$

Generalising the example, one shows that $\bigcap_{b\in C}\widehat{\CU}_{b}$
may be identified with $$\bigcup_{y \in Y_C} E(y,C),$$
where $E(y,C)$ is the open rectangle of $\BC$ defined by
$$\re(z) > R_+(y,C), \quad I_-(y,C) < \im(z) < I_+(y,C),$$
from which have been removed the possible points of $A_0$ and the
horizontal half-lines to their rights.

Let $\overline{E}(y,C)$ be the rectangle
$$\re(z) \geq R_+(y,C), \quad I_-(y,C) < \im(z) < I_+(y,C)$$
from which have been removed the possible points of $A_0$ and the
horizontal half-lines to their rights. 

A homotopy argument similar to the one used in the proof of Lemma
\ref{lemmacon} (or possibly a nicer argument from semialgebraic geometry)
shows that $\bigcup_{y \in Y_C} E(y,C)$ and $\bigcup_{y \in Y_C} \overline{E}
(y,C)$ are homotopy equivalent.
The latter may be retracted to the union of open intervals
$$\bigcup_{y\in Y_C} (I_-(y,C),I_+(y,C))$$
(on each rectangle, the retraction is 
$$\overline{E}(y,C) \times [0,1] \rightarrow  \overline{E}(y,C),
(z,t) \mapsto R_+(y,C) + t(\re(z) -R_+(y,C)) 
+ \sqrt{-1} \im(z)).$$

The union $\bigcup_{y\in Y_C} (I_-(y,C),I_+(y,C))$ may be retracted
to $\bigcup_{y\in Y_C} [I_0(y,C),I_+(y,C))$, where
$I_0(y,C):= \max(\frac{I_-(y,C) +I_+(y,C)}{2},I_+(y,C)-1)$. The latter
space is a fiber bundle over $Y_C$, with fibers intervals.
Since the basespace $Y_C$ is contractible (Proposition \ref{YCcontractible}),
this fiber bundle is contractible. So is $\bigcap_{b\in C}\widehat{\CU}_b$.
\end{proof}

\section{Centralizers of regular elements}
\label{section11}

This section is devoted to the proof of Theorem \ref{theo31}: if $W'$ is
the centralizer of a $d$-regular element in an irreducible well-generated
complex reflection group $W$, then the hyperplane complement of $W'$ is
$K(\pi,1)$. To prove this, we develop a \emph{relative version} of
the tools and constructions presented in the previous sections, following
the same generic pattern of proof, but
with a categorical twist.

Let $m$ be a positive integer.
As explained in Appendix \ref{garsidependix} (Theorem \ref{theodivided}),
starting from any Garside structure --
we will start from
the dual braid monoid $M$, with its Garside element $\delta$
and Garside automorphism $\phi$ -- there exists a 
groupoid $\CG_m$ equivalent
as a category to $B$, and admitting a Garside structure $(M_m,\delta_m,\phi_m)$
where the order of $\phi_m$ is $m$ times
the order of $\phi$. Because $\phi_m$
is compatible with the Garside structure, the fixed subgroupoid
$\CG_m^{\phi_m}$ admits a Garside
structure $(M_m^{\phi_m}, \delta_m, \phi_m)$.

Our main result so far is that $W\cq V^{\reg}$ is a $K(B,1)$.
While $\mu_d$ naturally acts on $W\cq V^{\reg}$, the dual braid monoid $M$
may not have an automorphism of order $d$.
By replacing $M$ by $M_d$, or by any $M_{m}$ where $m$ is such that $d|mh$,
we can obtain a Garside structure with a symmetry of order $d$.
So the categorical viewpoint provides the algebraic structure that we need.

As $B$ and $\CG_m$ are equivalent, a $K(B,1)$ space
is homotopy equivalent to a $K(\CG_m, 1)$ space.
If one could view
$W\cq V^{\reg}$ as a $K(\CG_m,1)$
in a $\mu_d$-equivariant fashion, then the relative $K(\pi,1)$ property
should follow by simply considering fixed points.
Although we do not attempt to systematically investigate this intuition (the connection between
the cyclic/helicoidal structure, the non-positively curved aspects, and
geometric constructions
such as the Milnor fiber, would be especially worth studying), it is the
true explanation behind the miracles of the current section.

\subsection{A Garside category with symmetry of order $d$}

Since $W\cq V$ and $Y$ are defined by graded subalgebras of $\CO_V$,
these varieties are equipped with quotient actions of $\BC^{\times}$.
There is also a natural $\BC^{\times}$-action on the configuration
spaces $E_n$ and $\overline{E}_n^{\circ}$, defined by $$\lambda \cdot \pmb{\{ } z_1,\dots,z_n\pmb{\} }
:= \pmb{\{ } \lambda z_1,\dots,\lambda z_n\pmb{\} }.$$

\begin{lemma}
\label{lemmaLLd}
For all $(y,z) \in W\cq V$, for all $\lambda \in \BC^{\times}$, we have
$\LLL(\lambda (y,z)) = \lambda^h \LLL((y,z)).$
\end{lemma}

\begin{proof}
Elementary.
\end{proof}

If $d$ divides $h$, then $(W\cq V^{\reg})^{\mu_d}$
is again the regular orbit space of a well-generated complex
reflection group (because of Theorem \ref{theoregular} (2), the
centralizer $W'$ is again a duality group, thus is well-generated).
Because we have already proved the $K(\pi,1)$ property in this case,
it would be sufficient to focus on the case when $d$ does not divide $h$.

\begin{convention}
\label{convention11}
{Until the end of Section \ref{section11},
we work under the following assumptions
and notations:

\begin{itemize}
\item $d$ is a fixed regular
number;
\item we set
$$h':=\frac{h}{h\wedge d} \quad \text{and} \quad d':=\frac{d}{h\wedge d},$$
and in particular we have
$$h'd = hd' = h'd' (h \wedge d) =  h \vee d.$$
\end{itemize}
}
\end{convention}

Note that we do not assume that $d'>1$. Actually, taking $d'=d=1$, we will
obtain an alternate proof that the arrangement of $W$ is $K(\pi,1)$. This
alternate proof is geometrically simpler, but a bit more abstract, compared to the argument
provided in
the previous two sections.

\begin{lemma}
We have $d' \mid n$.
\end{lemma}

\begin{proof}
Because $d$ is regular, $(W\cq V^{\reg})^{\mu_d}$ is non-empty.
Let $(y,z) \in (W\cq V^{\reg})^{\mu_d}$.
Applying Lemma \ref{lemmaLLd} to $\lambda = \zeta_d$, we see
that $\LLL((y,z))$ is
$\zeta_d^h$-invariant. Note that $\zeta_d^h$ is a primitive $d'$-th
root of unity; orbits under multiplication by primitive $d'$-th roots of unity have either cardinal $d'$, or consists only of $0$.
The orbit $\{0\}$ is not permitted as $\LLL((y,z))\in \overline{E}_n^{\circ}$
Thus $n$, the cardinality of the multiset $\LLL(y)$,
is a multiple of $d'$.
\end{proof}

\begin{remark}
In the above proof, another consequence of the fact that $\LLL((y,z)) = LL(y) - z$ is
$\zeta_d^h$-invariant is that its barycenter $-z$ must be $0$. We have
$\LLL((y,z)) = \LL(y)$, and $(W\cq V^{\reg})^{\mu_d}$ can be identified with its image
in $Y$.
\end{remark}

To construct a Garside category with a symmetry of order $d$, we start
with the dual braid monoid $M$, together with its Garside element $\delta$,
Garside automorphism $\phi$ and set of simple elements $S$.

As explained in Appendix \ref{garsidependix},
all the structure is nicely encoded in the Garside set $\CD_{\bullet}(\delta)$ of factorizations
of $\delta.$
Through the isomorphism of Theorem \ref{bulletiso}, we can work instead in
$\CD_{\bullet}(c)$,
whose degree $k$ component is
$$\CD_k(c) = \{ (c_1,\dots,c_k)\in W^k | c_1\dots c_k = c
\text{ and } l_R(c_1) + \dots l_R(c_k) = l_R(c)
\}
$$
By reminding us that
this is just about combinatorics in a finite group,
working in $\CD_{\bullet}(c)$ helps reducing the cognitive load.
The automorphism $\phi$ acts on $\CD_k(c)$ by $$(c_1,\dots,c_k)\cdot \phi = (c_2,\dots,c_m,c_1^{c}).$$

To distinguish it from the topologically-defined braid group
$B$, we denote by $G$ the group of fractions of $M$
(a consequence of Theorem \ref{theogarside} is that $G\simeq B$).
As explained in Appendix \ref{garsidependix}, the components of $\CD_{\bullet}(c)$
whose degrees are multiples of $d'$ form, collectively, the Garside set
$$\sqrt[d']{\CD}_{\bullet}(c)$$
of a Garside category $M_{d'}$, whose groupoid of fractions $G_{d'}$ is equivalent to
$G$, and whose Garside automorphism $\phi_{d'}$,
has order $d'$ times the order of $\phi$.

\begin{defi}[relative dual Garside set]
\index{$\CD_{\bullet}'$ (relative dual Garside set for $W'$)}
\index{$M'$ (relative dual braid category for $W'$)}
\label{defidualcategory}
We set $$\CD'_{\bullet}(c) := \left(\sqrt[d']{\CD}_{\bullet}(c)\right)^{\phi_{d'}^{h'}}.$$
The \emph{dual Garside category} for $W'$ relative
to the pair $(W,W')$ is
the fixed subcategory $$M':=M_{d'}^{\phi_{d'}^{h'}}.$$
We denote by $G'$ the groupoid obtained by formally inverting
the morphisms in $M'$.
\end{defi}

By functoriality, we have the following commutative diagram:
$$\xymatrix{ M' \ar @{^{(}->}[d] \ar @{^{(}->}[r]  &
M_{d'} \ar @{^{(}->}[d] \ar @{->>}[r]^{\kappa_m} &
M \ar @{^{(}->}[d] \\
G' \ar[r] & G_{d'} \ar[r]^{\sim} & G \ar[d]^{\simeq}\\
& & B
}$$

The ``monomorphism'' arrows denote functors that are faithful (injective on morphisms.)
The functor $\kappa_m$, introduced in Definition \ref{deficollapse}, is
\emph{full} (surjective on morphisms), as simple elements are in its image.
The ``$\sim$'' arrows denote equivalence of categories.
The morphism from $G$ to $B$ is a true isomorphism (Theorem \ref{theogarside}.)

In Definition \ref{defibraidcategory} below, we will define topological groupoids
to complete the bottom line of this diagram.

\begin{theo}
The category $M'$ is a Garside category, with Garside set $\CD'_{\bullet}(c)$. As a consequence, $G'$
is a Garside groupoid.
\end{theo}

\begin{proof}
See Appendix \ref{garsidependix} and \cite{cyclic}: $M_{d'}$ is
a Garside category, and the fixed subcategory by an automorphism
of a Garside category is a Garside category (whose Garside set is the fixed subset
for the automorphism acting on the original Garside set).
\end{proof}

As explained in Remark \ref{remarkrecover},
one can use the low degree components of
$\CD'_{\bullet}(c)$ to write a presentation by generators and relations for $M'$ and $G'$.
\begin{itemize}
\item There is just one element in $\CD'_0(c)$; this element corresponds to the unique
object in $M$, and hasn't much algebraic significance for $M'$.
\item Elements in $\CD'_1(c)$ are $\phi_{d'}^{h'}$-invariant
factorizations $(c_1,\dots,c_{d'})$ of $c$.
They correspond to objects in $M'$.
\item Elements in $\CD'_2(c)$ are $\phi_{d'}^{h'}$-invariant
factorizations $(c_1,\dots,c_{2d'})$ of $c$.
They correspond to simple elements in $M'$.
\end{itemize}

The centralizer $W'$ (see Theorem \ref{theoregular}) may be badly-generated,
but the category $M'$ provides a substitute for the dual braid monoid.

\begin{example}
\label{exampleG31category}
In the situation of Example \ref{exampleg31} ($W=G_{37}$, $d=4$, $W'=G_{31}$),
we have $d'=2$ and $h'=15$.
The category $M_2$ consists of factorizations of $c$ as a product $c=uv$ of two elements. The automorphism $\phi_2^{15}$ acts on such pairs by
$$(u,v) \mapsto (v^{c^7},u^{c^8}).$$
The object set $\CD'_1(c)$ of $M'$, the relative
dual Garside category for $G_{31}$, is indexed by the solutions $(u,v)$ to
the following system of equations:
\begin{eqnarray*}
uv & = & c \\
l_R(u) + l_R(v) & = & l_R(c)  \; \;  = \; \; 8\\
u & = & v^{c^7} \\
v & = & u^{c^8}
\end{eqnarray*}
or, equivalently (since $c^{15}$ is central, the last two equations are equivalent):
\begin{eqnarray*}
u u^{c^8} & = & c \\
l_R(u) & = & 4
\end{eqnarray*}
Using a computer to sift through all elements of length $4$, one finds $88$ solutions.

The Garside automorphism $\phi_2$ acts on $\CD'_1(c)$  by $(u, u^{c^8})\mapsto (u^{c^8}, u^{c})$.
The object set $\CD'_1(c)$ decomposes into $\phi_2$-orbits: one orbit of size $3$,
2 orbits of size 5 and 5 orbits of size 15.
\end{example}

Example \ref{exampleG31category}
is the only case needed in our proof of the $K(\pi,1)$ conjecture,
but it does not harm to study the relative situation in full generality.

\subsection{Cyclic labels}

Departing from the earlier convention initiated
in Definition \ref{defisupport}, we now enumerate the distinct points
of a configuration not containing $0$
clockwise, starting from 12 o'clock plus $\varepsilon$ seconds,
and for points with identical
argument we enumerate them by increasing modulus, as in the following
example:

$$
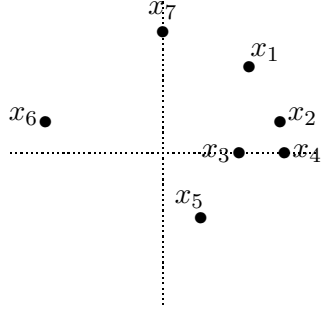
\xy {\xypolygon24"A"{~:{(10,0):}~>{}}},
{\xypolygon24"B"{~:{(16,0):}~>{}}},
{\xypolygon24"C"{~:{(7,0):}~>{}}},
{\xypolygon24"D"{~:{(19,0):}~>{}}},
{\xypolygon24"T"{~:{(20,0):}~>{}}},
"T0";"T1" **@{.},
"T0";"T7" **@{.},
"T0";"T19" **@{.},
"T0";"T13" **@{.},
"A1"*{\bullet},
"A21"*{\bullet},
"B4"*{\bullet},
"B12"*{\bullet},
"B2"*{\bullet},
"B1"*{\bullet},
"B7"*{\bullet},
"C1"*{x_3},
"C21"*{x_5},
"D4"*{x_1},
"D12"*{x_{6}},
"D2"*{x_2},
"D1"*{x_4},
"D7"*{x_{7}}
\endxy$$

\begin{defi}
Let $x\in \overline{E}_{n}^{\circ}$.
The above defined sequence $(x_1,\dots,x_{m})$ of distinct
points in $x$ is the \emph{cyclic support} of $x$.
For all $i$, we denote by $\theta_i$ the unique real number such
that  $0 < \theta_i \leq 2\pi$ and
$$e^{\sqrt{-1}\theta_i} x_i= \sqrt{-1} |x_i|$$
(in other words, a rotation with angle $\theta_i$ sends $x_i$ to the positive imaginary half-line). The non-decreasing
sequence $(\theta_1,\dots, \theta_{m})$
is the \emph{cyclic argument} of $x$.
The sequence $(n_1,\dots,n_m)$, where $n_i$ is the multiplicity of $x_i$, is the \emph{cyclic
multiplicity} of $x$.

The cyclic support (resp. argument, multiplicity) of
$(y,z)\in W\cq V^{\reg}$ is, by extension, defined as that
of $\LLL((y,z))$.
\end{defi}

In the above picture, we have $\theta_3=\theta_4=\pi/2$ and $\theta_7=2\pi$.

Let $x\in (W\cq V^{\reg})^{\mu_d}$ with associated
cyclic support $(x_1,\dots,x_{m})$ and cyclic argument
$(\theta_1,\dots, \theta_{m})$.
Let $\varepsilon > 0$.
By Lemma \ref{lemmaLLd}, the point
$e^{\sqrt{-1}(\theta_i - \varepsilon)} x_i$
lies in $\LLL(e^{\sqrt{-1}(\theta_i - \varepsilon)/h} x)$,
thus can be associated a simple element $s_{i,\varepsilon}$
that is part of the label of $e^{\sqrt{-1}(\theta_i - \varepsilon)/h} x$
(in the sense of Definition \ref{defilabel}).
As mentioned before, we prefer in this section to view simple elements in $W$
(via Proposition \ref{posetiso}).

Using the Hurwitz rule, one readily sees that
$s_{i,\varepsilon}$ does not depends on the choice
of a small enough $\varepsilon$.
Note for example that, for $i=1$, any value of $\varepsilon$ with $0 < \varepsilon
\leq \theta_1$ is suitable
and yields the same simple element. In particular, we have a well-defined sequence
$$(c_1,\dots, c_{m}) := (s_{1,\varepsilon},\dots, s_{m,\varepsilon})$$

\begin{defi}
\index{$\clbl$ (cyclic label)}
\index{label!cyclic label}
The \emph{cyclic label} of $(y,z)\in W\cq V^{\reg}$ is the sequence $\clbl(y,z) :=
(c_1,\dots, c_{m})$.
\end{defi}

\begin{lemma}
\label{lemmaclbl1}
For all $x\in W\cq V^{\reg}$, $\clbl(y,z)\in \CD_{\bullet}(c)$.
\end{lemma}

\begin{proof}
Consider the path:
$t \mapsto e^{\sqrt{-1}\frac{2\pi}{h}t}x$.
As explained after Definition \ref{defidelta}, this loop represents
the element $\delta \in B$.

This path is the concatenation of topologically trivial paths
and of paths of the form
\begin{eqnarray*}
\gamma_{y,i,\varepsilon}:[\theta_i -\varepsilon, \theta_i + \varepsilon]	& \to 	& W\cq V^{\reg} \\
t	& \mapsto & e^{\sqrt{-1}\frac{tL}{h}}x
\end{eqnarray*}
each of which, by Hurwitz rule, represents the simple element in $B$
corresponding to $c_i$.

The product of these simple elements is $\delta$, thus, after projecting
to $W$, $c_1\dots c_m = c$.
\end{proof}

\begin{lemma}
\label{lemmaonesteprotation}
Let $x\in W\cq V^{\reg}$ with cyclic support $(x_1,\dots,x_m)$,
cyclic argument $(\theta_1,\dots,\theta_m)$ and cyclic label
$(c_1,\dots,c_m)$. Let $l\leq m$ be the largest integer such that
$\theta_1 = \theta_2 = \dots = \theta_l$.
Then:
\begin{itemize}
\item[(i)] Let $(d_1,\dots,d_m)$ be the reduced label $\rlbl(x)$.
Let $\sigma\in \mathfrak{S}_m$
be the unique permutation such that the ordered
support of $\LLL(x)$ (in the sense of Definition \ref{defisupport})
is $(x_{\sigma(1)}, \dots, x_{\sigma(m)})$.
Then, for all integers $i$ from $1$ to $l$, we have $c_i = d_{\sigma^{-1}(i)}$.
\item[(ii)] $e^{\sqrt{-1}\frac{\theta_1}{h}} x$ has cyclic argument $(\theta_{l+1} - \theta_1,\dots,
\theta_m - \theta_1, 2\pi, \dots, 2\pi)$
\item[(iii)] $e^{\sqrt{-1}\frac{\theta_1}{h}} x$ has cyclic label
$(c_{l+1}, \dots, c_m, c_1^c, \dots, c_l^c)$
\end{itemize}
\end{lemma}

\begin{proof}
(i): as $\theta_1$ is the first element of the cyclic argument,
no point within $(x_1,\dots,x_l)$ is passed above by any other point in the support through
a rotation of angle $\theta_1 - \varepsilon$: using the Hurwitz rule, we see that
reduced labels of those points are preserved throughout the motion used to define the
first $l$ terms of the cyclic label.

(ii) is obvious by construction.

(iii): by construction and using (ii), we see that the last $l$ terms of the cyclic label
of $e^{\sqrt{-1}\frac{\theta_1}{h}} (y,z)$ are the first $l$ terms of the cyclic label of
$e^{\sqrt{-1}\frac{2\pi - \varepsilon}{h}}e^{\sqrt{-1}\frac{\theta_1}{h}} (y,z) = 
e^{\sqrt{-1}\frac{\theta_1 - \varepsilon}{h}}e^{\sqrt{-1}\frac{2\pi}{h}} (y,z) $.
We conclude by observing that, since
 the path $t \mapsto e^{\sqrt{-1}\frac{2\pi}{h} t} (y,z) $ represents $c$, the cyclic
label of of $e^{\sqrt{-1}\frac{2\pi}{h}} (y,z) $ is the conjugate label
$(c_1^c, \dots, c_m^c)$.
\end{proof}

\begin{defi}
\label{deficyclicallycompatible}
A configuration $x\in \overline{E}_{n}^{\circ}$ is \emph{cyclically compatible}
with $(c_1,\dots,c_m)\in \CD_{\bullet}(c)$ if the cyclic multiplicity of $x$ coincides with
$(l_R(c_1),\dots,l_R(c_m))$. We denote by 
$$\overline{E}_{n}^{\circ}
\boxdot \CD_{\bullet}(c) $$
the subspace of $\overline{E}_{n}^{\circ}\times \CD_{\bullet}(c)$ consisting of cyclically compatible pairs.
\end{defi}

\begin{prop}
\label{cycliccanonicalfull}
The map $(\LLL,\clbl)$ induces a bijection
\begin{eqnarray*}
W\cq V^{\reg}	& \longrightarrow 	& \overline{E}_{n}^{\circ} \boxdot \CD_{\bullet}(c) \\
x & \longmapsto & (\LLL(x), \clbl(x)).
\end{eqnarray*}
\end{prop}

\begin{proof}
That the image lies in $\overline{E}_{n}^{\circ} \boxdot \CD_{\bullet}(c)$ is obvious by construction.
Let $x\in W\cq V^{\reg}$ with cyclic argument $(\theta_1,\dots,\theta_m)$.
Consecutive applications of Lemma \ref{lemmaonesteprotation} to
$$y, e^{\sqrt{-1}\frac{\theta_1}{h}} x, \dots,e^{\sqrt{-1}\frac{\theta_{m-1}}{h}} x$$
 show
how to recover $\rlbl(x)$ from $\clbl(x)$, and vice-versa, via a sequence of permutations
and $c$-conjugacies
(note that the sequence of operations to get from $\rlbl(c)$ to $\clbl(c)$ depends on
$\LLL(x)$,
and it is much simpler to define it recursively using Lemma \ref{lemmaonesteprotation}
than it is to write down an explicit formula).
This provides a natural bijection $$\Phi: \overline{E}_{n}^{\circ} \boxtimes \CD_{\bullet}(c)\to
\overline{E}_{n}^{\circ} \boxdot \CD_{\bullet}(c)$$ that fits in a commutative diagram
$$
\xymatrix{
& \overline{E}_{n}^{\circ} \boxtimes \CD_{\bullet}(c) \ar[dd]^{\Phi} \\
 W\cq V^{\reg} \ar[ur]^{(\LLL,\rlbl)} \ar[dr]_{(\LLL,\clbl)} \\
& \overline{E}_{n}^{\circ} \boxdot \CD_{\bullet}(c) 
}
$$
Since the restriction of
$(\LLL,\rlbl)$ is bijective (Theorem \ref{nicetriv}),
the restriction of $(\LLL,\clbl)$ is the composition of two
bijections, hence is a bijection.
\end{proof}

In the above proof, the transition bijection
$\Phi: \overline{E}_{n}^{\circ} \boxtimes \CD_{\bullet}(c)\to
\overline{E}_{n}^{\circ} \boxdot \CD_{\bullet}(c)$ is a mere $\overline{E}_{n}^{\circ}$-dependent change of notation
on the  $\CD_{\bullet}(c)$-component. We equip $\overline{E}_{n}^{\circ} \boxdot \CD_{\bullet}(c)$
with the topology induced via $\Phi$ by that on $\overline{E}_{n}^{\circ} \boxtimes \CD_{\bullet}(c)$.
The trivialization of Proposition \ref{cycliccanonicalfull} is an homeomorphism.

\begin{lemma}
\label{lemmaclblm}
Let $x\in W\cq V^{\reg}$. Assume that $\LLL(x)\in (\overline{E}_{n}^{\circ})^{\mu_{d'}}$. Then
$\clbl(x)\in \sqrt[d']{\CD}_{\bullet}(c)$ and
$\clbl(\zeta_{d'h} x)=\clbl(x)\cdot{\phi_{d'}}.$
\end{lemma}

\begin{proof}
The path $t\mapsto e^{\sqrt{-1}\frac{2\pi}{d'h}t}x$ connects $x$ to $\zeta_{d'h} x$.
Since $\LLL(x)\in (\overline{E}_{n}^{\circ})^{\mu_{d'}}$, $\LLL(x)$ consists of $d'k$ distinct
points for some integer $k$.
The first $k$ points have cyclic argument less that $2\pi/{d'}$.
We have 
$$\clbl(x) = (c_1,\dots,c_{d'k})$$
and, by the very construction of $\clbl$, it is obvious that the first $d'(k-1)$ terms
of $\clbl(\zeta_{d'h} x)$
are
$(c_{k+1}, c_{k+2},\dots,c_{d'k}).$
The determination of the final $k$ terms
of $\clbl(\zeta_{d'h} x)$, and the equality $\clbl(\zeta_{d'h} x) = \clbl(x)\cdot
\phi_{d'}$, is an easy exercice using the Hurwitz rule and Lemma \ref{lemmaclbl1}.
\end{proof}

\begin{lemma}
\label{lemmadinvariant}
For all $x\in W\cq V^{\reg}$, the following assertions are equivalent:
\begin{itemize}
\item[(i)] $x\in (W\cq V^{\reg})^{\mu_d}$
\item[(ii)] $\LLL(x)\in (\overline{E}_{n}^{\circ})^{\mu_{d'}}$ and $\clbl(x) \in \CD'_{\bullet}(c)$.
\end{itemize}
\end{lemma}

\begin{proof}
Assume (i): we have $x = \zeta_d x$. By Lemma \ref{lemmaLLd},
$\LLL(\zeta_d x ) = \zeta_d^h \LLL(x)$,
so we have $\zeta_d^h \LLL(x) = \LLL(x)$. As $\zeta_d^h$ is a primitive root of order $d'$, we deduce
that $\LLL(x)\in (\overline{E}_{n}^{\circ})^{\mu_{d'}}$. We also have $\clbl(x) = \clbl(\zeta_d x) = \clbl(\zeta_{d'h}^{h'} x)$. Using Lemma \ref{lemmaclblm}, we deduce that
$\clbl(x) =\clbl(x)\cdot \phi_{d'}^{h'}$

Conversely, assuming $(ii)$, we conclude that $x$ and $\zeta_d x$
satisfy $\LLL(x) = \LLL(\zeta_d x)$
and, using Lemma \ref{lemmaclblm}, that
$\clbl(x) = \clbl(\zeta_d x)$. Using Proposition \ref{cycliccanonicalfull}, we obtain (i).
\end{proof}

Let us already note that, as a consequence, we obtain a very strong combinatorial property
of regular numbers:

\begin{coro}
\label{corodnonempty}
The object set $\CD'_{1}(c)$ is non-empty.
\end{coro}

By combining Proposition \ref{cycliccanonicalfull} and Lemma \ref{lemmadinvariant}, we obtain:

\begin{prop}
\label{cycliccanonical}
The map $(\LLL,\clbl)$ induces a bijection
\begin{eqnarray*}
(W\cq V^{\reg})^{\mu_d}	& \longrightarrow 	& (\overline{E}_{n}^{\circ})^{\mu_{d'}} \boxdot \CD'_{\bullet}(c) \\
y & \longmapsto & (\LLL(y), \clbl(y)).
\end{eqnarray*}
\end{prop}

\subsection{Basepoints and groupoids}

As explained in Definition \ref{fatty}, the
``fat basepoint'' $\CU$ is the set of pairs $(y,z)\in  Y\times\BC$
such that ``$z$ isn't below a point in $\LL(y)$''
or, equivalently, such that
$\LLL((y,z))$ does not contain any point in the closed
half-line $\sqrt{-1}\BR_{\geq 0}$:

$$\CU^{\mu_d} \simeq \{x \in (W\cq V^{\reg})^{\mu_d}| \LLL(x) \cap \sqrt{-1}\BR_{\geq 0} =
\varnothing\}.$$

This formula fails to capture the obvious symmetry constraints
deduced from Lemma \ref{lemmaLLd}:
if $x$ is $\mu_d$ invariant, then $\LLL(x)$ must be $(\mu_d)^h$-invariant,
hence $\mu_{d'}$-invariant. Hence:

$$\CU^{\mu_d} \simeq \{x \in (W\cq V^{\reg})^{\mu_d}
| \LLL(x) \cap \bigcup_{\zeta \in \mu_{d'}} \zeta
\sqrt{-1}\BR_{\geq 0} =
\varnothing\}.$$

We also introduce:
\index{$\CU_{d'}$}
$$\CU_{d'} := \{x\in W\cq V^{\reg} | \LLL(x) \cap \bigcup_{\zeta \in \mu_{d'}} \zeta
\sqrt{-1}\BR_{\geq 0} =
\varnothing\}.$$

Clearly, $$\CU^{\mu_d} \subseteq \CU_{d'} \subseteq \CU.$$

An illustration (with $d'=3$) of what $\LLL(x)$ may look like for some $x\in \CU^{\mu_d}$ (left)
and for some $x \in \CU_{d'} - \CU^{\mu_d}$ (right):

$$
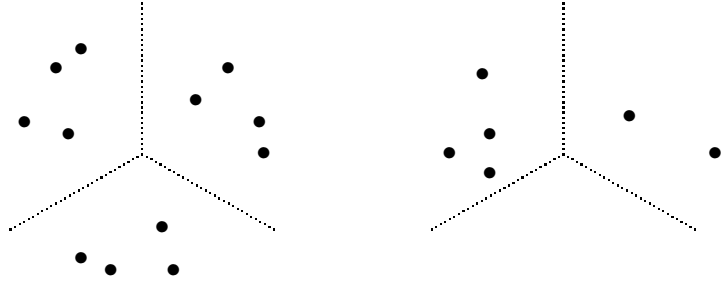
\xy {\xypolygon24"A"{~:{(10,0):}~>{}}},
{\xypolygon24"B"{~:{(16,0):}~>{}}},
{\xypolygon24"T"{~:{(20,0):}~>{}}},
"T0";"T23" **@{.},
"T0";"T7" **@{.},
"T0";"T15" **@{.},
"A4"*{\bullet},
"A12"*{\bullet},
"A20"*{\bullet},
"B4"*{\bullet},
"B12"*{\bullet},
"B20"*{\bullet},
"B2"*{\bullet},
"B10"*{\bullet},
"B18"*{\bullet},
"B1"*{\bullet},
"B9"*{\bullet},
"B17"*{\bullet}
\endxy
\qquad \qquad
\xy {\xypolygon24"A"{~:{(10,0):}~>{}}},
{\xypolygon24"B"{~:{(15,0):}~>{}}},
{\xypolygon24"C"{~:{(20,0):}~>{}}},
{\xypolygon24"T"{~:{(20,0):}~>{}}},
"T0";"T23" **@{.},
"T0";"T7" **@{.},
"T0";"T15" **@{.},
"A3"*{\bullet},
"C1"*{\bullet},
"A14"*{\bullet},
"A12"*{\bullet},
"B10"*{\bullet},
"B13"*{\bullet}
\endxy$$

By contrast with 
$\CU$, neither $\CU^{\mu_d}$  or $ \CU_{d'}$ is connected when $d'> 1$. However,
both happen to have contractible
connected components (see Lemma \ref{lemma_connected_objects} below), which makes them suitable as ``fat groupoid
basepoints'' (Definition \ref{fatgroupoid}).

Let $x\in \CU_{d'}$.
The support of $\LLL(x)$ can be partitioned into $d'$ (possibly empty) groups $B_1,\dots,B_{d'}$,
according to which $2\pi/d'$-sector they lie in: the
group $B_j$ consists of $m_j$ points 
$x_{\alpha_j},\dots,x_{\alpha_j + m_j}$ with
cyclic argument in $(\frac{2\pi}{d'}(j-1),\frac{2\pi}{d'}j)$.

\begin{defi}
\index{$\cc$ (cyclic content)}
\index{cyclic content}
The \emph{$d'$-cyclic content} of $x\in \CU_{d'}$
is the sequence $$\cc_{d'}(x) := (c_{\alpha_1} \dots c_{\alpha_1 + m_1},\dots,
c_{\alpha_{d'}} \dots c_{\alpha_{d'} + m_j}) \in \sqrt[d']{\CD}_{1}(c)$$
\end{defi}
Note that some of the $m_j$'s may be $0$.
Contrary to the cyclic label, which is always
non-degenerate, the cyclic content may contain trivial terms.

\begin{lemma}
\label{onemorelemma}
Let $x\in \CU_{d'}$.
Consider the path $[0,1]\to W\cq V^{\reg}, t\mapsto e^{\sqrt{-1}\frac{2\pi}{d'h}t}x$.
This path represents in $B$ the simple element associated with the first terms in
$\cc_{d'}(x)$. 
\end{lemma}

\begin{proof}
Obvious by construction of $\clbl$ and $\cc$.
\end{proof}

In the following definition, the existence and uniqueness of the standard image are
guaranteed by Proposition \ref{cycliccanonicalfull}:

\begin{defi}[standard image $x_{\sigma}$]
\label{standardimage}
Let $\sigma=(c_1,\dots,c_{d'k})\in \sqrt[d']{\CD}_{\bullet}(c)$, with $k\geq 1$.
The \emph{standard image $x_{\sigma}$} is the unique element of $\CU_{d'}$ such that,
for all $j$ such that $c_j\neq 1$, the point $e^{\sqrt{-1}\pi(\frac{1}{2} - \frac{2j-1}{d'k})}$
is in $\LLL(x_{\sigma})$ and the corresponding term in $\clbl(x_{\sigma})$ is $c_j$.
\end{defi}

Here are two examples with $d'=3$,
first with $k=1$ and $(c_1,c_2,c_3)$ not containing any trivial term, and then with $k=4$
and the trivial terms in $(c_1,\dots,c_{12})$ being $c_4$ and $c_7$:

$$
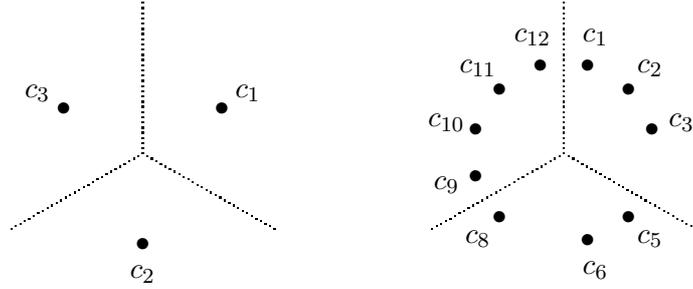
\xy {\xypolygon24"A"{~:{(12,0):}~>{}}},
{\xypolygon24"B"{~:{(16,0):}~>{}}},
{\xypolygon24"T"{~:{(20,0):}~>{}}},
"T0";"T23" **@{.},
"T0";"T7" **@{.},
"T0";"T15" **@{.},
"A3"*{\bullet},
"A11"*{\bullet},
"A19"*{\bullet},
"B3"*{c_1},
"B11"*{c_3},
"B19"*{c_2}
\endxy
\qquad \qquad
\xy {\xypolygon24"A"{~:{(12,0):}~>{}}},
{\xypolygon24"B"{~:{(16,0):}~>{}}},
{\xypolygon24"T"{~:{(20,0):}~>{}}},
"T0";"T23" **@{.},
"T0";"T7" **@{.},
"T0";"T15" **@{.},
"A2"*{\bullet},
"A4"*{\bullet},
"A6"*{\bullet},
"A8"*{\bullet},
"A10"*{\bullet},
"A12"*{\bullet},
"A14"*{\bullet},
"A16"*{\bullet},
"A20"*{\bullet},
"A22"*{\bullet},
"B2"*{c_3},
"B4"*{c_2},
"B6"*{c_1},
"B8"*{c_{12}},
"B10"*{c_{11}},
"B12"*{c_{10}},
"B14"*{c_9},
"B16"*{c_8},
"B20"*{c_6},
"B22"*{c_5}
\endxy$$

\begin{lemma}
Let $\sigma=(c_1,\dots,c_{d'k})\in \sqrt[d']{\CD}_{\bullet}(c)$, with $k\geq 1$. We have:
\begin{itemize}
\item[(i)] $\clbl(x_{\sigma})$ is obtained by removing trivial terms in $\sigma$,
\item[(ii)] $\cc_{d'}(x_\sigma) = (c_1\dots c_k,\dots,c_{k(d'-1)+1}\dots c_{kd'})$
\item[(iii)] if $\sigma\in\CD_{\bullet}'(c)$, then $x_{\sigma}\in \CU^{\mu_d}$.
\end{itemize}
\end{lemma}

\begin{proof}
Obvious by construction.
\end{proof}

\begin{lemma}
\label{lemma_connected_objects}
\label{lemma_connected_contractible}
\begin{itemize}
\item[(i)]
Let $C$ be a connected component of $\CU_{d'}$ (resp. $\CU^{\mu_d}$).
The map $\cc_{d'}$ is constant on $C$.
\item[(ii)] Let $C,C'$ be connected components of $\CU_{d'}$ (resp. $\CU^{\mu_d}$).
If $C\neq C'$, then $\cc_{d'}(C) \neq \cc_{d'}(C')$.
\item[(iii)] The map $\cc_{d'}$ restricts to
bijections $\pi_0(\CU_{d'})\simeq \sqrt[d']{\CD}_{1}(c)$ and
$\pi_0(\CU^{\mu_d})\simeq \CD'_{1}(c)$.
\item[(iv)] The connected components of $\CU_{d'}$ (resp. $\CU^{\mu_d}$)
are contractible.
\end{itemize}
\end{lemma}

\begin{proof}
(i). Let $x$ and $x'$ be two points in the same connected component of $\CU_{d'}$.
Starting from any path connecting $x$ to $x'$ within $\CU_{d'}$,
Lemma \ref{onemorelemma} produces an homotopy showing that
the first term of $\cc_{d'}(x)$ coincides with the first term of $\cc_{d'}(x')$.
Applying the same argument
to $\zeta x$ and $\zeta x'$
for $\zeta\in \mu_{d}$, we see that $\cc_{d'}(x) = \cc_{d'}(x')$.

(ii) and (iv): Let $x\in \CU_{d'}$.
Consider the standard element
$x_{\cc_{d'}(x)}$.
We construct a path $\gamma:[0,1]\to \CU_{d'}$ from $x$ to $x_{\cc_{d'}(x)}$
by sliding the points in $\LLL(x)$ in each
region $B_j$ affinely towards the corresponding point in $\LLL(x_{\cc_{d'}(x)})$, as
in the following picture:
$$
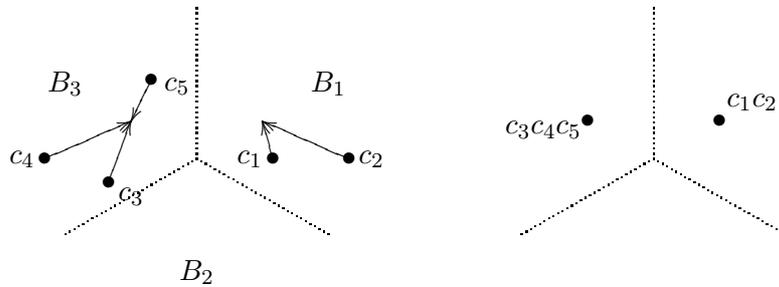
\xy {\xypolygon24"A"{~:{(10,0):}~>{}}},
{\xypolygon24"X"{~:{(7,0):}~>{}}},
{\xypolygon24"Y"{~:{(23,0):}~>{}}},
{\xypolygon24"B"{~:{(15,0):}~>{}}},
{\xypolygon24"U"{~:{(12,0):}~>{}}},
{\xypolygon24"C"{~:{(20,0):}~>{}}},
{\xypolygon24"T"{~:{(20,0):}~>{}}},
"T0";"T23" **@{.},
"T0";"T7" **@{.},
"T0";"T15" **@{.},
"A1"*{\bullet},
"X1"*{c_1},
"Y1"*{c_2},
"C1"*{\bullet},
"U14"*{\bullet},
"A15"*{c_3},
"U9"*{\bullet},
"A8"*{c_5},
"C13"*{\bullet},
"Y13"*{c_4},
"C3"*{B_1},
"C11"*{B_3},
"B19"*{B_2},
"A1";"A3" **@{-}, *\dir{>},
"C1";"A3" **@{-}, *\dir{>},
"U14";"A11" **@{-}, *\dir{>},
"U9";"A11" **@{-}, *\dir{>},
"C13";"A11" **@{-}, *\dir{>}
\endxy
\qquad \qquad
\xy {\xypolygon24"A"{~:{(10,0):}~>{}}},
{\xypolygon24"B"{~:{(15,0):}~>{}}},
{\xypolygon24"C"{~:{(20,0):}~>{}}},
{\xypolygon24"T"{~:{(20,0):}~>{}}},
"T0";"T23" **@{.},
"T0";"T7" **@{.},
"T0";"T15" **@{.},
"A3"*{\bullet},
"A11"*{\bullet},
"B3"*{c_1c_2},
"B12"*{c_3c_4c_5},
\endxy$$
This proves that $x$ and $x_{\cc_{d'}(x)}$ lie in the same connected component
of $\CU_{d'}$ (resp. $\CU^{\mu_d}$, since if $x$ is $\mu_d$-invariant, so is the path
$\gamma$.) (ii) follows: if $x,x'$ are such that $\cc_{d'}(x)=\cc_{d'}(x')$, then
$x$ and $x'$ must lie in the same connected component. Using (i), we also notice
that the path $\gamma$ is part of a deformation-retraction of the full connected
component of $x$ onto the single point $x_{\cc_{d'}(x)}$: this proves (iv).

(iii): combining (i) and (ii), we see that $\cc_{d'}$ induces a bijection
from $\pi_0(\CU_{d'})$ onto its image in $\sqrt[d']{\CD}_{1}(c)$ (resp. from
$\pi_0(\CU^{\mu_d})$ onto its image in $\CD'_{1}(c)$).
To prove (iii), it suffices to check that, for any decomposition in $\sigma\in\sqrt[d']{\CD}_{1}(c)$
(resp. $\CD'_{1}(c)$), there exists a point $x$ in $\CU_{d'}$ (resp. $\CU^{\mu_d}$)
such that $\cc_{d'}(x) = \sigma$; the standard image $x_{\sigma}$ provides a particular
example of such a point.
\end{proof}

A consequence of Lemma \ref{lemma_connected_contractible} (iv)
is that we can use $\CU_{d'}$ and
$\CU^{\mu_d}$ as a ``fat groupoid basepoint''
(see Definition \ref{fatgroupoid}):

\begin{defi}
\label{defibraidcategory}
\index{$\Braid '$ (relative braid groupoid for $W'$)}
The relative braid category associated with $W'$ is the
groupoid
$$B' := \pi_1((W\cq V^{\reg})^{\mu_d}, \CU^{\mu_d}).$$
We also set
$$B_{d'} := \pi_1(W\cq V^{\reg},\CU_{d'}).$$
\end{defi}

By functoriality of $\pi_1$, we have natural functors:
$$\xymatrix{ 
B' \ar[r] & B_{d'} \ar[r]^{\sim} & B }$$
Note that the functor $B_{d'}\to B$ is an \emph{equivalence of categories} and
not an isomorphism: it is not injective on objects.

\subsection{Circular tunnels and semitunnels}

\begin{defi}
\label{deficircularsemitunnel}
A \emph{circular semitunnel} is a
element $T=(x, L)$ in $$\CU_{d'}\times [0,\frac{2\pi}{d'h}].$$

We say that $L$ is the length of the semitunnel.

The path $\gamma_T$ associated with $T$ is the path
$[0,1] \to (W\cq V^{\reg})^{\mu_d}, t\mapsto e^{\sqrt{-1} t}x$.

The circular semitunnel $T$ is a \emph{circular tunnel} if it
satisfies the additional condition
$$e^{\sqrt{-1} L} x \in \CU_{d'}.$$
\index{tunnels!circular semitunnels}
\end{defi}

Let $T = (x, L)$ be a circular (semi)tunnel. Because $\LLL$ has degree $h$,
the angular rotation of
$\LL(\gamma_T)$ throughout the motion is $hL$. Here is an example of
tunnel of length $\frac{\pi}{3h}$:

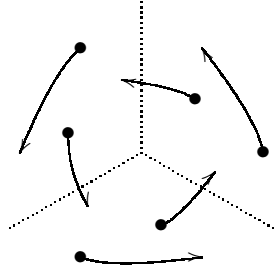
$$\xy {\xypolygon24"A"{~:{(10,0):}~>{}}},
{\xypolygon24"B"{~:{(16,0):}~>{}}},
{\xypolygon24"T"{~:{(20,0):}~>{}}},
"T0";"T23" **@{.},
"T0";"T7" **@{.},
"T0";"T15" **@{.},
"A4"*{\bullet},
"A12"*{\bullet},
"A20"*{\bullet},
"B1"*{\bullet},
"B9"*{\bullet},
"B17"*{\bullet},
"A12";"A16" **\crv{"A14"},*\dir{>},
"A20";"A24" **\crv{"A21"},*\dir{>},
"A4";"A8" **\crv{"A5"},*\dir{>},
"B1";"B5" **\crv{"B2"},*\dir{>},
"B9";"B13" **\crv{"B10"},*\dir{>},
"B17";"B21" **\crv{"B18"},*\dir{>},
\endxy$$

\begin{defi}
Let $\sigma = (c_1,\dots,c_{2d'})\in \sqrt[d']\CD_2(c)$.
The circular tunnel
$$(x_{\sigma}, \frac{\pi}{d'h})$$
defines an element of $B_{d'}$ which we denote by $b_{\sigma}$.

If $\sigma\in \CD'_2(c)$, then the circular tunnel also represents an
element of $B'$, which we denote by $b'_{\sigma}$.
\end{defi}

Later on, we will prove that the natural map $B'\to B_{d'}$ is injective,
which will allow us to identify $b'_{\sigma}$ with $b_{\sigma}$, but at this stage
we are not supposed to know this.

$$
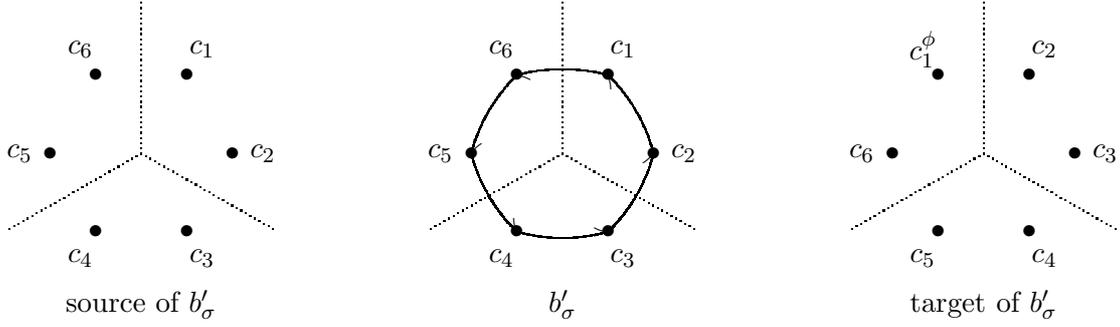
\xy {\xypolygon24"A"{~:{(12,0):}~>{}}},
{\xypolygon24"B"{~:{(16,0):}~>{}}},
{\xypolygon24"T"{~:{(20,0):}~>{}}},
"T0";"T23" **@{.},
"T0";"T7" **@{.},
"T0";"T15" **@{.},
"A1"*{\bullet},
"A5"*{\bullet},
"A9"*{\bullet},
"A13"*{\bullet},
"A17"*{\bullet},
"A21"*{\bullet},
"B1"*{c_2},
"B5"*{c_1},
"B9"*{c_6},
"B13"*{c_5},
"B17"*{c_4},
"B21"*{c_3},
"T19"*{\text{source of $b'_{\sigma}$}},
\endxy
\qquad \qquad\xy {\xypolygon24"A"{~:{(12,0):}~>{}}},
{\xypolygon24"B"{~:{(16,0):}~>{}}},
{\xypolygon24"T"{~:{(20,0):}~>{}}},
"T0";"T23" **@{.},
"T0";"T7" **@{.},
"T0";"T15" **@{.},
"A1"*{\bullet},
"A5"*{\bullet},
"A9"*{\bullet},
"A13"*{\bullet},
"A17"*{\bullet},
"A21"*{\bullet},
"B1"*{c_2},
"B5"*{c_1},
"B9"*{c_6},
"B13"*{c_5},
"B17"*{c_4},
"B21"*{c_3},
"T19"*{b'_{\sigma}},
"A5";"A9" **\crv{"A7"},*\dir{>},
"A9";"A13" **\crv{"A11"},*\dir{>},
"A13";"A17" **\crv{"A15"},*\dir{>},
"A17";"A21" **\crv{"A19"},*\dir{>},
"A21";"A1" **\crv{"A23"},*\dir{>},
"A1";"A5" **\crv{"A3"},*\dir{>}
\endxy
\qquad \qquad
\xy {\xypolygon24"A"{~:{(12,0):}~>{}}},
{\xypolygon24"B"{~:{(16,0):}~>{}}},
{\xypolygon24"T"{~:{(20,0):}~>{}}},
"T0";"T23" **@{.},
"T0";"T7" **@{.},
"T0";"T15" **@{.},
"A1"*{\bullet},
"A5"*{\bullet},
"A9"*{\bullet},
"A13"*{\bullet},
"A17"*{\bullet},
"A21"*{\bullet},
"B1"*{c_3},
"B5"*{c_2},
"B9"*{c_1^\phi},
"B13"*{c_6},
"B17"*{c_5},
"B21"*{c_4},
"T19"*{\text{target of $b'_{\sigma}$}},
\endxy$$

\begin{lemma}
\begin{itemize}
\item[(i)]
The map $\sqrt[d']{\CD}_2(c) \to B_{d'}, \sigma \mapsto b_{\sigma}$ extends to a groupoid
morphism $$\psi:G_{d'} \longrightarrow B_{d'}.$$
\item[(ii)]
The map ${\CD}'_2(c) \to B', \sigma \mapsto b'_{\sigma}$ extends to a groupoid
morphism $$\psi:G' \longrightarrow B'.$$
\end{itemize}
\end{lemma}

\begin{proof}
(i).
The groupoid $G_{d'}$ has a presentation with generators indexed by $\sqrt[d']{\CD}_2(c)$
and relations indexed by $\sqrt[d']{\CD}_3(c)$.
Let $\rho=(c_1,\dots,c_{3d'})\in \sqrt[d']{\CD}_3(c)$.
It corresponds to the relation $\sigma\tau = \mu$, where
$$\sigma = (c_1,c_2c_3, c_4, c_5c_6, \dots,c_{3d'-2}, c_{3d'-1}c_{3d'})$$
$$\tau = (c_2,c_3c_4, c_5, c_6c_7, \dots,c_{3d'-1},c_{3d'}c_1^\phi)$$
$$\mu =(c_1c_2, c_3, c_4c_5,c_6, \dots,c_{3d'-2}c_{3d'-1},c_{3d'})$$

Consider the standard element $x_{\rho}$:
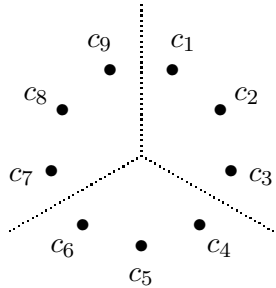
$$
\xy {\xypolygon36"A"{~:{(12,0):}~>{}}},
{\xypolygon36"B"{~:{(16,0):}~>{}}},
{\xypolygon24"T"{~:{(20,0):}~>{}}},
"T0";"T23" **@{.},
"T0";"T7" **@{.},
"T0";"T15" **@{.},
"A4"*{\bullet},
"A8"*{\bullet},
"A12"*{\bullet},
"A16"*{\bullet},
"A20"*{\bullet},
"A24"*{\bullet},
"A28"*{\bullet},
"A32"*{\bullet},
"A36"*{\bullet},
"B4"*{c_2},
"B8"*{c_1},
"B12"*{c_9},
"B16"*{c_8},
"B20"*{c_7},
"B24"*{c_6},
"B28"*{c_5},
"B32"*{c_4},
"B36"*{c_3}
\endxy
$$
The needed relation $b_{\sigma}b_{\tau} = b_{\mu}$
follows from the following easy observations:
\begin{itemize}
\item the circular tunnel $T_{\sigma}:=(x_{\rho},\frac{2\pi}{3d'h})$
represents $b_{\sigma}$,
\item the circular tunnel $T_{\tau}:=(\zeta_{3d'h} x_{\rho},\frac{2\pi}{3d'h})$
represents $b_{\tau}$,
\item the circular tunnel $T_{\mu}:=(x_{\rho},\frac{4\pi}{3d'h})$
represents $b_{\mu}$,
\item $T_{\mu}$ is the composition of $T_{\sigma}$ and $T_{\mu}$.
\end{itemize}

(ii). In the above construction, if $\rho\in \CD'_3(c)$, then $\sigma,\tau,\mu\in
\CD'_2(c)$ and the same argument applies. 
\end{proof}

\begin{lemma}
Let $\sigma = (c_1,\dots,c_{2d'})\in \sqrt[d']\CD_2(c)$.
Then the image of $b_{\sigma}$ via the natural functor $B_{d'}\to B$ is the
simple element corresponding to $c_1$ (using the notation introduced in
Definition \ref{deficollapse}: it is the image of $\kappa_m(\sigma)$ under the natural
embedding $M\hookrightarrow B$.)
\end{lemma}

\begin{proof}
This is obvious by definition of $\clbl$.
\end{proof}

We have a commutative diagram of functors:

$$\xymatrix{ M' \ar @{^{(}->}[d] \ar @{^{(}->}[r] &
M_{d'} \ar @{^{(}->}[d] \ar @{->>}[r]^{\kappa_m} &
M \ar @{^{(}->}[d] \\
G' \ar[d]^{\psi'} \ar @{^{(}->}[r]  &
G_{d'} \ar[r]^{\sim} \ar[d]^{\psi}  &  G \ar[d]^{\simeq} \\
B' \ar[r] & B_{d'} \ar[r]^{\sim}& B.}$$

\begin{theo}
\label{theopsi}
The morphism $\psi:G_{d'} \to B_{d'}$ is a groupoid isomorphism.
It restricts to a groupoid isomorphism $\psi': G' \to B'$.
\end{theo}

\begin{proof}
Let $o,o'$ be objects in $G_{d'}$. Since we have category equivalences $G_{d'}\sim G$ and
$B_{d'}\sim B$, we have a commutative diagram of set-theoretic maps:
$$\xymatrix{\Hom_{G_{d'}}(o,o') \ar[rr]^(.65){\simeq} \ar[d]^{\psi}  & & G \ar[d]^{\simeq} \\
\Hom_{B_{d'}}(\psi(o),\psi(o')) \ar[rr]^(.65){\simeq}& & B}$$
which shows that $\psi$ is an equivalence of categories.

By Lemma \ref{lemma_connected_objects} (iii), both $\psi$ and $\psi'$ are bijective
on objects. In particular, $\psi$ is an isomorphism of categories.

The subdiagram
$$\xymatrix{G' \ar[d]^{\psi'} \ar @{^{(}->}[r]  &
G_{d'}  \ar[d]^{\simeq}  \\
B' \ar[r] & B_{d'}
}$$
shows that $\psi'$ is faithful.

We are left with having to prove that $\psi'$ is full.
This follows from a generic position argument. Let $b'$ a morphism in $B'$. It can represented
by a path $\gamma: [0,1] \to (W\cq V^{\reg})^{\mu_d}$, and
we can assume that at any given $t\in [0,1]$, at most one point in $\LLL(\gamma(t))$ lies
on the vertical half-line $\sqrt{-1}\BR_{\geq 0}$
(as points not satisfying this form a subspace of real codimension $2$ in
$(W\cq V^{\reg})^{\mu_d}$); this expresses $\gamma$ as a concatenation of paths
homotopic to circular tunnel paths in $B'$.
\end{proof}

\begin{remark}
Using Theorem \ref{theopsi}, one can write down a presentation
by generators and relations for $B'$, and
deduce a presentation and relations for the braid group of $W'$.
\end{remark}

\subsection{The universal cover}

Let us choose a base object $o\in \CD'_1(c)$. By Lemma \ref{lemma_connected_objects} (iii), 
it corresponds to
a connected component of $\CU^{\mu_d}$ which, by abuse of notation, we still denote by $o$.
By Lemma \ref{lemma_connected_contractible} (iv), this component is contractible and, using
Definition \ref{fatcover}, we get a model $\UC((W\cq V^{\reg})^{\mu_d}, o)$ for the universal cover
of $(W\cq V^{\reg})^{\mu_d}$.

\begin{defi}
\label{boo}
\index{$\CVWr$ (basic patch in the universal cover, relative version)}
Let $g \in \obj(o \downarrow B')$ or, in other words,
let $g$ be a morphism in $B'$ with source $o$.
Let $o'$ be the target of $g$.
Let $\gamma$ be a path $(W\cq V^{\reg})^{\mu_d}$ representing $g$.
Let $T=(x,L)$ be a circular semitunnel such that $x$ lies in the connected component
(corresponding to) $o'$.
Then $\gamma$ and $\gamma_T$ can be composed (up to homotopically trivial glue binding them)
and $\gamma \gamma_T$ represents a point $p\in \UC((W\cq V^{\reg})^{\mu_d}, o)$.
We denote by
$$\CV_g$$ the subspace of $\UC((W\cq V^{\reg})^{\mu_d}, o)$
consisting of points $p$ that can be obtained this way, by concatenating a path representing
$g$ with a circular semitunnel.
\end{defi}

The following result is a variation on Proposition \ref{YCcontractible}.

\begin{lemma}
\label{lemmareprcontr}
Let $\sigma$ be a non-degenerate element of $\CD'_{\bullet}(c)$.
Then the subspace $U_{\sigma}\subseteq\CU^{\mu_d}$ consisting of points $x$ such that
$\clbl(x) \vdash \sigma$ is contractible.
\end{lemma}

The notation $\vdash$ was introduced in Definition \ref{defiface}.

\begin{proof}
Write $\sigma = (c_1,\dots,c_{md'})$.  Clearly, $x_{\sigma}\in U_{\sigma}$, so $U_{\sigma}\neq \varnothing$.

Let $x\in U_{\sigma}$. We have $\clbl(x) = (d_1,\dots,d_{pd'})$, with $p\geq m$.
Each $c_j$ is the
product of consecutive $d_{\alpha(j)}\dots d_{\alpha(j) + \beta(j)}$. This defines
a mapping $\{1,\dots,pd'\} \to \{1,\dots,md'\}$, mapping each $i$ to the unique $j$ such
that $\alpha(j)\leq i \leq \alpha(j) + \beta(j)$.
By moving the point $i$-th point in $\LLL(x)$ continuously towards the $j$-th point
in $\LLL(x_{\sigma})$, one obtains a deformation-retraction of $U_{\sigma}$ onto the single
point $x_{\sigma}$.

In the illustration below, we have chosen $p=5$ and $m=3$:

$$
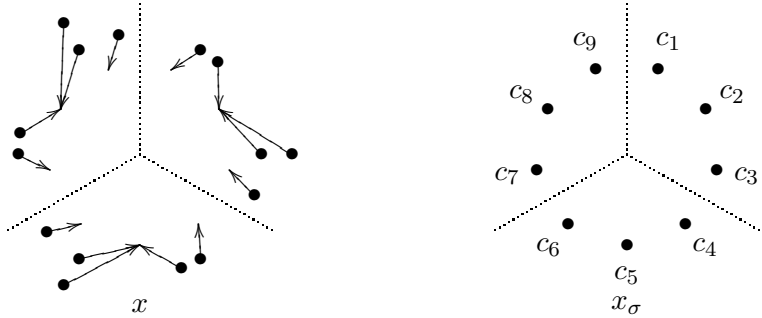
\xy {\xypolygon36"A"{~:{(12,0):}~>{}}},
{\xypolygon36"B"{~:{(16,0):}~>{}}},
{\xypolygon36"T"{~:{(20,0):}~>{}}},
"T0";"T10" **@{.},
"T0";"T22" **@{.},
"T0";"T34" **@{.},
"T1"*{\bullet},
"B1"*{\bullet},
"B6"*{\bullet},
"B7"*{\bullet},
"B35"*{\bullet},
"T1";"A4" **@{-}, *\dir{>},
"B1";"A4" **@{-}, *\dir{>},
"B6";"A4" **@{-}, *\dir{>},
"B7";"A8" **@{-}, *\dir{>},
"B35";"A36" **@{-}, *\dir{>},
"T13"*{\bullet},
"B13"*{\bullet},
"B18"*{\bullet},
"B19"*{\bullet},
"B11"*{\bullet},
"T13";"A16" **@{-}, *\dir{>},
"B13";"A16" **@{-}, *\dir{>},
"B18";"A16" **@{-}, *\dir{>},
"B19";"A20" **@{-}, *\dir{>},
"B11";"A12" **@{-}, *\dir{>},
"T25"*{\bullet},
"B25"*{\bullet},
"B30"*{\bullet},
"B31"*{\bullet},
"B23"*{\bullet},
"T25";"A28" **@{-}, *\dir{>},
"B25";"A28" **@{-}, *\dir{>},
"B30";"A28" **@{-}, *\dir{>},
"B31";"A32" **@{-}, *\dir{>},
"B23";"A24" **@{-}, *\dir{>},
"T28"*{x}
\endxy
\qquad \qquad \qquad
\xy {\xypolygon36"A"{~:{(12,0):}~>{}}},
{\xypolygon36"B"{~:{(16,0):}~>{}}},
{\xypolygon24"T"{~:{(20,0):}~>{}}},
"T0";"T23" **@{.},
"T0";"T7" **@{.},
"T0";"T15" **@{.},
"A4"*{\bullet},
"A8"*{\bullet},
"A12"*{\bullet},
"A16"*{\bullet},
"A20"*{\bullet},
"A24"*{\bullet},
"A28"*{\bullet},
"A32"*{\bullet},
"A36"*{\bullet},
"B4"*{c_2},
"B8"*{c_1},
"B12"*{c_9},
"B16"*{c_8},
"B20"*{c_7},
"B24"*{c_6},
"B28"*{c_5},
"B32"*{c_4},
"B36"*{c_3},
"T19"*{x_{\sigma}}
\endxy
$$

\end{proof}

\begin{theo}
\label{theonerve31}
\begin{itemize}
\item[(1)] For all $g\in \obj(o\downarrow B')$, the space $\CV_g$ is open.
\item[(2)] We have $$ \UC((W\cq V^{\reg})^{\mu_d}, o) = \bigcup_{g\in \obj(o\downarrow B')} \CV_g.$$
\item[(3)] The nerve of the open covering $(\CV_g)_{g\in  \obj(o\downarrow B')}$ is $\gar(B',S',o)$.
\item[(4)] For all simplex $\{g_0,\dots,g_k\}\in \gar(B',S',o)$, the intersection 
$\bigcap_{i} \CV_{g_i}$
is contractible.
\end{itemize}
\end{theo}

The complex $\gar$ is introduced in Definition \ref{gar2}.
Before proving the theorem, we observe that it implies Theorem \ref{theo31} and, by addressing
the remaining case of $G_{31}$, completes
the proof of the $K(\pi,1)$ conjecture.

\begin{coro}
The relative regular orbit space $(W\cq V^{\reg})^{\mu_d}$ is a $K(\pi,1)$ space.
\end{coro}

\begin{proof}[Proof of corollary]
Using \cite[4G.3]{hatcher},
we deduce from Theorem \ref{theonerve31} that $$\UC((W\cq V^{\reg})^{\mu_d}, o)$$ is homotopy
equivalent to $|\gar(B',S',o)|$. By Theorem \ref{gar3}, $|\gar(B',S',o)|$ is contractible.
\end{proof}

\begin{proof}[Proof of theorem]
(1):
consider a point in $\CV_g$ associated with $\gamma$ and $\gamma_T$, where $T=(x,L)$ is a circular semitunnel; for a neighborhood $\Omega$ of $x$ in $\CU^{\mu_d}$, concatenating
$\gamma$ with the
circular semitunnels $((x',L))_{x'\in \Omega}$ yields a neighborhood of the original point in 
$\CV_g$.

(2): consider a point $p$ in $\UC((W\cq V^{\reg})^{\mu_d}, o)$, associated with a path $\gamma$
in $(W\cq V^{\reg})^{\mu_d}$ with source in $o$. If $\gamma(1)\in \CU^{\mu_d}$,
then taking $T$ the trivial circular semitunnel of length $0$ starting at $\gamma(1)$ shows that
$p\in \CV_g$, where $g$ is the element associated with $\gamma$. In the general case, we may
find $\varepsilon>0$ arbitrarily small such that $e^{\sqrt{-1}\varepsilon} \gamma(1)$ lies in $\CU^{\mu_d}$; $p$ is then represented by $e^{\sqrt{-1}\varepsilon} \gamma$ concatenated with
the path associated with a small circular semitunnel $T_{\varepsilon}$: we see that
$p\in \CV_g$, where $g$ is the element associated with $e^{\sqrt{-1}\varepsilon} \gamma$.

(3) (rank $2$): assume that $\CV_g\cap \CV_{g'} \neq \varnothing$. Then $p \in \CV_g\cap \CV_{g'}$ can be represented by both $\gamma \gamma_T$ (where $\gamma$ represents $g$ and $T=(x,L)$ is a circular
semitunnel) and $\gamma' \gamma_{T'}$ (where $\gamma'$ represents $g'$ and $T'=(x',L')$ is a circular
semitunnel). Up to exchanging $g$ and $g'$, we may assume that $L\geq L'$. Let $T''$ be the shortened
circular semitunnel $(x,L-L')$, and $T'''$ the remaining chunk $(e^{\sqrt{-1}(L-L')}x, L')$.
The paths $\gamma\gamma_T$ and $\gamma'\gamma_{T'}$ are homotopic; so are
$\gamma_{T}$ and $\gamma_{T''}\gamma_{T'''}$. Because they have the same target
and because they both are scalar rotations, the paths $\gamma_{T'''}$ and $\gamma_{T'}$ are homotopic.
We conclude that $\gamma \gamma_{T''}$ and $\gamma'$ are homotopic. In particular, $T''$ is a circular
tunnel, representing a simple element $s''$, and we have $gs''=g'$. Conversely, it is clear
that if $s'$ is a simple element such that $gs''=g'$, we can explicitly construct a point
in $\CV_g\cap \CV_{g'}$.

(3) (higher rank): the same argument still works, after noting that a nerve simplex
$\{g_0,\dots,g_k\}$ can be ordered in such a way that each $g_i$ is represented
by $\gamma_i\gamma_{T_i=(x_i,L_i)}$ and $L_0\geq \dots \geq L_k$.
Considering for each $i \leq j$ the truncated circular
semitunnel $T''_{i,j} = (x_i,L_i-L_j)$ and remaining
chunk $T'''_{i,j} = (e^{\sqrt{-1}(L_i-L_j)}x_i, L_j)$, we see that $T''_{i,j}$ is actually
a circular tunnel, representing a simple element $s''_{i,j}$ such that $g_is''_{i,j} = g_j$.

(4) First, we consider indivual $\CV_g$. By retracting the circular
semitunnel part to length $0$,
we easily see that each $\CV_g$ is contractible (this is the analog of
Lemma \ref{lemmacon}).

Now we have to prove the analog of Proposition \ref{intersections}, and consider non-empty
non-trivial intersections $\bigcap_{j=0}^n \CV_{g_j}$.
The proof of (3) provides the basis for an explicit description of
these. We keep the same conventions and notations.
A point $p$ in $\bigcap_{j=0}^n \CV_{g_j}$ can be represented by
$\gamma_0 \gamma_{T}$, where $T=(x,L)$ is a circular semitunnel such that each truncation
$T''_{0,j}= (x_0,L_0 - L_j)$ is a circular tunnel representing
$s''_{0,j} := g_0^{-1}g_j$.
Interestingly, changing $\gamma_0$ to another path representing $g_0$ yields the same $p$,
because there is only one homotopy class of such paths.
In other words, $\bigcap_{j=0}^n \CV_{g_j}$ is indexed by circular
semitunnels $T=(x,L)$ such that there exists a non-decreasing sequence
$0  \leq L'_1 \leq \dots L'_k \leq L$ such that, for each $j$, the truncation
$(x,L'_j)$ represents $g_0^{-1}g_j$.
When $T$ satisfies this condition, we say that it \emph{represents} the sequence
$(g_0^{-1}g_1,\dots,g_0^{-1}g_k)$.
Actually, $\bigcap_{j=0}^n \CV_{g_j}$ is homotopy equivalent to the space of
circular semitunnels representing $(g_0^{-1}g_1,\dots,g_0^{-1}g_k)$, equipped with 
the product topology (this is the analog of Proposition \ref{inequalities}).

We say that $x\in \CU^{\mu_d}$ represents $(g_0^{-1}g_1,\dots,g_0^{-1}g_k)$
if there exists a non-decreasing sequence
$0  \leq L'_1 \leq \dots L'_k \leq \frac{2\pi}{d'h}$ such that, for each $j$, the truncation
$(x,L'_j)$ represents $g_0^{-1}g_j$.
An obvious retraction argument (onto the particular choice $L=\frac{2\pi}{d'h}$ for the length)
shows that the space of circular semitunnels
representing $(g_0^{-1}g_1,\dots,g_0^{-1}g_k)$ is homotopy-equivalent to
the space of points $x\in \CU^{\mu_d}$ representing
$(g_0^{-1}g_1,\dots,g_0^{-1}g_k)$.
 
To conclude, we are down to proving that the space of points $x\in \CU^{\mu_d}$ representing
$$(g_0^{-1}g_1,\dots,g_0^{-1}g_k)$$ is contractible. This is Lemma \ref{lemmareprcontr}.
\end{proof}

\begin{conj}
The methods of this section can be adapted to work in the context of Lehrer-Springer
theory (an extension of Springer theory, see \cite{lesp}.) The combinatorics of the associated
categorical Garside structure follow Armstrong's variant of the cyclic
sieving phenomenon (see \cite[Theorem 2, p. 204]{krmu}.)
\end{conj}

\section{Periodic elements in $B(W)$}
\label{section12}

As before,
$W$ is an irreducible well-generated
complex reflection group, 
$\tau\in P(W)$ is the full-twist and
$\delta\in B(W)$ is the Garside element of the dual braid monoid $M(W)$.
The image
of $\delta$ in $W$ is a Coxeter element $c$.

\begin{defi}
An element of $B(W)$ is
\emph{periodic} if it admits a central power.
\end{defi}

The goal of this section is to prove that the center $ZB(W)$ is
cyclic, and to establish a correspondence between periodic elements
in $B(W)$ and regular elements in $W$.
As with the previous section, the real substance of the arguments
lies more
in the algebraic tools from \cite{cyclic} than in the easy geometric
interpretation.

\begin{lemma}
\label{lemma1c}
The intersection of the subgroup $\left< c\right> \subseteq W$ with
the interval $[1,c]$ is $\{1,c\}$.
\end{lemma}

\begin{proof}
\cite[Proposition 4.2]{dimi} gives an argument for the real case that,
as pointed out by J. Michel, applies verbatim. We include
it for the convenience of the reader.
Let $2\leq d_1, \dots,d_n=h$ be the reflection degrees.
Because $c$ is regular, its eigenvalues are $\zeta_d^{1-d_i}$ 
(\cite[Theorem 4.2 (v)]{springer}).
Assume that some power $c^k$ lies in $[1,c]$, thus that
$$l_R(c^k)+l_R(c^{1-k})=n.$$
By Proposition \ref{propiiiii} and Lemma \ref{lemmaparacoxnew},
$n - l_R(c^k)$ is the number of eigenvalues of $c^k$ distinct from $1$, thus
$$l_R(c^k) = n- \#\{i | (1-d_i)k \equiv 0[h]\}.$$
Using the same formula for $c^{1-k}$, we deduce that
$$\#\{i | (1-d_i)k \equiv 0[h]\} + \#\{i | (1-d_i)(1-k) \equiv 0[h]\}=n.$$
For a given $i$, $(1-d_i)k$ and $(1-d_i)(1-k)$ cannot simultaneously
divide $h$. The identity then forces that, for all $i$, either
$(1-d_i)k$ or $(1-d_i)(1-k)$ is a multiple of $h$.
But, when $2\leq k \leq h-1$, 
neither $(1-h)k$ nor $(1-h)(1-k)$ is a multiple of $h$.
\end{proof}

\begin{theo}
\label{theoZ}
Let $h':= \frac{h}{\gcd (d_1,\dots,d_n)}$.
The centers of $B(W)$ and $W$ are cyclic, generated
respectively by $\delta^{h'}$ and $c^{h'}$. 
\end{theo}

\begin{proof}
That $ZW$ is cyclic of order $\gcd (d_1,\dots,d_n)$
is classical; any central element is regular;
because $c^{h'}$ is regular and has the right order, it must generate
$ZW$.

Let us study the conjugacy action of $\delta$.
Write a given $b\in B(W)$ in Garside normal form $b=\delta^k s_1\dots s_l$
where $k\in \BZ$ and $s_1,\dots,s_l \in [1,c] - \{1,c\}$.
The normal form of $b^{\delta}$ is $\delta^k s_1^{c}\dots s_l^{c}$.

Because $c^{h'}$ is the smallest central power of $c$ in $W$,
$\delta^{h'}$ is the smallest central power of $\delta$ in $B(W)$.
Moreover, if $b\in ZB(W)$, then it must commute with $\delta$ and
any simple term $s_i$ in its normal form $\delta^k s_1\dots s_l$
must commute with $c$.
Because $c$ is
Coxeter element, the centralizer of $c$ in $W$ is $\left< c\right>$
(this follows from Theorem \ref{theoregular} and, actually,
extends Corollary 4.4 in  \cite{springer}). Combining this with
 Lemma \ref{lemma1c}, we see that $l=0$, \ie, that
$b \in \left< \delta \right>$.
\end{proof}

Because of Theorem \ref{theoZ}, an element $\gamma\in B(W)$
is periodic if and only if it is \emph{commensurable} with $\tau$
(or $\delta$), \ie,
if there exists $p,q$ such that $$\gamma^q=\tau^p.$$
To simplify notations, we restrict our attention to the situation where
$$\gamma^d= \tau$$
and call such a periodic $\gamma$ a \emph{$d$-th root of $\tau$}.
The theory works the same way for other $(p,q)$.

When $d$ is regular and $x_0$ is a $\zeta_d$-regular eigenvector,
we may consider the \emph{standard $d$-th root of $\tau$} represented
by 
\begin{eqnarray*}
[0,1] & \longrightarrow  & W\cq V^{\reg}\\
t & \longmapsto &e^{2\pi\sqrt{-1}\frac{t}{d}} x_0
\end{eqnarray*}
and denoted by $\sqrt[d]{\tau}$.
This of course involves choosing a particular basepoint, but because
the statements below are ``up to conjugacy'', one should not worry too
much about this.

A particular case is $$\sqrt[h]{\tau}= \delta.$$

\begin{theo}[Springer theory in braid groups]
\index{Springer theory!braid group version}
\label{theoroots}
Let $d$ be a positive integer.
\begin{itemize}
\item[(i)] There exists $d$-th roots of $\tau$
if and only
if $d$ is regular.
\item[(ii)] When $d$ is regular,
there is a single conjugacy class of $d$-th roots of $\tau$ in $B(W)$.
In particular, all $d$-th roots of $\tau$ are
conjugate to $\sqrt[d]{\tau}$.
\item[(iii)] Let $\rho$ be a $d$-th root of $\tau$. 
Let $w$ be the image of $\rho$ in $W$. Then $w$ is
$\zeta_d$-regular, and the centralizer
$C_{B(W)}(\rho)$ is isomorphic to the braid group $B(W')$ of the centralizer
$W':=C_{W}(w)$.
\end{itemize}
\end{theo}

\begin{proof}
(i) A consequence of \cite[Corollary 10.4]{cyclic} is that, if
$\tau = \delta^h$ admits $d$-th roots, then $M_d^{\phi_d^h}$ is
non-empty. By Lemma \ref{lemmadinvariant}, this implies that $d$ is regular.
The converse is obvious (we may consider the particular
root $\sqrt[d]{\tau}$).

(ii) Using \cite[Corollary 10.4]{cyclic} 
and \cite[Proposition 9.8]{cyclic}, we see that
conjugacy classes of $d$-th roots of $\tau$ are in one-to-one
correspondence with connected components of the category
$G'$. Because it is equivalent to a group, this category is connected.

(iii) That $w$ is regular follows from (ii), because
it is conjugate to the image in $W$ of $\sqrt[d]{\tau}$,
whose image in $W$ has a $x_0$ as $\zeta_d$-regular
eigenvector. The assertion about the centralizer follows from
its categorical rephrasing in $M_d$, where it is trivial
(the conjugacy action being a power of the diagram automorphism
of the Garside structure).
\end{proof}

\begin{remark}
This answers many questions and conjectures by Brou\'e, Michel and others
(see \cite{boston} for more details). Particular cases of
(i) were obtained by Brou\'e-Michel and, independently, by Shvartsman,
\cite{brmi,shvartsman}.
Assertion (ii) can be viewed either as a Ker\'ekj\'art\'o type theorem
(``\emph{all periodic elements are conjugate to a rotation}'', see
\cite{cyclic}) or as a Sylow type theorem (``\emph{all $C_{B(W)}(\rho)$
are conjugate}'').
The type $A$ case of (ii) actually follows from Ker\'ekj\'art\'o's theorem
on periodic homeomorphism of the disk (\cite{kerekjarto}, see also \cite{cyclic} for a more
complete bibliography). The type $A$ case of (iii)
was proved in \cite{BDM}.

Of course, the most natural interpretation is
to view the theorem as providing a braid analog of Theorem \ref{theoregular}.
\end{remark}

\begin{remark}
Let $W'$ be the centralizer of a regular element in $W$. Let $W''$
be the centralizer of a regular element in $W'$.
In terms of orbit varieties, $W'\cq V'= (W\cq V)^{\mu_d}$ and
$W''\cq V''= (W\cq V)^{\mu_{de}}$. Regular
elements of $W'$ are regular in $W$. It should be possible,
by applying Theorem 
\ref{theoroots} to the pair $(W,W'')$, to generalize
the result to the pair $(W',W'')$.
\end{remark}

\begin{coro}
\label{coroZ31}
The center of the braid group $B(G_{31})$ is cyclic.
\end{coro}

\begin{proof}
Let $W$ be $G_{37}$, the well-generated reflection group
of type $E_8$. The degrees
are $$2, 8, 12, 14, 18, 20, 24, 30$$
and the codegrees are $$0, 6, 10, 12, 16, 18, 22, 28.$$

The number $4$ is regular, with centralizer $W'$ of type $G_{31}$.

By Theorem \ref{theoregular} (1), we see that $24$ is also regular.
Let $\rho$ be a $24$-th root of $\tau$. The centralizer is the rank
$1$ reflection group of type
$A_1$, with braid group $\mathbf{Z}$.

Applying Theorem \ref{theoroots} to $\rho^6$, we recognize
the braid group of $G_{31}$ as a centralizer in $B(W)$:
$$B(G_{31}) \simeq C_{B(W)}(\rho^6).$$

Applying Theorem \ref{theoroots} to $\rho$, we see that
$$C_{B(W)}(\rho)\simeq B(A_1) \simeq \mathbf{Z}.$$

Clearly, $\rho\in  C_{B(W)}(\rho^6)$.
In particular, any $z\in ZB(G_{31})$ must commute with $\rho$,
hence lie in $C_{B(W)}(\rho)\simeq \mathbf{Z}.$
\end{proof}

Combining Theorem \ref{theoZ} and Corollary \ref{coroZ31},
we complete the proof of:

\begin{theo}[Theorem \ref{theocenter}]
\label{theoZZ}
The center of the braid group of an irreducible complex reflection
group is cyclic.
\end{theo}

Indeed, Brou\'e-Malle-Rouquier conjectured this in \cite{bmr},
and proved it for all cases but 6 exceptional types: 5 of these exceptions
are well-generated and covered by Theorem \ref{theoZ}, the
remaining case being $G_{31}$.

\section{Generalized non-crossing partitions}

Here again, $W$ is an irreducible well-generated complex reflection 
group generated by reflections of order $2$.

When $W$ is of type $A_{n-1}$, the lattice $(S,\preccurlyeq)$ is isomorphic
to the lattice of non-crossing partitions of a regular $n$-gon
(\cite{BDM}, \cite{brady}).
Following \cite{reiner}, \cite{dualmonoid} and \cite{BC}, we call
\emph{lattice of generalized non-crossing partitions of type $W$}
the lattice $$(S,\preccurlyeq),$$
 and \emph{Catalan number of type $W$} 
\index{well-generated reflection groups!generalized Catalan numbers}
the number
$$\Cat(W):= \prod_{i=1}^n \frac{d_i+h}{d_i}.$$

The operation sending $s\preccurlyeq t$ to $s^{-1}t$ is an analogue
of the Kreweras complement operation. The map $s\mapsto s^{-1}\delta$ is
an anti-automorphism of the lattice.

In the Coxeter case,
Chapoton (see \cite{chapoton}) discovered a general formula for the number
$Z_W(N)$ of weak chains $s_1\preccurlyeq s_2 \preccurlyeq \dots
\preccurlyeq s_{N-1}$ of length $N-1$
in $(S,\preccurlyeq)$, or equivalently (via Lemma \ref{transferlemma}
for the cardinality of $\CD_N(c)$.
This formula continues to hold, though we are only able to prove this
case-by-case
(see \cite{athanasiadisreiner} and \cite{chapoton} for the Coxeter types;
the $G(e,e,n)$ case was done in \cite{BC};
the remaining types are done by
computer).

\begin{prop}
We have, for all $N$, 
$$Z_{W}(N)=\prod_{i=1}^n \frac{d_i+ (N-1)h}{d_i}.$$
\end{prop}

\begin{coro}
\index{simple elements!cardinality}
We have $|S| = \Cat(W)$.
\end{coro}

Another interesting numerical invariant is the \emph{Poincar\'e polynomial}
$$\Poin(S) := \sum_{s\in S} t^{l(s)}.$$

The numerical data for the exceptional types
(real and non-real) is summarised in Table \ref{table2}. The coefficient of $t$ in the Poincar\'e polynomial
is the cardinal of $R_c$. One observes that $$R=R_c \Leftrightarrow W
\text{ is real.}$$

In the Weyl group case, $\Poin(S)$ may be interpreted as the Poincar\'e
polynomial of the cohomology of a toric variety related to cluster algebras
(\cite{chapoton}).

\begin{table}
\begin{tabular}{|c|c|c|c|c|c|}
\hline
$W$ & degrees & $|R|$ & $\Cat(W)$ & $\Poin(S)$ & $|\Red_R(c)|$ \\
\hline
$G_{23}$ ($H_3$) & $2,6,10$  & $15$ & $32$ & $1 + 15t + 15t^2 + t^3$ & $50$ \\
\hline
$G_{24}$ & $4,6,14$  & $21$ & $30$ & $1+14t+14t^2+t^3$ & $49$ \\
\hline
$G_{27}$ & $6,12,30$  & $45$ & $42$ & $1+20t+20t^2+t^3$ & $75$ \\
\hline
$G_{28}$ ($F_4$) & $2,6,8,12$ &  $24$ & $105$ & 
$1+24t+55t^2+24t^3+t^4$ & $432$ \\ 
\hline
$G_{29}$ & $4,8,12,20$ &  $40$ & $112$ & $1+25t+60t^2+25t^3+t^4$
& $500$ \\
\hline
$G_{30}$ ($H_4$) & $2,12,20,30$  & $60$ & $280$ &
$1+60t+158t^2+60t^3+t^4$ & $1350$\\
\hline
$G_{33}$ & $4,6,10,12,18$ &  $45$ & $308$ &
$\begin{matrix}1+30t+123t^2 \\+123t^3+30t^4+t^5\end{matrix}$ & $4374$ \\
\hline
$G_{34}$ & $\begin{matrix} 6,12,18,\\ 24,30,42\end{matrix}$ &
$126$  &  $1584$ &
$\begin{matrix}1+56t+385t^2+ 700t^3 \\
+385t^4+56t^5+t^6\end{matrix}$ & $100842$ \\
\hline
$G_{35}$ ($E_6$) & $\begin{matrix}2,5,6, \\8,9,12\end{matrix}$ & 
 $36$  &  $833$ & 
$\begin{matrix}1+36t+204t^2+351t^3\\ +204t^4+36t^5+t^6\end{matrix}$ & 
$41472$ \\
\hline
$G_{36}$($E_7$) & $\begin{matrix} 2,6,8,10,\\ 12,14,18\end{matrix}$ & $63$ &
$4160$ &
$\begin{matrix}1+ 63t+ 546t^2+ 1470t^3\\
+ 1470t^4+ 546t^5+ 63t^6 + t^7 \end{matrix}$ & 
$1062882$ \\
\hline
$G_{37}$ ($E_8$)& $\begin{matrix}2,8,12,14, \\18,20,24,30 \end{matrix}$ &
$120$
& $25080$ &
$\begin{matrix} 1+120t+1540t^2 \\
+6120t^3+9518t^4+6120t^5 \\ +1540t^6+120t^7+t^8\end{matrix}$ &
$37968750$ \\
\hline
\end{tabular}
\smallskip
\caption{Numerical invariants of generalized non-crossing partitions.}
\label{table2}
\end{table}

When $W'$ is the (not necessarily well-generated) centralizer
of a $d$-regular element in a well-generated $W$, the natural
substitute for $Z_W(N)$
is the cardinality $Z'_{W'}(N)$ of $\CD'_N(c)$.
In a joint work with Vic Reiner,
we conjectured
that the $\mu_d$-action exhibits a \emph{cyclic sieving phenomenon}:
the number of fixed points should be the value at $\zeta_d$ of
a $q$-analog of the number of chains; the conjecture has now been proved
by Krattenthaler-M\"uller:

\begin{theo}[Conjecture 6.5 in \cite{BR}, proved in \cite{krmu}]
Let $q$ be an indeterminate. For any $a\in \BZ_{\geq 1}$, set
$[a]_q:=1+q+ \dots+q^{a-1}$.
Then $$\prod_{i=1}^n \frac{[d_i+(Nd'-1)h]_q}{[d_i]_q}$$
is a polynomial in $q$ whose value at $q=\zeta_d$ is
$Z'_{W'}(N)$.
\end{theo}

When $N=1$, the formula gives the number of objects in $M'$.
For $G_{31}$, we have $88$ objects (see Example \ref{exampleG31category}).

\appendix

\section{The fat basepoint trick}
\label{appendixA}

\subsection{Fundamental groupoids}
Let $E$ be a topological space. Let $\gamma$ be a path in $E$, \ie,
a continuous map $[0,1]\rightarrow E$. We say that $\gamma$ is a \emph{path
from $\gamma(0)$ to $\gamma(1)$} -- or that $\gamma(0)$ is the \emph{source}
and $\gamma(1)$ the \emph{target}. The concatenation rule is as follows:
if $\gamma,\gamma'$ are paths such
that $\gamma(1)=\gamma'(0)$,
the product $\gamma  \gamma'$ is the path mapping
$t\leq 1/2$ to $\gamma(2t)$ and $t\geq 1/2$ to $\gamma'(2t-1)$.

We denote by $\pi_1(E)$
\index{fundamental groupoid}
the fundamental groupoid of $E$: its elements are homotopy classes of
paths in $E$, with composition rule as above. As a category, its
object set is $E$.

For any ``basepoint'' $e\in E$, the \emph{fundamental
group of $E$ with respect to $e$} is

$$\pi_1(E,e) := \Hom_{\pi_1(E)}(e,e).$$

When $E$ is a \emph{path-connected} space, or equivalently when
$\pi_1(E)$ is a \emph{connected} groupoid, all fundamental groups
of $E$, with respect to all possible basepoints, are isomorphic.
\textbf{However, they are not canonically isomorphic:}
any element $g\in \Hom_{\pi_1(E)}(e,e')$ yields
an isomorphism:
\begin{eqnarray*}
\phi_g:\pi_1(E,e) & \stackrel{\sim}{\longrightarrow} & \pi_1(E,e')\\
f & \longmapsto & g^{-1} f g
\end{eqnarray*}
but there is no natural way to make consistent choices
of such $g$'s and build a \emph{transitive systems of
isomorphisms} connecting $$(\pi_1(E,e))_{e\in E}.$$ In other
words, there is no legitimate
way to drop the reference to a specific basepoint and talk about
\emph{the} fundamental group of $E$.

Let $A\subseteq E$. Consider the natural functor
$\iota:\pi_1(A) \to \pi_1(E)$.

\begin{lemma}
When $A$ is simply-connected, then $$(\phi_{\iota(g)})_{g\in \pi_1(A)}$$
is a transitive system of isomorphisms connecting
$(\pi_1(E,a))_{a\in A}$.
\end{lemma}

\begin{proof}
The space $A$ is simply-connected if and only if the category
$\pi_1(A)$ is equivalent to the trivial category.
In a category equivalent to the trivial category, transitivity
comes for free: $(\phi_g)_{g\in \pi_1(A)}$ is transitive,
and by functoriality so is $(\phi_{\iota(g)})_{g\in \pi_1(A)}$.
\end{proof}

This legitimates the following definition:

\begin{defi}[fat basepoint trick, group version]
\index{fat basepoint trick}
The fundamental group of $E$ with respect to a simply-connected
subspace $A$ is the transitive limit
$$\pi_1(E,A) = \lim_{\stackrel{\longrightarrow}{a\in A}} \pi_1(E,a)$$
with respect to the transitive system of isomorphisms 
$(\phi_{\iota(g)})_{g\in \pi_1(A)}$.
\end{defi}

\begin{remark}
Instead of constructing the group $\pi_1(E,A)$ as
a transitive limit of fundamental groups, one can choose to equip the set of relative
homotopy classes with a group structure. No difference, except in language.
\end{remark}

Practically speaking, $\pi_1(E,A)$ should be thought of
as any $\pi_1(E,a)$, for some $a\in A$, together with an unambiguous
recipe thanks to which any path in $A$
from any $a'\in A$ to any $a''\in A$ represents a unique element
of $\pi_1(E,a)$ -- moreover, one may
forget about which $a\in A$ was chosen and change it at our convenience.

In Section \ref{section11}, we use an extended version of the trick.
Let $C$ and $C'$ be path-connected components of $A$.
For any $g\in \Hom_{C}(c_1,c_2)$ and $g'\in \Hom_{C'}(c'_1,c'_2)$,
we have an isomorphism
\begin{eqnarray*}
\psi(g,g'): \pi_1(E,c_1,c'_1) & \stackrel{\sim}{\longrightarrow} &
\pi_1(E,c_2,c'_2)\\
f & \longmapsto & g^{-1} f g'
\end{eqnarray*}
When both $C$ and $C'$ are simply connected, this provides
a transitive system of isomorphisms $(\psi(g,g'))_{(g,g')\in \pi_1(C)
\times \pi_1(C')}$ thanks to which we can define
$$\Hom_E(C,C') := \lim_{\stackrel{\longrightarrow}{(c,c')\in
C\times C'}} \Hom_E(c,c')$$

\begin{defi}[fat basepoint trick, groupoid version]
\label{fatgroupoid}
The fundamental groupoid of $E$ with respect to a subspace
$A$ whose path-connected components are simply connected
is the groupoid
$$\pi_1(E,A)$$
whose objects are path-connected components of $A$ and such
that 
$$\Hom_{\pi_1(E,A)}(C,C') = \Hom_E(C,C').$$
\end{defi}

When $A$ is connected, we recover the group version.

\subsection{Universal covers}
\label{subunivoc}
Assume that $E$ is path connected.
By Galois theory, there is a
correspondence between subgroups of the fundamental
group and topological coverings of $E$.

Universal covers can be constructed as soon
as $E$ is reasonably healthy:
e.g., it is enough to assume that $E$ is locally simply-connected
(which is trivially satisfied by all spaces considered here).

A good reference is Hatcher, \cite[Section 1.3]{hatcher}.
The construction starts with the choice of a basepoint $e\in E$.
As a set, the universal cover $\widetilde{E}_e$ has one point
per element in $\pi_1(E)$ with source $e$.
The fundamental groupoid $\pi_1(\widetilde{E}_e)$
coincides with the category $(e\downarrow \pi_1(E))$ of \emph{objects under $e$ in $\pi_1(E)$},
in the sense of Mac Lane, \cite[II.6, Comma Categories]{maclane}.
This interpretation actually clarifies why $\widetilde{E}_e$
is simply connected: the category $(e\downarrow \pi_1(E))$
is equivalent to the trivial category (this is the categorical way of being contractible).

The universal cover construction can thus be rephrased as follows:
the object set of $(e\downarrow \pi_1(E))$
can be equipped with a natural topology (this is where
the locally simple-connectedness of $E$ comes into play).

Actually, when $\CG$ is a groupoid and $o$ is an object of $\CG$, then $(o\downarrow \CG)$
should be viewed as a
categorical universal cover for $\CG$.

When $A\subseteq E$ is simply-connected, there is
a natural transitive system of isomorphisms connecting
$$\left( (a\downarrow \pi_1(E)) \right)_{a\in A}$$
and thus a transitive system of bijections connecting
$$(\widetilde{E}_a)_{a\in A}.$$
It is not hard to check that these bijections are homeomorphisms
and are compatible with the covering maps.

\begin{defi}[fat basepoint trick, universal cover version]
\label{fatcover}
\index{$\UC$ (universal cover model)}
The universal cover of $E$ with respect to a simply-connected
subspace $A$ is
$$\UC(E,A) := \lim_{\stackrel{\longrightarrow}{a\in A}} \widetilde{E}_a.$$
\end{defi}

Clearly, $\UC(E,a) \simeq \widetilde{E}_a$.

\begin{lemma}
The natural left-action of $\pi_1(E,a)$ on $\widetilde{E}_a$
is compatible with the transitive system and gives rise
to a left-action of $\pi_1(E,A)$ on $\UC(E,A)$.
\end{lemma}

In real life, this means: up to topologically trivial paths within $A$,
any path connecting two points of $A$
(representing an element of $\pi_1(E,A)$) can be 
concatenated with any path with source in $A$
(representing a point in $\UC(E,A)$) to
unambiguously yield another
point in $\UC(E,A)$. We don't care about which exact points
in $A$ were chosen,
we don't have to give them names and we can forget
about them.

When $A$ isn't contractible but has simply connected components,
the \emph{groupoid cover} version, blending Definitions \ref{fatgroupoid} and \ref{fatcover},
is a bit more tedious to formulate:
\begin{itemize}
\item we get a model $\UC(E,o)$ for each connected component $o\subseteq A$;
\item any $g$ in $\pi_1(E,A)$ with source $o$ and target $o'$ induces an homeomorphism
from $$\UC(E,o) \to \UC(E,o')$$ subject to obvious compatibility rules.
\end{itemize}

The latter data is what deserved to be called
a \emph{groupoid action on} $(\UC(E,o))_{o\in \pi_0(A)}$.
This situation is implicit behind Definition \ref{boo}.

\section{Garside structures}
\label{garsidependix}

\subsection{Cohomology of groups and groupoids}

A groupoid is (small) category where all morphisms are invertible. A group is a groupoid with a single object.

A \emph{simplicial complex} $X$ is a family of subsets of an ambient space $S$, such that
whenever $A\in X$, any subset $B\subseteq A$ also lies in $X$. Simplicial complexes form
a category, equipped with a \emph{geometric realization} functor to the category of
topological spaces.

As explained in \cite{maclane}, a
\emph{simplicial set}
is a contravariant
functor from the simplicial category to the category of sets,
or equivalently a collection $(X_n)_{n\in \mathbf{Z}_{\geq 0}}$ of sets,
together with \emph{face maps} and \emph{degeneracy maps} respectively shifting dimensions by $-1$ and $+1$, and subject to obvious compatibility rules.

Simplicial complexes naturally give rises to simplicial sets, but not all simplicial sets
can be obtained that way.

The \emph{nerve} \index{nerve!of a category}
of a (small) category $\CC$ is a simplicial
set $\CN\CC$ whose $0$-skeleton is the object set of $\CC$
and whose $n$-simplices, $n\geq 1$, are composable sequences
$(f_1,\dots,f_n)$ of 
$\CC$-morphisms:
$$\xymatrix{
x_0 \ar[r]^{f_1} & x_1 \ar[r]^{f_2} & x_2 \ar[r]^{f_3}  & x_3  \ar@{..>}[r]
& x_{n-1} \ar[r]^{f_{n}} & x_{n}}.$$
In the simplicial structure on $\CN\CC$, face maps correspond to removing objects (and composing or dropping
morphisms
accordingly) and degeneracy
maps correspond to inserting identity morphisms at a given object.
An element of $\CN\CC$ is \emph{non-degenerate} if it does not contain
any identity morphism.

Simplicial sets naturally form a category, and the nerve construction is functorial from the category $\mathbf{Cat}$ of small categories to the category $\mathbf{SimpSet}$ of simplicial sets. There is a standard \emph{geometric realization} functor $\mathbf{SimpSet}\to \mathbf{Top}, X\mapsto |X|$, where $\mathbf{Top}$ category of topological spaces.
See for example \cite[Appendix, Simplicial CW structures]{hatcher}.
The construction actually provides us, for each 
\emph{abstract} simplex $x\in X$, with a singular simplex in $|X|$, \ie
a continuous from a standard affine simplex to $X$.
Actually, the $0$-skeleton $X_0$ is mapped injectively into $|X|$.
We can use $X_0$ as a groupoid fat basepoint (Definition \ref{fatgroupoid}).

\begin{defi}
Let $\CG$ be a groupoid.
A \emph{simplicial $K(\CG,1)$} is a simplicial set $X$ such that
$\pi_1(|X|, X_0)\simeq \CG$, and such that connected components
of $|X|$ have no higher homotopy groups.
\end{defi}

In particular, $X_0$ must be in bijection with
the object set of $\CG$.

\begin{theo}[``bar'' resolution, quotient version]
\label{bar1}
The nerve $\CN\CG$ of a groupoid $\CG$ is a simplicial $K(\CG,1)$.
\end{theo}

The beauty of the categorical viewpoint is that the theorem comes as
a mostly free by-product of the observation that the
\emph{nerve realization} functor $\mathbf{Cat} \to \mathbf{Top}$
extends to a functor of $2$-categories (mapping natural transformations
to homotopies). See Lemma 7.1 and Proposition 7.3 in
\cite{cyclic}.
The underlying combinatorics coincide with that of the
standard ``bar'' resolution of group cohomology:

\begin{defi}[``bar'' simplicial set]
\label{bar2}
Let $o$ be an object of $\CG$. Let $(g_0,g_1,\dots,g_k)$ be a sequence of composable morphisms
in $\CG$ (\ie, an element of $\CN_{k+1}\CG$) such that the source of $g_0$ is $o$.
We use the bar symbol $g_0[g_1|\dots|g_k]$ to express that $(g_0,g_1,\dots,g_k)$ is a tuple
satisfying these properties.
We denote by $$\bar(\CG, o)$$ the simplicial set
whose $k$-skeleton consists of bar symbols $g_0[g_1|\dots|g_k]$, subject to the faces and degeneracy maps of $\CN_{k+1}\CG$ (except those involving $g_0$, as we view $g_0[g_1|\dots|g_k]$ as 
a $k$-simplex, not a $(k+1)$-simplex).
\end{defi}

\begin{theo}[``bar'' resolution, universal cover version]
\label{bar3}
Let $\CG$ be a groupoid, let $o$ be an object of $\CG$. The geometric realization of
$\bar(\CG, o)$ is contractible.
\end{theo}

Theorems \ref{bar1} and \ref{bar3} express two flavors of the same result: $\UC(|\CN\CG|,o)
\simeq |\bar(\CG, o)|$. Note that we do not have to assume that $\CG$ is connected:
Theorem \ref{bar1} expresses something about each connected component, 
while Theorem \ref{bar3} sees only one connected component per choice of $o$.

\subsection{Garside structures}
Garside's approach \cite{garside}
to the word and conjugacy problem in the classical braid group $B_n$
was a key ingredient in Deligne's paper \cite{deligne}.
It was later
axiomatized as a generic combinatorial group theory
notion \cite{dehpa}, rephrased with a geometric group theory
viewpoint \cite{bestvina,cmw} and generalized to groupoids
\cite{krammer}.

There are many ways to tell the story, and the most general
setup involves quite a lot of technicalities.
Under favorable conditions, typically when
there is a natural homogeneous length function (which is the
case here), the whole story could
probably fit in a 100 pages graduate-level textbook.
In this absence of this yet-to-be-written account,
the only detailed reference at hand is
the much longer book \cite{ddgm}, that focuses on word-theoretic
axiomatic aspects and does not cover all aspects explained here.

As far as the current paper is concerned,
the notations and results listed in \cite{cyclic} are more than sufficient.
Especially, we only consider Garside structure that are
homogeneous, which removes a lot of the technicalities addressed in \cite{ddgm}.

Let $\CC$
be a (small) category equipped with
\begin{itemize}
\item an endofunctor
$\phi$, which we write with \emph{right-conjugacy} notation $x\mapsto x^{\phi}$, $f\mapsto f^{\phi}$,
\item a natural transformation $\Delta$ from the identity functor to $\phi$.
\end{itemize}

\begin{example}
When $\CC$ is a monoid $M$, this simply means that $\Delta\in M$ and
$\phi$ is the right conjugacy action $f\mapsto f^{\phi} =
\Delta^{-1} f\Delta$.
\end{example}

In general, there is a morphism $\Delta_x:x\to x^{\phi}$
for each object $x$,
and for each morphism $f$ from $x$ to $y$
the following diagram is commutative:
$$\xymatrix{ x \ar[r]^{\Delta_x} \ar[d]_{f} &
x^{\phi} \ar[d]^{f^{\phi}} \\
y \ar[r]_{\Delta_{y}} & y^{\phi} }$$

As $\Delta$ is a collection of morphisms, one for each source
object, it makes sense to write ``$f\in \Delta$'' instead of ``$f$
is the morphism in $\Delta$ whose source is the source of $f$''.
It is even tempting to write ``$f = \Delta$'' instead,
an abusive yet convenient notation inspired by
the monoid case (where $\Delta$ consists of a single element), just like
we write ``$f = 1$'' instead of ``$f$
is the identity morphism whose source is the source of $f$''.

\begin{defi}
An element of $f\in\CC$ is \emph{simple} with respect to $\Delta$ if 
there exists $g\in \CC$ such that $fg= \Delta$.

An \emph{atom} is an element $a\in \CC$ such that, for all $f,g\in\CC$, $a=fg \Rightarrow
f=1 \text{ or } g=1$.

The category $\CC$ is \emph{homogeneous} if there exists a functor $l$ from $\CC$
to the monoid $(\BZ_{\geq 0}, +)$ such that $\forall f\in \CC, l(f) = 0 \Rightarrow f=1$.

The category $\CC$ is \emph{cancellative} if $\forall f,g,h\in \CC, (fh=gh \text{ or } hg=hf)
\Rightarrow f=g$.

The category $\CC$ is a \emph{lattice} if it admits pullbacks and pushouts.
\end{defi}

Pushouts and pullbacks are classical concepts from category theory (see
for example \cite{maclane}). The existence of pushouts means that any two morphisms
$f,g$ with common source admit a right least common multiple. In poset language,
this means that $f$ and $g$ admit a least upper bound for the prefix ordering.
Pullback is the dual concept.

\begin{defi}[Definition 2.4 in \cite{cyclic}]
A \emph{Garside structure} is a triple $(\CC,\Delta,\phi)$
satisfying:
\begin{itemize}
\item[(i)] $\CC$ is a category, $\phi$ an automorphism of $\CC$ and
$\Delta$ a natural transformation from the identity functor to $\phi$,
\item[(ii)] $\CC$ is homogeneous and cancellative,
\item[(iii)] all atoms are simple with respect to $\Delta$,
\item[(iv)] $\CC$ is a lattice.
\end{itemize}
\end{defi}

This axiom set is more restrictive than necessary, but was chosen in \cite{cyclic}
to be on the safe side when claiming that existing proofs in the context of Garside monoids still
worked in the categorical contexts.

The crux is axiom (iv). This is where concrete examples of Garside structures
encapsulate deep geometric/topological/combinatorial miracles
(such as Lemma \ref{lemmalattice} or, in Deligne's paper, properties of galleries
that are specific to \emph{simplicial} arrangements).

Let $(\CC,\Delta,\phi)$ be a Garside structure. The set $\CS$ of simple
elements with respect to $\Delta$, seen either as an abstract set together
with a partial product structure obtained by restricting the category structure to $\CS$ 
(a Garside \emph{germ}, in the sense of \cite{cyclic}), or as a subset of $\CC$
(a Garside \emph{family}, in the sense of \cite{ddgm}), is enough to recover the
whole structure $(\CC,\Delta,\phi)$.

As explained in \cite{cyclic, ddgm}, the axiom set can be rewritten in terms
of axioms involving only $\CS$. Another interesting invariant of Garside structures
is the Garside set $\CD_{\bullet}$ introduced below.

\begin{defi}
A \emph{Garside category} is a category $\CC$ that can be equipped with a Garside structure.
A \emph{Garside groupoid} is a groupoid $\CG$ that is the groupoid of fractions of a
Garside category.
\end{defi}

Note that Garside groupoids can appear as groupoids of fractions of several
Garside categories that are not equivalent (for example, finite type Artin groups
admit both the classical and dual Garside structures).

Let $\CC$ be a Garside category, with groupoid of fractions $\CG$. Key properties include:
\begin{itemize}
\item the natural function $\CC \to \CG$ is faithful,
\item the word and conjugacy problems in $\CC$ and $\CG$ can be solved,
\item $\CG$ has finite cohomological dimension (Theorem \ref{gar1} below),
\item in particular, $\CG$ is torsion-free.
\end{itemize}

\subsection{Cohomology of Garside groups and groupoids}

\begin{defi}[Garside nerve]
\label{defigarnerve}
A simplex $(f_1,\dots,f_n)\in\CN\CC$ is \emph{simple}
with respect to $\Delta$ if the product $f_1\dots f_n$ is simple.
The set of simple simplices in $\CN\CC$ forms a simplicial set
which we denote by $\CN\CS$, the \emph{Garside nerve
of $\CC$ \index{nerve!of a Garside structure}
with respect to $\CS$.}
\end{defi}

\begin{theo}[``gar'' resolution, quotient version]
\label{gar1}
Let $\CG$ be a groupoid, let $\CS$ be a Garside structure
on $\CG$. The Garside nerve $\CN\CS$ is a simplicial $K(\CG,1)$.
\end{theo}

\begin{defi}[``gar'' flag complex]
\index{$\gar(\CG,\CS,o)$ (``gar'' flag complex)}
\label{gar2}
Let $\CG$ be a groupoid, let $\CS$ be a Garside structure
on $\CG$, let $o$ be an object of $\CG$. We denote by
$$\gar(\CG,\CS,o)$$ the simplicial complex with underlying space $\obj(o\downarrow \CG)$ (\ie, the set of morphisms in $\CG$ with source $o$), and such
that $\{g_0,\dots,g_k\}$ is a $k$-simplex if and only if, for all $i,j$, $g_i^{-1}g_j\in \CS$ or $g_j^{-1}g_i \in \CS$.

When $\CG$ is a group, we shorten the notation to $\gar(\CG,\CS)$ as there is only one
possible choice for $o$.
\end{defi}

In particular, a subset of $\CG$ spans a simplex if and only if
every subpair spans an edge. Such a simplicial complex is called a \emph{flag complex}. A flag complex is uniquely determined
by its $1$-skeleton.

A variant of Definition \ref{gar2} that more closely resembles Definition \ref{bar2} is as follows:

\begin{defi}[``gar'' simplicial set]
We denote by
$\gar'(\CG,\CS,o)$ be the subcomplex of $\bar(\CG,o)$ consisting of bar symbols $g_0[g_1|\dots|g_k]$ such that $$\forall i>0, g_i\in \CS \quad \text{and} \quad  g_1\dots g_k\in \CS.$$
\end{defi}

\begin{lemma}
The map $g_0[g_1|\dots|g_k] \mapsto \{g_0,g_0g_1,g_0g_1g_2,\dots,g_0\dots g_k\}$ induces
a bijection between non-degenerate simplices in $\gar'(\CG,\CS,o)$ and simplices
in $\gar(\CG,\CS,o)$.
\end{lemma}

\begin{proof}
Because of the natural length function on $\CS$, when $g,h\in\CG$ are distinct, we
cannot have both $g^{-1}h\in \CS$ and $h^{-1}g \in \CS$. Thus any simplex $A\subseteq\CG$ of
$\gar(\CG,\CS,o)$ admits a unique ordering $A = \{h_0,\dots,h_k\}$ such that
$i \leq j \Leftrightarrow h_i^{-1}h_j\in \CS$. The only bar symbol in the preimage of $A$
is $h_0[h_0^{-1}h_1|h_1^{-1}h_2| \dots | h_{k-1}^{-1}h_k]$.
\end{proof}

By comparing Definitions \ref{bar2} and \ref{gar2}, we see that Garside structures allow
a twofold gain:
\begin{itemize}
\item by contrast with $\bar(\CG,o)$, the complex $\gar(\CG, \CS, o)$ is finite-dimensional:
Garside structures on groupoids bound their cohomological dimension;
\item we can replace the abstract simplicial set by a very concrete simplicial complex
(that is actually a flag complex): in this paper, proving the $K(\pi,1)$ property
involves interpreting this simplicial complex as the nerve of an open covering.
\end{itemize}

\begin{theo}[``gar'' resolution, universal cover version]
\label{gar3}
Let $\CG$ be a groupoid, let $\CS$ be a Garside structure
on $\CG$, let $o$ be an object of $\CG$. The geometric realization of
$\gar(\CG, \CS, o)$ is contractible.
\end{theo}

We have slighly departed from results and phrasings that can be found in the litterature, but
Theorems \ref{gar1} and \ref{gar3} are easy categorical variants
of the main results in \cite{cmw}, which themselves are Garside group versions of results by Bestvina
about Artin groups, \cite{bestvina, brady}; these variants
can be proved using the same exact strategy.

\subsection{Garside sets}

\begin{defi}
\label{deficddelta}
\index{$\CDbullet$!Garside set}
Let $(\CC,\Delta,\phi)$ be a Garside structure with set of simple elements $\CS$.
The associated \emph{Garside set} is the collection
$$\CD_{\bullet}(\Delta) := (\CD_k(\Delta))_{k\in \BZ_{\geq 0}}$$
where
$$\CD_{k}(\Delta) := \{ (f_1,\dots,f_k)\in \CS^k | f_1\dots f_k = \Delta\},$$
together with, for all $k$, the following structure:
\begin{itemize}
\item the $k$ face maps $d_1,\dots,d_k:\CD_{k}(\Delta) \to \CD_{k-1}(\Delta)$,
such that for all $\sigma=(f_1,\dots,f_k)\in \CD_{k}(\Delta)$,
we have $d_1(\sigma)=(f_1f_2,f_3,\dots,f_k)$, $\dots$, $d_{k-1}(\sigma)=(f_1,\dots,f_{k-2},f_{k-1}f_k)$ and $d_k(\sigma)=(f_2,\dots,f_kf_1^{\phi})$
\item the $k$ degeneracy maps
$s_1,\dots,s_k:\CD_{k-1}(\Delta) \to \CD_{k}(\Delta)$ obtained by inserting
identities at the $k$-possible locations in $(f_1,\dots,f_{k-1})$,
\item the ``screwdriver'' map $\rho:\CD_{k}(\Delta) \to \CD_{k}(\Delta), (f_1,\dots,f_k)\mapsto
(f_2,\dots,f_k,f_1^{\phi})$
\end{itemize}
When there is no ambiguity on $\Delta$, we write $\CD_{\bullet}$ instead of
$\CD_{\bullet}(\Delta)$.
\end{defi}

Note the analogy with $\CD_{\bullet}(c)$ (Definition \ref{deficdc}) which,
by Theorem \ref{bulletiso}, happens to be isomorphic to the Garside set of the dual braid
monoid.

Together, the face and degeneracy maps form a (degree-shifted) simplicial structure,
which is isomorphism to that on $\CN\CS$ via the following 
trivial lemma:

\begin{lemma}[Kreweras map]
\label{transferlemma}
For all $k\geq 1$, the map
\begin{eqnarray*}
(\CN\CS)_{k-1} & \longrightarrow & \CD_k(\Delta)\\
(f_1,\dots,f_{k-1}) & \longmapsto & (f_1,\dots,f_{k-1}, (f_1\dots f_{k-1})^{-1}
\Delta)
\end{eqnarray*}
is bijective.
\end{lemma}

\begin{remark}[Helicoidal structure]
\label{helicoremark}
The are obvious compatibility axioms between the screwdriver map 
$\rho$ and the simplicial structure.
If $\phi$ were to act trivially on simple elements, these axioms would be that of 
a \emph{cyclic set}, in the sense of \cite{connes}.
B\"okstedt-Hsiang-Madsen provide variant axioms
for the case when $\rho$ has finite order (\emph{$|\rho|$-cyclic sets}, \cite{bhm}).
As it is interesting to consider Garside structures where $\rho$ has infinite order (see
for example \cite{free}), imposing conditions on the order of $\rho$ seems a bit artificial.
Relaxing the finite order condition yields a natural notion of \emph{helicoidal set}, which
I haven't found in the litterature (maybe I didn't search well enough).
The geometric realization of an helicoidal set is equipped with a natural $\mathbf{R}$-action, which
for cyclic sets and $d$-cyclic sets factors through a natural $S^1$-action.
This $S^1$-action, and its compatibility
with the scalar action on $V$ is the true explanation for the miracles
of Section \ref{section11}.
\end{remark}

\begin{remark}[Garside structure from Garside set]
\label{remarkrecover}
Note that the whole Garside
structure can be recovered from $\CD_{\bullet}$.
Consider the following category presentation:

\begin{itemize}
\item[($1$)] The set $\CD_{1}$ is, litteraly, the Garside family $\Delta$. It contains
one element per object in $\CC$, and serves as an abstract object set.
\item[($2$)] The set $\CD_{2}$ is in bijection with $\CS$.
We take $\CD_{2}$ as a formal set of generators.
\item[($3$)] The elements in $\CD_{3}$ are the defining relations of the presentation:
the triple $(f_1,f_2,f_3)$ expresses the relation $(f_1,f_2f_3)(f_2,f_3f_1^{\phi})=(f_1f_2,f_3)$
(indeed: via the Kreweras bijection between $\CD_{2}$ and $\CS$, this is
the relation $f_1\cdot f_2=f_1f_2$).
\end{itemize}

This presentation defines an abstract category $\CC(\CD_{\bullet})$ that is isomorphic to
$\CC$. In other words, the category $\CC$, together with $\Delta$ and $\CS$,
can be functorially retrieved from the helicoidal set $\CD_{\bullet}$. Note that
higher degree elements in $\CD_{\bullet}$ express further syzygies but, because
of Theorem \ref{gar1}, they only contain homotopically trivial stuff.
\end{remark}

\begin{exercise}
Write down the axioms for an helicoidal set.
\end{exercise}

\begin{question}
Is there a pleasant way to phrase axioms for
\emph{abstract Garside sets} (helicoidal sets $H$ such
that $\CC(H)$ is a Garside category)?
\end{question}

\subsection{Divided Garside structures}
The main construction in \cite{cyclic} is a kind of "barycentric subdivision" functor
for Garside categories. At the level of
cyclic sets, it coincides with an earlier construction by B\"okstedt-Hsiang-Madsen, 
\cite{bhm}. Thanks to our index-shifting Lemma \ref{transferlemma}, we can describe
this construction in remarkably simple terms (note how the Kreweras map
allows much simpler notations compared to those of \cite{bhm} and \cite{cyclic}):

\begin{defi}[divided Garside set]
\label{defidivided}
Let $(\CC,\Delta,\phi)$ be a Garside structure.
Let $m$ be a positive integer.
The \emph{$m$-divided Garside set} is the graded set
$$\sqrt[m]{\CD}_{\bullet}(\Delta) = (\sqrt[m]{\CD}_{k}(\Delta))_{k \in \BZ_{\geq 0}}$$
where 
$$\sqrt[m]{\CD}_{k}(\Delta):= \CD_{mk}(\Delta),$$
equipped with:
\begin{itemize}
\item faces $d'_1,\dots,d'_k:\sqrt[m]{\CD}_{k}(\Delta) \to \sqrt[m]{\CD}_{k-1}(\Delta)$
defined by $$d'_i = d_id_{i+k}\dots d_{i+(d'-1)k}$$ (composed from right to left),
\item degeneracy maps
$s'_1,\dots,s'_k:\sqrt[m]{\CD}_{k-1}(\Delta) \to \sqrt[m]{\CD}_{k}(\Delta)$
defined by $$s'_i = s_is_{s+k}\dots s_{i+(d'-1)k}$$ (composed from right to left),
\item the screwdriver map is the restriction of that of $\CD_{\bullet}$.
\end{itemize}
\end{defi}

\begin{theo}[after Section 9 in \cite{cyclic}]
\label{theodivided}
Let $(\CC,\Delta,\phi)$ be a Garside structure with set of simple elements $\CS$.
Let $m$ be a positive integer.
Then $\sqrt[m]{\CD}_{\bullet}(\Delta)$ is the Garside set of Garside structure
$(\CC_m,\Delta_m,\phi_m)$ such that $\CG_m$, the groupoid of fractions of $\CC_m$,
is equivalent as a category to $\CG$, the groupoid of fractions of $\CC$.
\end{theo}

Following Remark \ref{remarkrecover}, the Garside structure $(\CC_m,\Delta_m,\phi_m)$
is uniquely determined by its Garside set $\sqrt[m]{\CD}_{\bullet}(\Delta)$, and we can
write presentations by generators and relations for $\CC_m$ and $\CG_m$ as follows:

The underlying object is $\sqrt[m]{\CD}_{1}(\Delta)  = \CD_{m}(\Delta)$
and the set of simple elements $\CS_m$ is in bijection with $\sqrt[m]{\CD}_{2}(\Delta) =
\CD_{2m}(\Delta)$.
It is better to understand everything in terms of commutative diagrams.
An object $(f_1,\dots,f_m)$ is viewed as a commutative diagram

$$
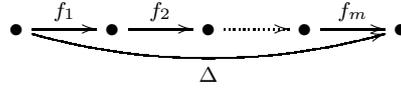
\xymatrix{
\bullet \ar[r]^{f_1} \ar@/_1em/[rrrr]_{\Delta}
& \bullet \ar[r]^{f_2} &  \bullet  \ar@{..>}[r]
& \bullet \ar[r]^{f_{m}} & \bullet}$$
where all arrows are in $\CS$, whereas a
simple morphism $(f_1,\dots,f_{2m})$ in $\CS_m\subseteq \CC_m$ is 
viewed as a commutative diagram

$$
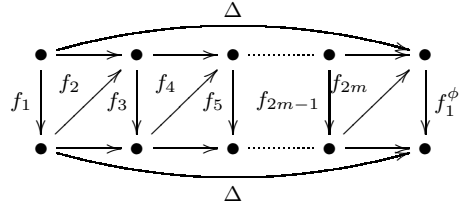
\xymatrix{
\bullet \ar[r] 
\ar@/^1em/[rrrr]^{\Delta}
\ar[d]_{f_1} & \bullet \ar[r] \ar[d]_{f_3}
& \bullet \ar[d]_{f_5}
\ar@{..}[r] & \bullet \ar[d]_{f_{2m-1}}
\ar[r] & \bullet  \ar[d]^{f_1^{\phi}} \\
\ar@/_1em/[rrrr]_{\Delta}
\bullet \ar[r] \ar[ur]^{f_2} & \bullet  \ar[ur]^{f_4}
 \ar[r] &
\bullet \ar@{..}[r] 
& \bullet \ar[r] \ar[ur]^{f_{2m}} & \bullet 
 }$$
where all arrows are in $S$; its source is the object
$(f_1f_2,f_3f_4,\dots,f_{2m-1}f_{2m})$
and its target is the object $(f_2f_3,f_4f_5,\dots,f_{2m}f_1^{\phi})$:

$$\xymatrix{
\bullet \ar[r]^{f_1f_2} \ar@/_1em/[rrrr]_{\Delta}
& \bullet \ar[r]^{f_3f_4} & \bullet  \ar@{..>}[r]
& \bullet \ar[r]^{f_{2m-1}f_{2m}} & \bullet
& \ar[r]^{(f_1,f_2,\dots,f_{2m})} & & 
\bullet \ar[r]^{f_2f_3} \ar@/_1em/[rrrr]_{\Delta}
& \bullet \ar[r]^{f_4f_5} & \bullet  \ar@{..>}[r]
& \bullet \ar[r]^{f_{2m}f_1^{\phi}} & \bullet}$$

These morphisms are subject to defining relations indexed by
$\sqrt[m]{\CD}_{3}(\Delta) = {\CD}_{3m}(\Delta)$.
The relation $(f_1,\dots,f_{3m})$ can be visualized as:

$$
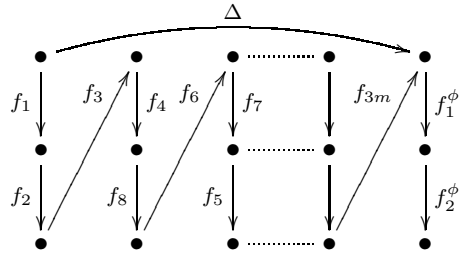
\xymatrix{
\ar@/^1em/[rrrr]^{\Delta}
\bullet  
\ar[d]_{f_1} & \bullet \ar[d]^{f_4}
& \bullet \ar[d]^{f_7}
\ar@{..}[r] & \bullet \ar[d]
 & \bullet  \ar[d]^{f_1^{\phi}} \\
\bullet
\ar[d]_{f_2} & \bullet  \ar[d]_{f_8}
& \bullet \ar[d]_{f_5}
\ar@{..}[r] & \bullet \ar[d]
 & \bullet  \ar[d]^{f_2^{\phi}} \\
\bullet \ar[uur]^(.75){f_3} & \bullet  \ar[uur]^(.75){f_6}
 &
\bullet \ar@{..}[r] 
& \bullet \ar[uur]^(.75){f_{3m}} & \bullet 
 }$$
and expresses the defining relation
$$(f_1,f_2f_3,f_4,f_5f_6,\dots)\cdot (f_2,f_3f_4,f_5,f_6f_7,\dots) = (f_1f_2,f_3,f_4f_5,f_6,\dots).$$

The automorphism $\phi_m$ acts on objects by
$$(f_1,f_2,\dots,f_m) \mapsto (f_2,\dots,f_m,f_1^{\phi})$$
and on simple morphisms by
$$(f_1,f_2,f_3,f_4,\dots,f_{2m-1},f_{2m})
\mapsto (f_3,f_4,\dots,f_{2m-1},f_{2m},f_1^{\phi},f_2^{\phi})$$

The Garside element $\Delta_m$ with source $(f_1,f_2,\dots,f_{m-1},f_m)$
has target $(f_2,f_3\dots,f_m,f_1^{\phi})$
and corresponds to the commutative diagram:
$$
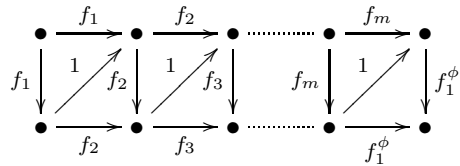
\xymatrix{
\bullet \ar[r]^{f_1} 
\ar[d]_{f_1} & \bullet \ar[r]^{f_2} \ar[d]_{f_2}
& \bullet \ar[d]_{f_3}
\ar@{..}[r] & \bullet \ar[d]_{f_m}
\ar[r]^{f_m} & \bullet  \ar[d]^{f_1^{\phi}} \\
\bullet \ar[r]_{f_2} \ar[ur]^{1} & \bullet  \ar[ur]^{1}
 \ar[r]_{f_3} &
\bullet \ar@{..}[r] 
& \bullet \ar[r]_{f_1^{\phi}} \ar[ur]^{1} & \bullet 
}$$
 
\begin{defi}
\label{deficollapse}
The \emph{collapse map} is
\begin{eqnarray*}
\kappa_m: \sqrt[m]{\CD}_{2}(\Delta) & \longrightarrow & \CD_{2}(\Delta)\\
(f_1,f_2,\dots,f_{2m})  & \longmapsto & (f_1,f_2\dots f_{2m})
\end{eqnarray*}
\end{defi}

By inspecting the defining relations of $\CC$ (Remark \ref{remarkrecover}) and $\CC_m$ (just above), one sees that $\kappa_m$ extends to a collapse functor $\CC_m\to \CC$.
In \cite[Section 9]{cyclic} is defined a less trivial functor
$\Theta_m: \CC\to \CC_m$.

\begin{theo}
\label{theocollapse}
The functors $\kappa_m:\CC_m\to \CC$ and $\Theta_m:\CC\to \CC_m$ are
such that $\kappa_m\circ \Theta_m = 1_{\CC}$.
They induce equivalences of categories $\kappa_m:\CG_m\to \CG$ and $\Theta_m:\CG\to \CG_m$.
\end{theo}

\begin{proof}
That $\kappa_m\circ \Theta_m = 1_{\CC}$ is obvious by construction. That $\Theta_m:\CG\to \CG_m$
is an equivalence of categories is \cite[Theorem 9.5]{cyclic}.
\end{proof}

\section*{Acknowledgements (2007)}
Work on this project started during
stimulating visits to KIAS (Seoul) and RIMS (Kyoto), in the summer of 2003,
when I realised that the dual braid monoid construction
could be generalized to well-generated
groups. I thank Sang Jin Lee and Kyoji Saito for their hospitality and
for their interests in discussing these topics.
After this initial progress, I remained stuck for many months, trying
to construct the open covering using convex geometry in $V^{\reg}$
(refining \cite[Section 4]{dualmonoid}).
Two observations were crucial to figuring out that working in the quotient space
was more appropriate. 
First, Kyoji Saito pointed out that the starting point for the
construction of the \emph{flat structure} (or \emph{Frobenius manifold
structure}, see \cite{saitoorbifold}) on real reflection orbifolds 
was precisely 
the duality between degrees and codegrees.
The intuition that the flat structure has something to do with the 
$K(\pi,1)$ property is explicitly mentioned as a motivation for
 \cite{saitoorbifold} (see also \cite{saitopolyhedra}).
The second useful discussion was with Fr\'ed\'eric Chapoton, who
pointed out the numerological coincidence, in the Coxeter case, between the
degree of the Lyashko-Looijenga covering and the number of
maximal chains in the lattice of
non-crossing partitions.

I thank Pierre Deligne, Eduard Looijenga,
Jean Michel, Vic Reiner and Vivien Ripoll
for comments, critics and suggestions.

Theorem \ref{theoroots} is an overdue answer to a question Michel Brou\'e
asked me in 1997, when I was a graduate student under his
supervision. Most of my work on braid groups was motivated by this
problem.
I am glad to have been able to solve it just in time for
his sixtieth birthday.

\section*{Postscript (2013)}

This article, based on results obtained between 2003 and 2006,
was first circulated as a preprint in October 2006.
Having spent three years focused on a single theorem,
I liked the idea of writing the proof as a single article.
That was naive, and I now understand that
a series of smaller papers would have been a much wiser approach to publication.

More than 5 years had passed between the initial submission, in April 2007, and
the beginning of the revision work. In the meantime, I had quit
academic life, and engaged in a new project that was both very demanding and
impossible to put on hold.
The revision work, from August 2012 to September 2013, took place during week-ends
and a few dedicated day offs.
A consequence is the inevitable Harlequin pattern of styles, not just because the proof
borrows from different areas of mathematics,
but because it was written over so many years.

The introduction, as well as Sections 1--6, 8--10 and 12--13, are from 2007
(with due corrections, naturally). Section 11 and Appendix B are brand new
from 2013. None of the referees could get through the old Section 11 and
my advice is to burn any remaining copy;
the new Section 11 follows the same argument, but in a cleaned-up,
clarified and hopefully intelligible manner.
It actually handles the well-generated case too ($d=1$), so
Sections 9--10 could be removed without logical harm.
Appendix A is expanded from an old
``Notation'' section. As for Section 7, it retains its original content, complemented
with a much clarifying factorization theorem (Theorem \ref{theolabel}) and an even
better one (Theorem \ref{nicetriv}). Least to say, the old sections now look incredibly
clumsy and immature.

I thank the four referees for their time and effort, especially Referee \#3,
for his enthusiastic comments, and even more so Referee \#4: his
careful and detailed 2012 report allowed the editorial process to resume,
effectively saving the paper. I thank David Gabai for his editorial patience and tenacity.

Two letters commenting my initial draft, sent by Pierre Deligne in 2006,
remained for many years the only sign that \emph{at least one person
seemed to believe my proof.} I thank him for his generosity.
I also thank Jean Michel, who demonstrated an indefatigable curiosity for this
paper, challenging me until he would understand it.

I would have renounced trying to get this work published,
if it wasn't for the comforting support
of my mathematical friends, especially Emmanuel Breuillard, Michel Brou\'e and Rapha\"el Rouquier.
 
\end{document}